\documentclass[a4paper]{article}
\usepackage{mypaper}

\usepackage[english]{babel}
\usepackage[utf8]{inputenc}
\usepackage[T1]{fontenc}

\usepackage{subcaption}
\usepackage{graphicx}
\usepackage{epstopdf}

\usepackage[mode=buildnew]{standalone}
\usepackage{tikz}
\usetikzlibrary{shapes,arrows,positioning,external,calc,math}
\usepackage{pgfplots}
\pgfplotsset{compat=1.16}

\usepackage{url}
\usepackage{verbatim}
\usepackage{multirow, booktabs, adjustbox} 
\usepackage{algorithm}
\usepackage{algorithmic}
\usepackage[shortlabels]{enumitem}

\newenvironment{enum_roman}{\begin{enumerate}[(i),nosep]}{\end{enumerate}}
\newenvironment{enum_alpha}{\begin{enumerate}[(a),nosep]}{\end{enumerate}}

\makeatletter
\@ifpackageloaded{ntheorem}{}{\usepackage{amsthm}}
\makeatother

\usepackage{amssymb}
\makeatletter
\@ifpackageloaded{amsmath}{}{\usepackage[nosumlimits]{amsmath}}
\makeatother
\usepackage{fixmath} 
\usepackage{xfrac} 
\usepackage[makeroom]{cancel} 
\usepackage{dsfont}
\usepackage{mathrsfs}
\usepackage{mathtools}
\usepackage{centernot}

\newcommand{\norm}[2][]{{#1\|}{#2}{#1\|}}
\newcommand{\inprod}[3][]{{#1\langle}{#2},{#3}{#1\rangle}}

\newcommand{\weakto}{\rightharpoonup}

\DeclareMathOperator{\proj}{\Pi}
\DeclareMathOperator{\Id}{Id}
\DeclareMathOperator{\resolv}{J}

\DeclareMathOperator{\dom}{dom}

\DeclareMathOperator{\fix}{fix}
\DeclareMathOperator{\zer}{zer}
\DeclareMathOperator{\ran}{ran}
\DeclareMathOperator{\gra}{gra}
\DeclareMathOperator{\rank}{rank}

\newcommand{\R}{\mathbb{R}}

\newcommand{\N}{\mathbb{N}}
\newcommand{\Hprim}{\mathcal{H}}

\newcommand{\defeq}{:=}

\usepackage[numbers,square,sort&compress]{natbib}
\usepackage{hyperref} 
\usepackage[capitalise]{cleveref} 

\crefformat{equation}{(#2#1#3)}
\crefrangeformat{equation}{(#3#1#4) to~(#5#2#6)}
\crefmultiformat{equation}{(#2#1#3)}{ and~(#2#1#3)}{, (#2#1#3)}{ and~(#2#1#3)}

\crefformat{enumi}{#2#1#3}
\crefrangeformat{enumi}{#3#1#4 to~#5#2#6}
\crefmultiformat{enumi}{#2#1#3}{ and~#2#1#3}{, #2#1#3}{ and~#2#1#3}

\newtheorem{thm}{Theorem}[section]
\newtheorem{lem}[thm]{Lemma}
\newtheorem{cor}[thm]{Corollary}
\newtheorem{prop}[thm]{Proposition}

\newtheorem{ass}[thm]{Assumption}
\newtheorem{defin}[thm]{Definition}

\makeatletter
\@ifpackageloaded{cleveref}{
	\crefname{thm}{Theorem}{Theorems}
	\crefname{lem}{Lemma}{Lemmas}
	\crefname{cor}{Corollary}{Corollaries}
	\crefname{prop}{Proposition}{Propositions}

	\crefname{ass}{Assumption}{Assumptions}
	\crefname{defin}{Definition}{Definitions}
	\crefname{rem}{Remark}{Remarks}
	\crefname{ex}{Example}{Examples}
}{}
\makeatother

\newcommand{\ops}{\mathcal{A}}

\newcommand{\pdop}{\Phi}
\newcommand{\pdskew}{\Gamma}
\newcommand{\pddiag}{\Delta}

\newenvironment{fmatrix}[1]{%
	\left[
		\setlength{\tabcolsep}{20pt}
		
		\begin{array}{#1}
			}{
		\end{array}
	\right]
}

\title{Frugal Splitting Operators: Representation, Minimal Lifting, and Convergence}
\author{
	Martin Morin
	\thanks{Department of Automatic Control, Lund University, (%
		\url{martin.morin@control.lth.se},
		\url{sebastian.banert@control.lth.se},
		\url{pontus.giselsson@control.lth.se}%
		)
	}
	\And
	Sebastian Banert\footnotemark[1]
	\And
	Pontus Giselsson\footnotemark[1]
}
\date{}

\begin{document}

\maketitle

\begin{abstract}
    We investigate frugal splitting operators for finite sum monotone inclusion problems. These operators utilize exactly one direct or resolvent evaluation of each operator of the sum, and the splitting operator's output is dictated by linear combinations of these evaluations' inputs and outputs.
To facilitate analysis, we introduce a novel representation of frugal splitting operators via a {\emph{generalized primal-dual resolvent}}.
The representation is characterized by an index and four matrices, and we provide conditions on these that ensure equivalence between the classes of frugal splitting operators and generalized primal-dual resolvents.
Our representation paves the way for new results regarding lifting numbers and the development of a unified convergence analysis for frugal splitting operator methods, contingent on the directly evaluated operators being cocoercive. 
The minimal lifting number is $n-1-f$ where $n$ is the number of monotone operators and $f$ is the number of direct evaluations in the splitting.
Notably, this lifting number is achievable only if the first and last operator evaluations are resolvent evaluations.
These results generalize the minimal lifting results by Ryu and Malitsky--Tam that consider frugal resolvent splittings. 
Building on our representation, we delineate a constructive method to design frugal splitting operators, exemplified in the design of a novel, convergent, and parallelizable frugal splitting operator with minimal lifting.
\end{abstract}
\section{Introduction}
A zero of a maximally monotone operator can be determined through a fixed point iteration of its resolvent \cite{rockafellarMonotoneOperatorsProximal1976}. However, when an operator is a finite sum of maximally monotone operators with individually readily computable resolvents, its overall resolvent may not be directly computable. Splitting methods address this by evaluating each resolvent individually and merging the results to form a convergent fixed point iteration.
In this work, we consider a comprehensive class of splitting operators (that give rise to a splitting method when iterated) for maximally monotone finite sum problems that, inspired by the terminology in \cite{ryuUniquenessDRSOperator2020,malitskyResolventSplittingSums2021}, we call \emph{frugal splitting operators}.
Informally, the class contains all operators whose fixed points encode the zeros of the sum of the monotone operators and can be computed with exactly one evaluation of each operator, either directly or via a resolvent.
Apart from the operator evaluations, only predetermined linear combinations of the input and operator evaluations are allowed.
The class encompasses the Douglas--Rachford \cite{lionsSplittingAlgorithmsSum1979} and forward-backward \cite{goldsteinConvexProgrammingHilbert1964,levitinConstrainedMinimizationMethods1966} operators along with many others, for instance \cite{ryuUniquenessDRSOperator2020,malitskyResolventSplittingSums2021,chambolleFirstOrderPrimalDualAlgorithm2011,davisThreeOperatorSplittingScheme2017,vuSplittingAlgorithmDual2013,condatPrimalDualSplitting2013,malitskyForwardBackwardSplittingMethod2020,alvarezInertialProximalMethod2001,moudafiConvergenceSplittingInertial2003,ryuFindingForwardDouglasRachfordForward2020,botInertialDouglasRachford2015,raguetGeneralizedForwardBackwardSplitting2013}.

We introduce the notion of a \emph{generalized primal-dual resolvent} and provide an equivalence between the class of frugal splitting operators and generalized primal-dual resolvents.
The latter is inspired by the works of \cite{latafatAsymmetricForwardBackward2017,giselssonNonlinearForwardBackwardSplitting2021,morinNonlinearForwardBackwardSplitting2022,buiWarpedProximalIterations2020} that propose other resolvent generalizations leading to powerful algorithm modeling tools, covering a wide range of algorithms. However, our resolvent generalization stands out in that it is the first to provably fully cover the class of frugal splitting operators.
This novel representation is a key difference to the related works of \cite{ryuUniquenessDRSOperator2020,malitskyResolventSplittingSums2021} that examine frugal resolvent splitting operators, i.e., frugal splitting operators that only use resolvents and no direct evaluations.
These works focus on the \emph{lifting number} of frugal splitting operators. While we also provide minimal lifting results, our representation allows to easily design new splitting operators and analyze the convergence of their fixed point iterations in a general setting.
The representation also allows us to relax one assumption of \cite{ryuUniquenessDRSOperator2020,malitskyResolventSplittingSums2021} regarding the existence of a \emph{solution map} since the existence of such a solution map for any frugal splitting operator is evident directly from our new representation.
The representation also directly reveals which individual resolvent/direct evaluations can be performed in parallel.

The general idea behind lifting is to trade computational complexity for storage complexity by creating an easier to solve problem in some higher dimensional space.
In the context of frugal splitting operators, this means that, while the monotone inclusion problem lies in some real Hilbert space $\Hprim$, the splitting operator maps to and from $\Hprim^d$ for some non-zero natural number $d$ which we call the lifting number.
For example, the Douglas--Rachford and forward-backward splitting operators have lifting number one while the splitting operator of the primal-dual method of Chambolle and Pock \cite{chambolleFirstOrderPrimalDualAlgorithm2011} has lifting number two%
\ifdefined\shortversion
	{}.
\else
	\footnote{
		The Chambolle--Pock method considers monotone inclusion problems that allow for compositions with linear operators.
		We will implicitly assume that the linear operators are the identity operator when discussing Chambolle--Pock or other similar methods for monotone inclusions with compositions.
	}.
\fi

In the three operator case, Ryu \cite{ryuUniquenessDRSOperator2020} showed that the smallest possible lifting number for a frugal resolvent splitting operator is two.
Malitsky and Tam \cite{malitskyResolventSplittingSums2021} later expanded this result to sums of $n$ maximally monotone operators and found that the minimal lifting is $n-1$.
In this paper, we show that this lifting number can be reduced to $n-1-f$ where $f$ is the number of operators that are evaluated directly and not via resolvents.
We also show that the minimal lifting number is dependent on the order of the direct and resolvent evaluations.
In particular, if the first or last operator is evaluated directly, the minimal lifting is $n-f$ instead of $n-1-f$.
An example of this can be found in the three operator splitting of Davis and Yin \cite{davisThreeOperatorSplittingScheme2017}, which places its only direct evaluation second and hence can achieve a lifting number of one.
This is to be compared to the primal-dual method of V\~u and Condat \cite{vuSplittingAlgorithmDual2013,condatPrimalDualSplitting2013}, which performs a direct evaluation first and hence requires a higher lifting number of two in the three operator case.

We provide sufficient conditions for weak sequence convergence of a fixed point iteration of a frugal splitting operator in which all operators evaluated via direct evaluation are cocoercive. The convergence analysis entails establishing Fej\'er monotonicity of the generated sequence w.r.t.\ to the fixed points of the splitting operator.
There are a number of different general or unified approaches for analyzing algorithm classes \cite{lessardAnalysisDesignOptimization2016,latafatAsymmetricForwardBackward2017,giselssonNonlinearForwardBackwardSplitting2021,morinNonlinearForwardBackwardSplitting2022,buiWarpedProximalIterations2020,droriPerformanceFirstorderMethods2014,taylorExactWorstCasePerformance2017,taylorSmoothStronglyConvex2017,ryuOperatorSplittingPerformance2020}.
Many of these approaches can be applied to a frugal resolvent splitting algorithm.
Our analysis has the advantage of being performed directly on the generalized primal-dual resolvent representation, leading to conditions for convergence and the conditions for a generalized primal-dual resolvent to be a frugal resolvent splitting on the same set of matrices. This simplifies the design of new frugal resolvent splittings.
As an example, we construct a new convergent and parallelizable frugal splitting operator and briefly discuss its relation to the existing splitting operators with minimal lifting in \cite{ryuUniquenessDRSOperator2020,malitskyResolventSplittingSums2021,aragon-artachoDistributedForwardBackwardMethods2022,condatProximalSplittingAlgorithms2021,campoyProductSpaceReformulation2022}.
Note that the works in \cite{ryuUniquenessDRSOperator2020,malitskyResolventSplittingSums2021} show that their algorithms have minimial lifting, which is not the case for \cite{aragon-artachoDistributedForwardBackwardMethods2022,condatProximalSplittingAlgorithms2021,campoyProductSpaceReformulation2022}. However, these can retroactively be shown to have minimal lifting using the results of \cite{malitskyResolventSplittingSums2021} and this paper.


Some of our proofs require the underlying real Hilbert space $\Hprim$ of the monotone inclusion problem to have dimension greater than one, which we will make a blanket assumption.
This is not a significant restriction, partly because the zero-dimensional and one-dimensional cases are of limited practical interest and partly because many of our results still hold in these cases.
For instance, although the proof for the necessary conditions of our representation theorem no longer holds when $\dim \Hprim \leq 1$, the proof for the sufficient conditions still holds.
Hence, if we find a representation that yields a frugal splitting operator when $\dim \Hprim \geq 2$, then it also yields a frugal splitting operator when $\dim \Hprim \leq 1$.
All frugal splitting operators presented in this paper are therefore also applicable to the case $\dim \Hprim \leq 1$.

\ifdefined\shortversion
\else
	\subsection{Outline}
	In \cref{sec:prel}, we introduce some preliminary notation and results together with the main monotone inclusion problem.
	We define and discuss our definition of a frugal splitting operator in \cref{sec:fosplit}.
	\Cref{sec:gen-pd-resolv} contains the definition of a generalized primal-dual resolvent which we will use to represent frugal splitting operators.
	The lemmas that prove our representation results can be found in \cref{sec:rep-lemmas} and they are summarized in our main convergence theorem in \cref{sec:representation} which also contains some general remarks on the representation.
	For instance, the relationship between fixed points of the splitting and solutions to the monotone inclusion is proven and how the parallelizable resolvent/direct evaluations can be identified is demonstrated.
	We also show how a representation can be derived via an example.
	The minimal lifting results in terms of a rank bound on a structured matrix can be found in \cref{sec:minlift}.
	This result covers the setting of Malitsky and Tam \cite{malitskyResolventSplittingSums2021} and shows that the minimal lifting number depends on whether the first or last operator evaluation is a forward or backward evaluation.
	Convergence under cocoercive forward evaluations is proven in \cref{sec:conv}.
	In \cref{sec:conv-examples} we apply the convergence theorem and reestablish the convergence criterion of the three operator splitting of Davis and Yin \cite{davisThreeOperatorSplittingScheme2017} as well as provide conditions for the convergence of a fixed point iteration of forward-backward splitting with Nesterov-like momentum.
	The last part of the paper, \cref{sec:new-method}, contains the construction of a new frugal splitting operator with minimal lifting and parallelizable forward/backward evaluations.
	The paper ends with a short conclusion in \cref{sec:conclusion}.
	In the accompanying supplement many more examples of representations of frugal splitting operators and application of our convergence theorem can be found.
\fi

\section{Preliminaries}\label{sec:prel}
Let $\N = \{0,1,\dots\}$ be the set of natural numbers and $\N_+ = \{1,2,\dots\}$ be the set of non-zero natural numbers.
The cardinality of a set $A$ is denoted by $|A|$.
A subset $B$ of $A$ is denoted by $B \subseteq A$, while a strict subset is denoted by $B \subset A$.
Let $\R$ be the set of real numbers.
We refer to the range and kernel of a matrix $A \in \R^{n\times m}$ as $\ran A$ and $\ker A$ respectively.
These are linear subspaces of $\R^n$ and $\R^m$ respectively and their orthogonal complements with respect to the standard Euclidean inner product are denoted by $(\ran A)^\perp$ and $(\ker A)^\perp$ respectively, which also are linear subspaces.
With $\Hprim$ being a real Hilbert space, the set of all subsets of $\Hprim$ is denoted $2^\Hprim$.
The inner product and norm on $\Hprim$ are denoted by $\inprod{\cdot}{\cdot}$ and $\norm{\cdot}$ respectively.
Given a bounded, self-adjoint, and positive linear operator $H\colon\Hprim\to\Hprim$, we define $\inprod{\cdot}{\cdot}_H = \inprod{H(\cdot)}{\cdot}$ and $\norm{\cdot}_H^2 = \inprod{H(\cdot)}{\cdot}$.
If $H$ is self-adjoint and strongly positive then $\inprod{\cdot}{\cdot}_H$ is an inner product and $\norm{\cdot}_H$ is a norm.
\ifdefined\shortversion
	For notation and properties regarding set-valued operators $A\colon \Hprim \to 2^\Hprim$ and their resolvents $\resolv_A = (\Id + A)^{-1}$, we refer to \cite{bauschkeConvexAnalysisMonotone2017}.
\fi

Let $U$, $V$ and $W$ be sets and $A\colon U\times W \to V$ an operator.
With the notation $A_wu = v$ we mean $A(u,w) = v$ and we define $A_w = A(\cdot,w)$.
Since $A_w$ is an operator from $U$ to $V$, instead of writing $A\colon U\times W \to V$ we will say that $A_{(\cdot)} \colon U \to V$ is \emph{parameterized} by $W$.

\ifdefined\shortversion
\else
	Let $A\colon \Hprim \to 2^\Hprim$ be a \emph{set-valued operator} on $\Hprim$, i.e., $A$ is an operator that maps any point in $\Hprim$ to a subset of $\Hprim$.
	The \emph{graph} of $A$ is $\gra A = \{(x,u) \in \Hprim\times\Hprim \mid u \in Ax\}$.
	The \emph{range} of $A$ is $\ran A = \{u\in\Hprim \mid \exists x\in \Hprim \text{ s.t. } (x,u) \in \gra A\}$.
	The \emph{domain} of $A$ is $\dom A = \{x \in\Hprim \mid Ax \neq \emptyset\}$.
	The operator $A$ is said to have \emph{full domain} if $\dom A = \Hprim$.
	If $Ax$ is a singleton for all $x\in\Hprim$ then it is said to be \emph{single-valued}.
	A single-valued operator has full domain and we will make no distinction between single-valued operators $A\colon \Hprim\to 2^\Hprim$ and mappings $A\colon \Hprim \to \Hprim$.

	An operator $A \colon \Hprim \to 2^\Hprim$ is \emph{monotone} if
	\[
		\inprod{u-v}{x-y} \geq 0
	\]
	for all $(x,u) \in \gra A$ and all $(y,v) \in \gra A$.
	A monotone operator is \emph{maximal} if its graph is not contained in the graph of any other monotone operator.
	The \emph{resolvent} of a maximally monotone operator $A$ is $\resolv_A = (\Id + A)^{-1}$.
	An operator $A\colon \Hprim \to 2^\Hprim$ is \emph{$\mu$-strongly monotone} where $\mu > 0$ if
	\[
		\inprod{u-v}{x-y} \geq \mu\norm{x-y}^2
	\]
	for all $(x,u) \in \gra A$ and all $(y,v) \in \gra A$.
	An operator $A \colon \Hprim \to 2^\Hprim$ is \emph{$\beta$-cocoercive} if it is single-valued and
	\[
		\inprod{Ax-Ay}{x-y} \geq \beta\norm{Ax-Ay}^2
	\]
	for all $x,y\in\Hprim$.
	The inverse of a $\beta$-cocoercive operator is $\beta$-strongly monotone.
\fi

\begin{lem}\label{lem:eq-range-right-fact}
	Let $A\in\R^{n\times m}$ and $B\in\R^{n\times d}$.
	If $\ran A \subseteq \ran B$, there exists a unique $S\in\R^{d\times m}$ such that $A = BS$ and $\ran S \subseteq (\ker B)^\perp$.
	If $\ran A = \ran B$, such an $S$ satisfies $\ran S = (\ker B)^\perp$.

	\begin{proof}
		Since $\ran A \subseteq \ran B$, the columns of $A$ lie in the span of the columns of $B$.
		This means that there exists a matrix $S'\in\R^{d\times m}$ such that $A = BS'$.
		Let $\proj_\perp$ be the orthogonal projection onto $(\ker B)^\perp$ in the standard Euclidean inner product and define $S = \proj_\perp S'$.
		It is clear that $\ran S \subseteq (\ker B)^\perp$ and, since $B = B\proj_\perp$, it also holds that $A = BS' = B\proj_\perp S' = BS$.

		To show uniqueness, let $S$ and $S'$ be such that $A = BS = BS'$ and $\ran S \subseteq (\ker B)^\perp$ and $\ran S' \subseteq (\ker B)^\perp$ and let $x\in\R^m$ be such that $Sx \neq S'x$.
		Since $Sx-S'x \neq 0$ and $Sx-S'x \in(\ker B)^\perp$ it must hold that $Sx-S'x \notin \ker B$ and we have
		\[
			0 \neq B(Sx-S'x) = BSx - BS'x = Ax - Ax = 0
		\]
		which is a contradiction.

		To show the last statement, assume $\ran A = \ran B$ and let $S$ be the unique matrix that satisfies $A = BS$ and $\ran S \subseteq (\ker B)^\perp$.
		Assume $\ran S \subset (\ker B)^\perp$, then $\rank S < \dim (\ker B)^\perp = \rank B$ and $\rank A \leq \min(\rank B,\rank S) < \rank B$.
		However, this contradicts $\ran A = \ran B$ and hence $\ran S = (\ker B)^\perp$.
	\end{proof}
\end{lem}

\begin{lem}\label{lem:eq-ker-left-fact}
	Let $A\in\R^{m\times n}$ and $B\in\R^{d\times n}$.
	If $\ker A \supseteq \ker B$, there exists a unique $S\in\R^{m\times d}$ such that $A = SB$ and $\ker S \supseteq (\ran B)^\perp$.
	If $\ker A = \ker B$ such an $S$ satisfies $\ker S = (\ran B)^\perp$.

	\begin{proof}
	    This is the dual to \cref{lem:eq-range-right-fact}.
		If $\ker A \supseteq \ker B$ then $\ran A^T \subseteq \ran B^T$ and \cref{lem:eq-range-right-fact} then implies the existence of a unique $S^T\in\R^{d\times m}$ with $\ran S^T \subseteq (\ker B^T)^\perp$ such that $A^T = B^TS^T$ or equivalently $A = SB$.
		Since $\ran S^T = (\ker S)^\perp$ and $(\ker B^T)^\perp = \ran B$ we have $(\ker S)^\perp \subseteq \ran B$ or equivalently $\ker S \supseteq (\ran B)^\perp$.
		If $\ker A = \ker B$ then $\ran A^T = \ran B^T$ and \cref{lem:eq-range-right-fact} then yields $\ker S = (\ran B)^\perp$.
	\end{proof}
\end{lem}

\subsection{Problem and Notation}\label{sec:problem}
For the remainder of this paper, let $\Hprim$ be a real Hilbert space with $\dim \Hprim \geq 2$.
We consider the problem of finding a zero of a finite sum of operators,
\begin{equation}\label{eq:prob}
	\text{find $x^\star \in \Hprim$ such that}\; 0 \in \sum_{i=1}^n A_ix^\star,
\end{equation}
where $A_i\colon\Hprim\to 2^\Hprim$ is maximally monotone for all $i\in\{1,\dots,n\}$ and $A_i$ is single-valued for $i\in F$ for some $F \subseteq \{1,\dots,n\}$.
Instead of solving \cref{eq:prob} directly, we will work with a family of primal-dual problems.
Let $p \in \{1,\dots,n\}$. A primal-dual problem associated with \cref{eq:prob} is then
\begin{equation}\label{eq:primdual-base}
	\text{find $(y_1^\star,\dots,y_n^\star)\in\Hprim^n$ such that}
	\begin{cases}
		0 \in A_{i}^{-1}y_i^\star - y_p^\star  \text{ for all } i \in  \{1,\dots,n\} \setminus \{p\}, \\
		0 \in A_{p} y_p^\star + \sum_{j \in \{1,\dots,n\} \setminus \{p\}} y_j^\star.
	\end{cases}
\end{equation}
We call $p$ the primal index since the corresponding variable, $y_p^\star$, solves the primal problem \cref{eq:prob}.
The equivalence between \cref{eq:primdual-base} and \cref{eq:prob} is straightforward to show and holds in the sense that if $y_p^\star \in \Hprim$ is a solution to \cref{eq:prob} then there exists $y_i^\star \in \Hprim$ for all $i \in \{1,\dots,n\}\setminus \{p\}$ such that $(y_1^\star,\dots,y_{n}^\star)$ solves \cref{eq:primdual-base}.
Conversely, if $(y_1^\star,\dots,y_{n}^\star) \in \Hprim^n$ is a solution to \cref{eq:primdual-base} then $y_p^\star$ solves \cref{eq:prob}.

The aim of this paper is to examine a class of iterative methods for solving problems \cref{eq:prob,eq:primdual-base}.
We will give the exact definition of the considered class of solution methods in \cref{sec:fosplit} and in the remainder of this section we introduce some notation in order to simplify its description.
\begin{defin}[Operator Tuples]
	With $F \subseteq \{1,\dots,n\}$, $\ops_n^F$ is the set of operator tuples where $A = (A_1,\dots,A_n) \in \ops_n^F$ if and only if
	\begin{enum_roman}
	\item $A_i \colon \Hprim \to 2^{\Hprim}$ is maximally monotone for all $i\in\{1,\dots,n\}$,
	\item $A_i$ is single-valued for all $i\in F$.
	\end{enum_roman}
	We further define $\ops_n = \ops_n^\emptyset$ and note that $\ops_n^F \subseteq \ops_n$ for all $F \subseteq \{1,\dots,n\}$.
\end{defin}
This allows us to associate a tuple $A = (A_1,\dots,A_n) \in\ops_n^F$ with each monotone inclusion problem of the form \cref{eq:prob}, or \cref{eq:primdual-base}, and vice versa.
To simplify the notation further we will also make extensive use of the following (abuse of) notation.
\begin{defin}[Matrix as Operator]
	Let $B \in \R^{n \times m}$ and $z = (z_1,\dots,z_m) \in \Hprim^m$.
	We define $Bz$ as%
	\ifdefined\shortversion
		{}
	\else
		\footnote{
			This operator could more accurately be represented by $B\otimes \Id$ where $\otimes$ is the tensor product.
			However, this notation will quickly become tedious.
		}
	\fi
	\begin{align*}
		Bz =
		\begin{bmatrix}
			B_{11}z_1 + \dots + B_{1m}z_m \\
			\vdots \\
			B_{n1}z_1 + \dots + B_{nm}z_m \\
		\end{bmatrix}.
	\end{align*}
\end{defin}
With this notation in place, we can identify a \emph{primal-dual operator} $\pdop_{A,p}\colon\Hprim^n \to 2^{\Hprim^n}$ with $A\in\ops_n^F$ and $p\in\{1,\dots,n\}$ that allows us to write a primal-dual problem \cref{eq:primdual-base} associated with $A$ as
\begin{equation}\label{eq:primdual}
	\text{find $y^\star\in \Hprim^n$ such that}\; 0 \in \pdop_{A,p} y^\star = \pddiag_{A,p}y^\star + \pdskew_p y^\star.
\end{equation}
The operator $\pddiag_{A,p}\colon \Hprim^n \to 2^{\Hprim^n}$ and the matrix $\pdskew_p \in \R^{n\times n}$ are defined as
\begin{equation}\label{eq:Delta-and-hats}
	\pddiag_{A,p}(y_1,\dots,y_n) = \widehat{A}_1y_1\times \dots \times \widehat{A}_ny_n
	\quad\text{and}\quad
	\pdskew_p = R_p - R_p^T,
\end{equation}
where $\widehat{A}_p = A_p$, $\widehat{A}_i = A_i^{-1}$ for all $i\in\{1,\dots,n\}\setminus \{p\}$, and $R_p \in \R^{n\times n}$ is the  matrix with ones on the $p$th row and zeros in all other positions, rendering $\Gamma_p$ skew-symmetric.
\ifdefined\shortversion
\else
	For illustration, in the case when $p=n$, the operator $\pddiag_{A,n}$ and the matrix $\pdskew_n$ have the following structures
	\[
		\pddiag_{A,n}(y_1,\dots,y_{n}) =
		\begin{bmatrix}
			A_1^{-1}y_1 \\
			\vdots \\
			A_{n-1}^{-1}y_{n-1} \\
			A_{n}y_{n} \\
		\end{bmatrix}
		\quad\text{and}\quad
		\pdskew_n =
		\begin{bmatrix}
			0 & \cdots & 0 & -1 \\
			\vdots & \ddots & \vdots & \vdots \\
			0 & \cdots & 0 & -1 \\
			1 & \cdots & 1 & 0
		\end{bmatrix}.
	\]
\fi

All these operators are maximally monotone; $\pdskew_p$ as an operator on $\Hprim^n \to \Hprim^n$ is bounded linear and skew-adjoint, $\pddiag_{A,p}$ is separable w.r.t.\ the components with each component operator being maximally monotone, and $\pdop_{A,p}$ is the sum of two maximally monotone operators with one having full domain.
The primal-dual operator $\pdop_{A,p}$ will feature extensively in the rest of the paper.

\section{Frugal Splitting Operators}\label{sec:fosplit}
A common way of solving \cref{eq:prob} associated with some $A\in \ops_n^F$ is with a fixed point iteration.
These are methods that, given some initial iterate $z_0\in\Hprim^d$ and operator $T_A\colon\Hprim^d \to \Hprim^d$, iteratively perform
\begin{equation}\label{eq:fp-iteration}
	z_{k+1} = T_Az_k
\end{equation}
for $k\in\N$.
The operator $T_{A}$ is such that the sequence $\{z_k\}_{k\in\N}$ converges to a fixed point from which a solution to \cref{eq:prob} can be recovered.
In this paper, this means that $T_{(\cdot)}$ is a \emph{frugal splitting operator} and the main focus will be on examining the representation and properties of such an operator.
\begin{defin}[Frugal Splitting Operator]\label{def:fo-split}%
	Let $d\in\N_+$ and $T_{(\cdot)} \colon \Hprim^d \to \Hprim^d$ be parameterized by $\ops_n^F$.
	We say that $T_{(\cdot)}$ is a \emph{frugal splitting operator} over $\ops_n^F$ if, for all $i\in\{1,\dots,n\}$, there exists $\tau_{i,(\cdot)} \colon \Hprim^{d_i} \to \Hprim^{d_{i+1}}$ parameterized by $\ops_n^F$ such that, for all $A = (A_1,\dots,A_n) \in\ops_n^F$,
	\begin{enum_roman}
	\item\label{def:fo-split:fix} $\fix T_A \neq \emptyset \iff \zer \sum_{i=1}^n A_i \neq \emptyset$,
	\item\label{def:fo-split:comp} $T_A = \tau_{n,A}\circ \tau_{n-1,A} \circ \cdots \circ \tau_{1,A}$.
	\end{enum_roman}
	Furthermore, for each $i\in\{1,\dots,n\}$ there exist $B_i \in \real^{d_{i+1} \times d_i}$, $C_i \in \real^{d_{i+1} \times 1}$, $D_i \in \real^{1 \times d_i}$ and $\gamma_i > 0$
	such that
	\[
		\tau_{i,A} z_i
		=
			\begin{cases}
				B_i z_i + C_i \resolv_{\gamma_i A_i}(D_i z_i) & \text{if } i \notin F \\
				B_i z_i + C_i A_i(D_i z_i)& \text{if } i \in F \\
			\end{cases}
	\]
	for all $A = (A_1,\dots,A_n) \in \ops_n^F$ and all $z_i \in\Hprim^{d_i}$.
	Note that $\tau_{i,A}$ only uses $A_i$ and $d_1 = d_{n+1} = d$.

	A frugal splitting operator over $\ops_n$ is called a \emph{resolvent splitting operator}.
\end{defin}
Following the terminology introduced in \cite{ryuUniquenessDRSOperator2020}, we call the class frugal since the computational requirement in \cref{def:fo-split}\cref{def:fo-split:comp} specifies exactly one evaluation of each operator in $A$, either directly or via a resolvent.
We will also refer to a direct evaluation as a forward evaluation while a resolvent evaluation will be referred to as a backward step or backward evaluation.
Apart from forward and backward evaluations, only predetermined vector additions and scalar multiplications of the inputs and results of the operator evaluations are allowed.
\cref{def:fo-split}\cref{def:fo-split:comp} also specifies evaluation order in that the operators $A_i$ in the tuple $A$ must be used in the order they appear in $A$ when evaluating $T_{A}$.
Since the operators of the original monotone inclusion problem can be arbitrarily rearranged, this entails no loss of generality.

\Cref{def:fo-split} is functionally the same as the definitions of Ryu and Malitsky--Tam \cite{ryuUniquenessDRSOperator2020,malitskyResolventSplittingSums2021} but with the addition of forward evaluations being allowed.
However, both Ryu \cite{ryuUniquenessDRSOperator2020} and Malitsky--Tam \cite{malitskyResolventSplittingSums2021} assume the existence of a computationally tractable \emph{solution map} $S_{A}\colon \Hprim^d \to \Hprim$ that maps fixed points of $T_A$ to solutions of \cref{eq:prob}.
One of our results, \cref{prop:sol-map}, shows that this is not needed since \cref{def:fo-split:fix,def:fo-split:comp} together imply the existence of such a solution map.
\ifdefined\shortversion
	Note that \cref{def:fo-split} only considers the solution encoding and the computational requirements of a splitting; it does not make any assumptions regarding the convergence of its fixed point iteration as in \cite{ryuUniquenessDRSOperator2020}.
\fi

The class of frugal splitting operators in \Cref{def:fo-split} covers many splitting operators such as the operators used in the forward-backward method \cite{goldsteinConvexProgrammingHilbert1964,levitinConstrainedMinimizationMethods1966} and Douglas--Rachford method \cite{lionsSplittingAlgorithmsSum1979}, but also operators of more recent methods such as the three operator splitting method of Davis and Yin \cite{davisThreeOperatorSplittingScheme2017} and various primal-dual methods such as the Chambolle and Pock method \cite{chambolleFirstOrderPrimalDualAlgorithm2011}.
However, it does not allow for multiple evaluations of an operator and therefore does not cover the forward-backward-forward method of Tseng \cite{tsengModifiedForwardBackwardSplitting2000} although it covers the forward-reflected-backward method of Malitsky and Tam \cite{malitskyForwardBackwardSplittingMethod2020}.

\ifdefined\shortversion
\else
	Note that \cref{def:fo-split} only considers the solution encoding and the computational requirements of a splitting and does not make any assumptions regarding the convergence of its fixed point iteration.
	For instance, consider a fixed point iteration of the forward-backward splitting operator $T_{(A_1,A_2)} = \resolv_{\gamma A_2}\circ(\Id - \gamma A_1)$ where $\gamma > 0$.
	This is a frugal splitting operator over $\ops_{2}^{\{1\}}$ but it is well known that its fixed point iteration can fail to converge without further assumptions on $A_1$ and/or $A_2$.
	The standard assumption is cocoercivity of $A_1$ but even then the fixed point iteration can fail to converge if the step-size $\gamma$ is too large.
	Further examples exist with both \cite{ryuUniquenessDRSOperator2020,malitskyResolventSplittingSums2021} providing resolvent splitting operators whose fixed point iterations fail to converge in general.
	Since we consider the question of convergence as separate to the definition of a frugal splitting operator, we will treat convergence separately in \cref{sec:conv} where we provide sufficient convergence conditions that are applicable to any fixed point iteration of a frugal splitting operator.
\fi

\section{Generalized Primal-Dual Resolvents}\label{sec:gen-pd-resolv}
Although \cref{def:fo-split} fully defines all operators we aim to consider, we do not find it conducive for analysis.
In this section, we will therefore develop an equivalent representation in the form of what we call a generalized primal-dual resolvent.
\begin{defin}[Generalized Primal-Dual Resolvent]\label{def:gen-pd-resolv}%
	Let $d\in\N_+$ and $T_{(\cdot)} \colon \allowbreak \Hprim^d \to \Hprim^d$ be parameterized by $\ops_n^F$.
	We say that $T_{(\cdot)}$ is a \emph{generalized primal-dual resolvent} if there exist
	$p \in \{1,\dots,n\}$, $M\in\R^{n\times n}$, $N\in\R^{n\times d}$, $U\in\R^{d\times d}$ and $V\in\R^{d\times n}$
	such that, for all $A\in\ops_n^F$ and all $z\in\Hprim^d$,
	$(M+\pdop_{A,p})^{-1}\circ N$ is single-valued at $z$ and
	\begin{equation}\label{eq:asym-res-op}
		\begin{aligned}
			y &= (M + \pdop_{A,p})^{-1}Nz, \\
			T_Az &= z - Uz + Vy,
		\end{aligned}
	\end{equation}
	where $y\in\Hprim^n$ and $\pdop_{A,p}$ is defined in \cref{eq:primdual}.

	We call the tuple $(p,M,N,U,V)$ a \emph{representation} of $T_{(\cdot)}$.
\end{defin}

When $M = N = \gamma^{-1}I$ and $U = V = \theta I$, where $\theta \in (0,2]$, $\gamma > 0$ and $I\in\R^{n\times n}$ is the identity matrix,
the generalized primal-dual resolvent becomes a standard relaxed resolvent operator of the primal-dual operator $\pdop_{A,p}$,
\[
	(1- \theta)\Id + \theta(\gamma^{-1}\Id + \pdop_{A,p})^{-1}\circ \gamma^{-1}\Id
	=
	(1- \theta)\Id + \theta\resolv_{\gamma \pdop_{A,p}}.
\]
While the standard relaxed resolvent operator is single-valued, the operator $(M+\pdop_{A,p})^{-1}\circ N$ used in \cref{eq:asym-res-op} is not necessarily single-valued, let alone computationally tractable, for an arbitrary choice of $M\in\R^{n\times n}$, $N\in\R^{n\times d}$ and $A\in\ops_n^F$.
However, for the purposes of this paper, we do not need to find general conditions for when $(M+\pdop_{A,p})^{-1}\circ N$ is single-valued or readily evaluated.
The objective is to parameterize frugal splitting operators which are computationally tractable if the evaluations of the forward and backward steps are tractable.
As it turns out, the \emph{kernel} $M$ in a representation of a frugal splitting operator is highly structured, making the single-valuedness of $(M+\pdop_{A,p})^{-1}\circ N$ easy to establish.
\begin{defin}[$p$-Kernel over $\ops_n^F$]\label{def:p-kernel}%
	We call a matrix $M \in \R^{n\times n}$ a \emph{$p$-kernel over $\ops_n^F$} if $p\notin F$ and $M+\pdskew_p$ is lower triangular with $M_{i,i} \geq 0$ for $i\in \{1,\dots,n\}$ and $M_{i,i} = 0$ if and only if $i\in F$ where $M_{i,i}$ denotes the $i$th element on the diagonal of $M$.
	The matrix $\pdskew_p$ is defined in \cref{eq:primdual}.
\end{defin}
Requiring \(M + \pdskew_p\) to be lower triangular with non-negative diagonal entries enables a substitution that resembles Gaussian elimination. This allows us to reduce the question on single-valuedness of generalized primal-dual resolvents to the single-valuedness of the constituent resolvent and direct evaluations.
\begin{prop}\label{prop:asym-res-sv}
	Let $M\in\R^{n\times n}$ be a $p$-kernel over $\ops_n^F$. Then
	\[
		(M + \pdop_{A,p})^{-1} \colon \Hprim^n \to 2^{\Hprim^n}
	\]
	is single-valued---and hence has full domain---for all $A \in\ops_n^F$.
	\begin{proof}
		Let $A\in\ops_n^F$, $z = (z_1,\dots,z_n) \in\Hprim^n$, $y = (y_1,\dots,y_n) \in\Hprim^n$, and $L = M + \pdskew_p$.
		By definition we have that $y \in (M + \pdop_{A,p})^{-1}x$ is equivalent to
		\[
			x \in (M + \pdop_{A,p})y = (L + \pddiag_{A,p})y.
		\]
		Since $L$ is lower triangular, see \cref{def:p-kernel}, this can be written as
		\begin{align*}
			x_1 &\in L_{1,1}y_1 + \widehat{A}_1y_1 \\
			x_2 &\in L_{2,1}y_1 + L_{2,2}y_2 + \widehat{A}_2y_2 \\
			&\,\,\,\vdots \\
			x_n &\in L_{n,1}y_1 + L_{n,2}y_2 + \dots + L_{n,n}y_n + \widehat{A}_ny_n
		\end{align*}
		where $L_{i,j}$ is the $i,j$th element of $L$ and $\widehat{A}_i = A_i^{-1}$ for all $i\in\{1,\dots,n\}\setminus \{p\}$ and $\widehat{A}_p = A_p$.
		Note that $\widehat{A}_i$ is maximally monotone for all $i\in\{1,\dots,n\}$ since $A\in\ops_n^F$.
		We get
		\begin{equation}\label{eq:back-solv-system}
			\begin{aligned}
				y_1 &\in (L_{1,1}\Id + \widehat{A}_1)^{-1}x_1 \\
				y_2 &\in (L_{2,2}\Id + \widehat{A}_2)^{-1}(x_2 - L_{2,1}y_1) \\
				y_3 &\in (L_{3,3}\Id + \widehat{A}_3)^{-1}(x_3 - L_{3,1}y_1 - L_{3,2}y_2) \\
				&\,\,\,\vdots \\
				y_n &\in (L_{n,n}\Id + \widehat{A}_n)^{-1}(x_n - \sum_{j=1}^{n-1}L_{n,j}y_j).
			\end{aligned}
		\end{equation}
		For all $i\in F$, $L_{i,i} = 0$ and $i \neq p$ by \cref{def:p-kernel} and hence
		$
			(L_{i,i}\Id + \widehat{A}_i)^{-1} = (A_i^{-1})^{-1} = A_i,
		$
		which is single-valued since $A \in \ops_n^F$.
		For all $i\in\{1,\dots,n\}\setminus F$, $L_{i,i} > 0$ and
		$
			(L_{i,i}\Id + \widehat{A}_i)^{-1}  = \resolv_{L_{i,i}^{-1}\widehat{A}_i}\circ(L_{i,i}^{-1}\Id)
		$
		which also is single-valued since $\widehat{A}$ is maximally monotone.
		Therefore, regardless of $x_1,y_1, \ldots,x_n,y_n\in\Hprim$, the right-hand sides of all lines in \cref{eq:back-solv-system} are always singletons and the inclusions can be replaced by equalities.
		Furthermore, in \cref{eq:back-solv-system} we see that $x_1$ uniquely determines $y_1$, which in turn implies that $x_1$ and $x_2$ uniquely determine $y_2$, and so forth.
		Hence, for all $x = (x_1,\dots,x_n)\in\Hprim^n$ there exists a unique $y = (y_1,\dots,y_n)\in\Hprim^n$ such that \cref{eq:back-solv-system} holds.
		Since \cref{eq:back-solv-system} is equivalent to $y \in (M + \pdop_{A,p})^{-1}x$ we conclude that $(M + \pdop_{A,p})^{-1}$ is single-valued.
	\end{proof}
\end{prop}

\section{Representation of Frugal Splitting Operators}\label{sec:representation}
This section will be devoted to the equivalence between frugal splitting operators and generalized primal-dual resolvents whose representations satisfy certain conditions.
The following representation theorem presents our results.
\begin{thm}\label{thm:iff-representation}
	Let $F \subset \{1,\dots,n\}$ and $p \in \{1,\dots,n\} \setminus F$.
	An operator $T_{(\cdot)}\colon \allowbreak \Hprim^d \to \Hprim^d$ parameterized by $\ops_n^F$ is a frugal splitting operator over $\ops_n^F$ if and only if it is a generalized primal-dual resolvent with representation $(p,M,N,U,V)$ where
	\begin{enum_roman}
	\item\label{thm:iff:kernel} $M$ is a $p$-kernel over $\ops_n^F$,
	\item\label{thm:iff:fixedpoint-matrix} $\ker \begin{bmatrix} N & -M\end{bmatrix} \supseteq \ker \begin{bmatrix} U & - V\end{bmatrix}$,
	\item\label{thm:iff:zero-matrix} $\ran U \supseteq \ran V$,
	\end{enum_roman}
	and $M\in\R^{n\times n}$, $N\in\R^{n\times d}$, $U\in\R^{d\times d}$, and $V\in\R^{d\times n}$.

	\begin{proof}
%
		The proof will be carried out in several steps.
		\begin{enumerate}
			\item \label{itm:proof:thm:iff-representation:rep} We will start by showing that any frugal splitting operator is a generalized primal-dual resolvent that admits a representation.

				Roughly speaking, the procedure is as follows: select a primal index; apply Moreau's identity to all backward steps that do not correspond to the primal index; explicitly define the results of the forward evaluations and backward steps; rearrange and identify the primal-dual operator and the generalized primal-dual resolvent representation.
			\item \label{itm:proof:thm:iff-representation:kernel} Then, we will show that the \(M\) in the representation that we constructed in step \cref{itm:proof:thm:iff-representation:rep} is a \(p\)-kernel. This will show \cref{thm:iff:kernel}.
			\item \label{itm:proof:thm:iff-representation:fixedpoint-matrix} We will show the stronger statement that \emph{any} representation of a generalized primal-dual resolvent satisfies \cref{thm:iff:fixedpoint-matrix}, not restricting ourselves to the construction in Step \cref{itm:proof:thm:iff-representation:rep}.
			\item \label{itm:proof:thm:iff-representation:zero-matrix} We will show that \emph{any} representation of a generalized primal-dual resolvent satisfies \cref{thm:iff:zero-matrix}.
			\item \label{itm:proof:thm:iff-representation:surjective}Finally, we will reverse the construction in Step~\cref{itm:proof:thm:iff-representation:rep} and show that any generalized primal-dual resolvent that satisfies \cref{thm:iff:kernel,thm:iff:fixedpoint-matrix,thm:iff:zero-matrix} is a frugal splitting operator.
		\end{enumerate}

		\paragraph{Step~\cref{itm:proof:thm:iff-representation:rep}. A frugal splitting operator admits a representation}
		Let \(T_{(\cdot)} \colon \Hprim^d \to \Hprim^d\) be a frugal splitting operator over \(\ops_n^F\), let $p \in \{1,\dots,n\} \setminus F$, let $A = (A_1,\dots,A_n) \in\ops_n^F$, let $z\in \Hprim^d$ be arbitrary, and consider the evaluation of $T_Az$.

		From \cref{def:fo-split}\cref{def:fo-split:comp} we know that $T_{A} = \tau_{n,A} \circ \cdots \circ \tau_{1,A}$ where $\tau_{i,A} \colon \Hprim^{d_i}\to \Hprim^{d_{i+1}}$ and $d_1 = d_{n+1} = d$.
		If we introduce variables for the intermediate results, i.e., $z_1 = z$ and $z_{i+1} = \tau_{i,A}z_i$ for all $i\in\{1,\dots,n\}$, we can write $T_Az = z_{n+1}$.
		Note that the variables $z_i$ of this proof should not be confused with the variables of a fixed point iteration \cref{eq:fp-iteration} of $T_A$.
		Furthermore, from the definition of $\tau_{i,A}$ in \cref{def:fo-split}\cref{def:fo-split:comp} we can conclude that for all $i\in\{1,\dots,n\}$ there exist matrices $B_i\in\R^{d_{i+1}\times d_i}$, $C_i \in \R^{d_{i+1}\times 1}$ and $D_i \in \R^{1 \times d_i}$ (that do not depend on \(A\)) such that
		\[
			z_{i+1} = \tau_{i,A}z_i =
			\begin{cases}
				B_iz_i + C_i \resolv_{\gamma_i A_i} D_i z_i & \text{if } i \notin F, \\
				B_iz_i + C_i A_i D_i z_i & \text{if } i \in F.
			\end{cases}
		\]
		Now, for all $i \in \{1,\dots,n\}\setminus (F \cup \{p\})$, we apply the Moreau identity, $\resolv_{\gamma_i A_{i}} = \Id - \gamma_i \resolv_{\gamma_i^{-1}A_{i}^{-1}} \circ \gamma_i^{-1}\Id$, and rewrite $z_{i+1} = \tau_{i,A}z_i$ as
		\begin{align*}
			z_{i+1}
			&= B_iz_i + C_i \resolv_{\gamma_i A_i} D_i z_i \\
			&= (B_i + C_iD_i)z_i - (\gamma_iC_i) \resolv_{\gamma_i^{-1}A_i^{-1}} (\gamma_{i}^{-1}D_i) z_i \\
			&= \widetilde{B}_iz_i + \widetilde{C}_i (l_{i,i}\Id + \widehat{A}_i)^{-1} \widetilde{D}_i z_i
		\end{align*}
		where $\widetilde{B}_i = B_i + C_iD_i$, $\widetilde{C}_i = -\gamma_i C_i$, $\widetilde{D}_i = D_i$, $l_{i,i} = \gamma_i$ and $\widehat{A}_i = A_i^{-1}$ (see also \cref{eq:Delta-and-hats}).
		We can write $\tau_{i,A}z_i$ in a similar form for all $i\in F$ as well,
		\begin{align*}
			z_{i+1}
			&= B_iz_i + C_i A_i D_i z_i \\
			&= \widetilde{B}_iz_i + \widetilde{C}_i (l_{i,i}\Id + \widehat{A}_i)^{-1}\widetilde{D}_i z_i
		\end{align*}
		where $\widetilde{B}_i = B_i$, $\widetilde{C}_i = C_i$, $\widetilde{D}_i = D_i$, $l_{i,i} = 0$ and $\widehat{A}_i = A_i^{-1}$.
		Finally, since $p \notin F$, we can write $\tau_{p,A}z_p$ as
		\begin{align*}
			z_{p+1}
			&= B_pz_p + C_p \resolv_{\gamma_p A_p} D_p z_p \\
			&= B_pz_p + C_p (\gamma_p^{-1}\Id +  A_p)^{-1} (\gamma_p^{-1}D_p) z_p \\
			&= \widetilde{B}_pz_p + \widetilde{C}_p (l_{p,p}\Id + \widehat{A}_p)^{-1} \widetilde{D}_p z_p
		\end{align*}
		where $\widetilde{B}_p = B_p$, $\widetilde{C}_p = C_p$, $\widetilde{D}_p = \gamma_p^{-1}D_p$, $l_{p,p} = \gamma_p^{-1}$ and $\widehat{A}_p = A_p$.
		With these notations in place, we define $y_i \in \Hprim$ as
		\begin{equation}\label{eq:re:normalized-composition}
			y_i = (l_{i,i} \Id + \widehat{A}_i)^{-1}\widetilde{D}_iz_i
			,\quad
			\forall i \in \{1,\dots,n\}
		\end{equation}
		which gives
		\[
			z_{i+1} = \widetilde{B}_iz_i + \widetilde{C}_iy_i
			,\quad
			\forall i\in \{1,\dots,n\}
		\]
		and for clarity we note that $\widetilde{B}_i\in\R^{d_{i+1}\times d_i}$, $\widetilde{C}_i \in \R^{d_{i+1}\times 1}$ and $\widetilde{D}_i \in \R^{1 \times d_i}$ for all $i\in \{1,\dots,n\}$.
		Unrolling this iteration gives
		\[
			z_{i} = \widetilde{B}_{(i-1)\dots 1} z + \sum_{j=1}^{i-1}\widetilde{B}_{(i-1)\dots(j+1)}\widetilde{C}_jy_j
			,\quad
			\forall i\in \{1,\dots,n+1\}
		\]
		where $\widetilde{B}_{a\dots b} = \widetilde{B}_a\widetilde{B}_{a-1}\dots \widetilde{B}_{b}$ and $\widetilde{B}_{a\dots b} = I$ if $a < b$ where $I\in\R^{d_b\times d_b}$ is the identity matrix.
		The expression for $y_i$ in \cref{eq:re:normalized-composition} can then be rewritten as
		\[
			l_{i,i}y_i + \widehat{A}_iy_i \ni \widetilde{D}_i\widetilde{B}_{(i-1)\dots 1}z + \sum_{j=1}^{i-1}\widetilde{D}_i\widetilde{B}_{(i-1)\dots(j+1)}\widetilde{C}_jy_j
			,\quad
			\forall i\in \{1,\dots,n\}.
		\]
		If we define $N_i = \widetilde{D}_i\widetilde{B}_{(i-1)\dots 1} \in \R^{1\times d}$, $l_{i,j} = -\widetilde{D}_i\widetilde{B}_{(i-1)\dots(j+1)}\widetilde{C}_j \in \R$ for $j < i$ and rearrange this expression we get
		\[
			\sum_{j=1}^i l_{i,j}y_j + \widehat{A}_iy_i \ni N_i z
			,\quad
			\forall i\in \{1,\dots,n\}.
		\]
		If we define $y = (y_1,\dots,y_n) \in \Hprim^n$, the lower triangular matrix $L \in\R^{n\times n}$ with elements $L_{i,j} = l_{i,j}$ for all $j \leq i$ and the matrix $N \in \R^{n\times d}$ whose $i$th row is given by $N_i$, this can be written as $Ly + \pddiag_{A,p}y \ni Nz$ or equivalently
		\[
			y = (L + \pddiag_{A,p})^{-1}Nz = (M + \pdop_{A,p})^{-1}Nz
		\]
		where $M = L - \pdskew_p$.

		Finally, if we define the matrices $U = I - \widetilde{B}_{n\dots 1} \in \R^{d\times d}$ and $V\in \R^{d\times n}$ where the $j$th column of $V$ is given by $\widetilde{B}_{n\dots(j+1)}\widetilde{C}_j$ we can write the expression of $z_{n+1}$ as
		\[
			T_{A}z = z_{n+1} = z - Uz + Vy
		\]
		which shows that $(p,M,N,U,V)$ is a representation of $T_{(\cdot)}$.

		\paragraph{Step~\cref{itm:proof:thm:iff-representation:kernel}. \(M\) is a \(p\)-kernel}
		By construction, \(M + \pdskew_p = L\) is a lower triangular matrix.
		Moreover, since the diagonal elements of $M$ are $M_{i,i} = l_{i,i} \geq 0$ and $l_{i,i} = 0$ if and only if $i \in F$, from \cref{def:p-kernel} we conclude that $M$ is a $p$-kernel over $\ops_n^F$, i.e., that \cref{thm:iff:kernel} holds.

		\paragraph{Step~\cref{itm:proof:thm:iff-representation:fixedpoint-matrix}. Any representation satisfies \cref{thm:iff:fixedpoint-matrix}}
		Let $K\colon\Hprim \to \Hprim$ be a nonzero bounded linear skew-adjoint operator, and let $v \in \Hprim$ be such that $Kv \neq 0$ and hence $v \neq 0$.
		Such an operator $K$ and point $v$ exist due to the assumption that $\dim \Hprim \geq 2$.

		We show the necessity of \cref{thm:iff:fixedpoint-matrix} by providing a counterexample.
		In particular we will show that, if $U\hat{z}-V\hat{y} = 0$ while $N\hat{z}-M\hat{y} \neq 0$ for some $\hat{z}\in\R^d$ and $\hat{y}\in\R^n$,
		we can always construct some operator tuple $A \in \ops_n^F$ such that the equivalence $\fix T_A \neq \emptyset \iff \zer \sum_{i=1}^n A_i \neq \emptyset$ of \cref{def:fo-split}\cref{def:fo-split:fix} leads to a contradiction.

		Assume that $\hat{z}\in \R^d$ and $\hat{y}\in \R^n$ are such that $U\hat{z}-V\hat{y} = 0$ but $N\hat{z}-M\hat{y} \neq 0$.
		Then $\tilde z = (\hat{z}_1v,\dots,\hat{z}_dv) \in \Hprim^d$ and $\tilde y = (y_1,\dots,y_n) = (\hat{y}_1v,\dots,\hat{y}_nv) \in \Hprim^n$ are such that $U\tilde z-V\tilde y = 0$ and $N\tilde z-M\tilde y\neq 0$.
		Define $(\tilde a_1,\dots,\tilde a_n) = N\tilde z - M\tilde y \neq 0$ and note that all of $\tilde a_1,\dots,\tilde a_n$ are parallel to $v$.
		Let $l \in\{1,\dots,n\}\setminus \{p\}$ and define $A = (A_1,\dots,A_n)\in\ops_n^F$ where
		\[
			A_ix =
			\begin{cases}
				K(x-\tilde y_p) + \tilde a_p - \sum_{j\in \{1,\dots,n\} \setminus \{p\} } \tilde y_j & \text{if } i = p, \\
				\tilde y_l - K(x-\tilde y_p-\tilde a_l) & \text{if } i = l, \\
				\tilde y_i & \text{otherwise}
			\end{cases}
		\]
		for all $x\in \Hprim$ and all $i\in\{1,\dots,n\}$.
		Note that both $K$ and $-K$ are maximally monotone since $K$ is skew-adjoint, hence, $A_i$ is maximally monotone for all $i\in\{1,\dots,n\}$.
		The primal-dual operator $\pdop_{A,p}$ of this tuple evaluated at the $\tilde y$ specified above satisfies
		\begin{multline*}
			N\tilde z-M\tilde y = (\tilde a_1,\dots,\tilde a_n) \in \pdop_{A,p}\tilde y = \prod_{i = 1}^n \left\lbrace \begin{aligned}
					&\lbrace \tilde a_p \rbrace&&\text{if } i = p \\
					&\tilde a_l + \ker K &&\text{if } i = l \\
					&\Hprim &&\text{otherwise}
				\end{aligned}\right\rbrace
				\\
			\text{and hence}
			\quad
			\tilde y = (M + \pdop_{A,p})^{-1}N\tilde z
		\end{multline*}
		since $(M+ \pdop_{A,p})^{-1}$ is single-valued by \cref{def:gen-pd-resolv}.
		We then have $T_{A}\tilde z = \tilde z - U\tilde z + V\tilde y$ and since $U\tilde z-V\tilde y=0$ by assumption we have $T_A\tilde z = \tilde z$ and $\fix T_A \neq \emptyset$.
		\cref{def:fo-split}\cref{def:fo-split:fix} then implies $\zer \sum_{i=1}^nA_i \neq \emptyset$ and since
		\[
			\sum_{i=1}^nA_ix = \tilde a_p + K\tilde a_l
		\]
		for all $x\in\Hprim$ must $K\tilde a_l = -\tilde a_p$.
		Since $\tilde a_p$ and $\tilde a_l$ are parallel to $v$, there exist $\lambda_p,\lambda_l\in\R$ such that $\tilde a_p = \lambda_p v$ and $\tilde a_l = \lambda_l v$.
		Furthermore, since $K$ is skew-adjoint---hence $\inprod{Kx}{x} = 0$ for all $x\in\Hprim$---
		\[
			0 = \inprod{K\tilde a_l}{\tilde a_l} = \inprod{-\tilde a_p}{\tilde a_l} = \inprod{-\lambda_pv}{\lambda_l v} = -\lambda_p\lambda_l\norm{v}^2
		\]
		and, since $v\neq 0$, at least one of $\lambda_p$ and $\lambda_l$ is zero.
		In fact, both must be zero since
		\[
			-\lambda_p v = -\tilde a_p = K\tilde a_l = \lambda_l Kv
		\]
		and both $v \neq 0$ and $Kv \neq 0$.
		This means that $\tilde a_p = \tilde a_l = 0$ but, since $l \in\{1,\dots,n\}\setminus \{p\}$ was arbitrary, we must have $\tilde a_i = 0$ for all $i\in\{1,\dots,n\}$ and hence $(\tilde a_1,\dots,\tilde a_n) = 0$ which is a contradiction.
		This concludes the necessity of \cref{thm:iff:fixedpoint-matrix}.

		\paragraph{Step~\cref{itm:proof:thm:iff-representation:zero-matrix}. Any representation satisfies \cref{thm:iff:zero-matrix}}
		For the necessity of \cref{thm:iff:zero-matrix}, let $\tilde y = (\tilde y_1,\dots,\tilde y_n) \in \Hprim^n$ be arbitrary and define $A = (A_1,\dots,A_n)\in\ops_n^F$ where
		\begin{align*}
			A_ix =
			\begin{cases}
				x - \tilde y_p + \tilde y_i & \text{if } i \neq p, \\
				x - \sum_{j=1}^n \tilde y_j & \text{if } i = p
			\end{cases}
		\end{align*}
		for all $x\in\Hprim$ and all $i\in\{1,\dots,n\}$.
		These operators satisfy $\{\tilde y_p\} = \zer \sum_{i=1}^n A_i$ and $\{\tilde y\} = \zer \pdop_{A,p}$.
		Since $T_{(\cdot)}$ is a frugal splitting operator we must then have $\emptyset \neq \fix T_{A}$.
		Let $z^\star \in \fix T_{A}$, \cref{eq:asym-res-op} then implies the existence of some $y'\in\Hprim^d$ such that $Uz^\star - Vy' = 0$ and $Nz^\star-My' \in \pdop_{A,p}y'$.
		But \cref{thm:iff:fixedpoint-matrix} (which has been shown in Step~\cref{itm:proof:thm:iff-representation:fixedpoint-matrix}) implies that $0 = Nz^\star-My' \in \pdop_{A,p}y'$ and, therefore, $y' \in \zer \pdop_{A,p} = \{\tilde y\}$, i.e., $y' = \tilde y$ and $Uz^\star = V\tilde y$.
		Since the choice of $\tilde y \in \Hprim^n$ is arbitrary this implies that for all $y \in \Hprim^n$ there exists $z\in\Hprim^d$ such that $Uz = Vy$ which then implies the existence of $\hat{z}\in \R^d$ for each $\hat{y}\in\R^n$ such that $U\hat{z} = V\hat{y}$, i.e., $\ran U \supseteq \ran V$.



		\paragraph{Step~\cref{itm:proof:thm:iff-representation:surjective}. Any representation satisfying \cref{thm:iff:kernel,thm:iff:fixedpoint-matrix,thm:iff:zero-matrix} represents a frugal splitting operator}
		Let $T_{(\cdot)}\colon \Hprim^d \to \Hprim^d$ be a generalized primal-dual resolvent over $\ops_n^F$ with representation $(p,M,N,U,V)$ and let $A\in \ops_n^F$.
		For each $z\in\Hprim^d$ there exists $y\in\Hprim^n$ such that
		\begin{equation}\label{eq:rs:opdef}
			\left\{
				\begin{aligned}
					y &= (M + \pdop_{A,p})^{-1}Nz, \\
					T_{A}z &= z - Uz + Vy
				\end{aligned}
			\right.
			\quad\text{or equivalently}\quad
			\left\{
				\begin{aligned}
					Nz-My &\in \pdop_{A,p}y, \\
					Uz - Vy &= z - T_{A}z.
				\end{aligned}
			\right.
		\end{equation}
		Consider \cref{def:fo-split}\cref{def:fo-split:fix}.
		If $z\in\fix T_A$, then $Uz-Vy = 0$ and \cref{thm:iff:fixedpoint-matrix} implies $0 = Nz - My \in \pdop_{A,p}y$ and hence $y \in \zer \pdop_{A,p} \neq \emptyset$.
		This proves the right implication of \cref{def:fo-split}\cref{def:fo-split:fix}.

		For the left implication of \cref{def:fo-split}\cref{def:fo-split:fix}, let $y^\star\in \zer \pdop_{A,p}$.
		Then \cref{thm:iff:zero-matrix} implies that there exists $z\in\Hprim^d$ such that $Uz-Vy^\star = 0$ and \cref{thm:iff:fixedpoint-matrix} then implies $Nz-My^\star = 0$.
		Since $y^\star \in \zer \pdop_{A,p}$ we then have $0 = Nz - My^\star \in \pdop_{A,p}y^\star$ and $z$ and $y^\star$ then satisfy \cref{eq:rs:opdef} which proves that $z\in\fix T_A \neq \emptyset$.

		To show that \cref{def:fo-split}\cref{def:fo-split:comp} is implied by item \cref{thm:iff:kernel} we first introduce some notation.
		Let $V_i \in \R^{d\times 1}$ denote the $i$th column of $V$, $N_i \in \R^{1\times d}$ denote the $i$th row of $N$, and $l_{i,j}$ denote the $i,j$-element of the matrix $L = M +\pdskew_p$.
		Define the matrices
		\begin{gather*}
			\widetilde{B}_1 =
			\big[
				\underbrace{I}_{\R^{d\times d}} \mid
				\underbrace{I}_{\R^{d \times d}} \mid
				\underbrace{\mathbf{0}}_{\R^{d\times d}} \mid
				\underbrace{\mathbf{0}}_{\R^{d\times n}}
			\big]^T
			,\quad
			\widetilde{B}_n =
			\big[
				\underbrace{I-U}_{\R^{d\times d}} \mid
				\underbrace{\mathbf{0}}_{\R^{d \times d}} \mid
				\underbrace{I}_{\R^{d\times d}} \mid
				\underbrace{\mathbf{0}}_{\R^{d\times n}}
			\big],
		\end{gather*}
		and $\widetilde{B}_{i} = I \in \R^{(3d + n) \times (3d+n)}$ for all $i \in \{2,\dots,n-1\}$ where $I$ and $\mathbf{0}$ denote identity and zero matrices of appropriate sizes, respectively.
		Further define the matrices
		\begin{align*}
			\widetilde{C}_i =
			\big[
				\underbrace{\mathbf{0}}_{\R^{1\times d}} \mid
				\underbrace{\mathbf{0}}_{\R^{1\times d}} \mid
				V_i^T \mid
				\underbrace{0\dots 0}_{\R^{1\times (i-1)}}\;1 \; \underbrace{0\dots 0}_{\R^{1\times (n-i)}}
			\big]^T
			&\quad\text{for all } i\in \{1,\dots,n-1\}
			,\\
			\widetilde{D}_i =
			\big[
				\underbrace{\mathbf{0}}_{1\times d} \mid
				N_i \mid
				\underbrace{\mathbf{0}}_{\R^{1\times d}} \mid
				{-l_{i,1} \dots -l_{i,n}}
			\big]
			&\quad\text{for all } i\in \{2,\dots,n\},
		\end{align*}
		$\widetilde{C}_n = V_n$ and $\widetilde{D}_1 = N_1$.
		Now, for $z_1\in\Hprim^d$ define
		\[
			z_{i+1} = \widetilde{B}_iz_i + \widetilde{C}_i(l_{i,i}\Id + \widehat{A}_i)^{-1}\widetilde{D}_iz_i
		\]
		for $i\in\{1,\dots,n\}$ where $\widehat{A}_i = A_i^{-1}$ for all $i\in\{1,\dots,n\}\setminus \{p\}$ and $\widehat{A}_p = A_p$.
		By following the same procedure as in the proof of step~\cref{itm:proof:thm:iff-representation:rep}, starting at \cref{eq:re:normalized-composition}, it can be verified that $z_{n+1} = T_{A}z_1$ if \cref{thm:iff:kernel} holds.
		Furthermore, by inverting the relationship between \(B_i\), \(C_i\), and \(D_i\) on the one hand and \(\widetilde B_i\), \(\widetilde C_i\), and \(\widetilde D_i\) on the other in the proof of step~\cref{itm:proof:thm:iff-representation:rep} that led to \cref{eq:re:normalized-composition} we can conclude that there exist real matrices $B_i$, $C_i$, $D_i$ such that
		\[
			z_{i+1} =
			\begin{cases}
				B_iz_i + C_i \resolv_{\gamma_i A_i} D_i z_i & \text{if } i \notin F, \\
				B_iz_i + C_i A_i D_i z_i & \text{if } i \in F,
			\end{cases}
		\]
		for all $i\in\{1,\dots,n\}$.
		This proves that \cref{def:fo-split}\cref{def:fo-split:comp} holds.
	\end{proof}
\end{thm}
\ifdefined\shortversion
	Note that we require $F$ to be a strict subset of $\{1,\dots,n\}$ which implies that there always exists $p \in \{1,\dots,n\} \setminus F$.
	This generalized primal-dual resolvent representation is advantageous to work with over the original frugal splitting operator definition since the primal-dual operator appears in the representation and one only needs to consider range and structure constraints of a handful of matrices.
	In \cref{sec:minlift,sec:conv,sec:new-method}, we will use this to analyze the lifting and convergence of general splittings as well as constructing a new splitting operator.
	In the remainder of this section, we will explore some of the more direct consequences of \cref{thm:iff-representation}.
\else
	Note that we require $F$ to be a strict subset of $\{1,\dots,n\}$ which implies there always exists $p \in \{1,\dots,n\} \setminus F$.
	It is also worth remembering that it is assumed that $\dim \Hprim \geq 2$.
	However, this assumption is only ever required for the proof of \cref{lem:rep-mat-prop} and, since this lemma only concerns the ``only if'' part of \cref{thm:iff-representation}, the sufficient conditions are not affected.
	Hence, any representation $(p,M,N,U,V)$ that satisfies \cref{thm:iff-representation} yields a frugal splitting operator in the $\dim \Hprim \leq 1$ setting as well and all examples of frugal splitting operators in this paper are frugal splitting operators regardless of the dimension of $\Hprim$.
	We have succeeded in finding replacements for the counterexample in the proof of \cref{lem:rep-mat-prop} that require $\dim \Hprim \geq 2$, but not without relaxing some other assumption.
	For instance, when $\dim \Hprim = 1$ we have been able to construct sufficient counterexamples if we instead of single-valued operators allow for at most single-valued operators in the tuples of $\ops_n^F$.
	However, this relaxation complicates questions regarding the domain of our generalized primal-dual resolvent and we believe the $\dim \Hprim = 1$ case is of too limited practical interest to warrant these complications.
	Similarly, the trivial $\dim \Hprim = 0$ case is also not worth handling.

	The fact that the primal-dual operator appears directly in the representation and one only needs to consider simple range and structure constraints of a handful of matrices makes the representation easy to work with, see for instance \cref{sec:minlift,sec:conv,sec:new-method} where we analyze general splittings and construct a new splitting.
	The representation also makes the relationship between fixed points of a splitting and solutions to \cref{eq:prob} clearly visible, something we will illustrate in the remainder of this section.
	There we will also discuss and illustrate different properties of the representation, such as how the step-sizes and parallelizability of the forward and backward evaluations are encoded.
	We will also demonstrate an approach for deriving a representation of a frugal splitting operator and provide an alternative factorization of a representation $(p,M,N,U,V)$.
\fi

\subsection{Alternative Factorization}
\Cref{thm:iff-representation} gives conditions on all four matrices of a representation $(p,M,N,U,V)$.
However, it turns out that these conditions make the kernel $M$ uniquely defined given $N$, $U$ and $V$.\footnote{%
	If \(SU = N = N' = S' U'\), \(U = U'\), and \(UP = V = V' = U' P' \), then \(M = SUP\) and \(M' = S' U' P'\) implies
	\[
		M' = S' U' P' = S' UP = S' U' P = SUP = M.
	\]
}
This leads to the following corollary that will be used both in the examination of minimal lifting in \cref{sec:minlift} and in the convergence analysis in \cref{sec:conv}.
\begin{cor}\label{cor:factor-rep}
	Let $F\subset \{1,\dots,n\}$ and $p \in \{1,\dots,n\} \setminus F$.
	A generalized primal-dual resolvent with representation $(p,SUP,SU,U,UP)$ where $U\in\R^{d\times d}$, $S\in\R^{n\times d}$ and $P\in\R^{d\times n}$ is a frugal splitting operator if
	\begin{enum_roman}
	\item \label{thm:fr:kernel} $SUP$ is a $p$-kernel over $\ops_n^F$,
	\item \label{thm:fr:ran-ker} $\ran P \subseteq (\ker U)^\perp$ and $\ker S \supseteq (\ran U)^\perp$.
	\end{enum_roman}
	Furthermore, for any frugal splitting operator over $\ops_n^F$ with representation $(p,M,N,\allowbreak U,V)$ that satisfies \cref{thm:iff-representation}, there exist matrices $S\in\R^{n\times d}$ and $P\in\R^{d\times n}$ such that $M = SUP$, $N = SU$ and $V = UP$ and the conditions \cref{thm:fr:kernel,thm:fr:ran-ker} are satisfied.

	\begin{proof}
		It is straightforward to verify that a representation $(p,SUP,SU,U,UP)$ that satisfies the conditions of the theorem satisfies the conditions of \cref{thm:iff-representation} and hence is a representation of a frugal splitting operator over $\ops_n^F$.

		Assume $T_{(\cdot)}\colon\Hprim^d\to\Hprim^d$ is a frugal splitting operator over $\ops_n^F$ with representation $(p,M,N,U,V)$ that satisfies \cref{thm:iff-representation}.
		The theorem states that $\ran U \supseteq \ran V$ and \cref{lem:eq-range-right-fact} then proves the existence of a unique $P\in\R^{d\times n}$ with $\ran P \subseteq (\ker U)^\perp$ such that
		\[
			V = UP.
		\]
		\Cref{thm:iff-representation} also states $\ker \begin{bmatrix}U & -V\end{bmatrix} \subseteq \ker \begin{bmatrix}N & -M\end{bmatrix}$ and \cref{lem:eq-ker-left-fact} yields the existence of a unique $S\in\R^{n\times d}$ with $\ker S \supseteq (\ran \begin{bmatrix}U&-V\end{bmatrix})^\perp$ such that
		\[
			N = SU
			\quad\text{and}\quad
			M = SV = SUP.
		\]
		However, since $\ran U \supseteq \ran V$ we have $\ran \begin{bmatrix}U&-V\end{bmatrix} = \ran U$, this concludes the proof.
	\end{proof}
\end{cor}
We provide a similar factorization in \cref{prop:kernel2split} that may be more useful when constructing new splittings since it avoids the search for matrices $S$, $U$, and $P$ such that $SUP$ is a $p$-kernel.
However, unlike \cref{cor:factor-rep}, it is not guaranteed that all representations of frugal splitting operators have a factorization of the form in \cref{prop:kernel2split}.

\subsection{Solution Map and Fixed Points}\label{sec:sol-and-fix-point}
The following two propositions reveal the relationship between fixed points of the frugal splitting operator and solutions to the monotone inclusion problem.
They also clarify why the assumption of the existence of a solution map is not needed in \cref{def:fo-split}. The reason is that there always exists a mapping from fixed points of a frugal splitting operator to solutions of the monotone inclusion problem \cref{eq:prob}.
Furthermore, this mapping is always evaluated within the evaluation of a frugal splitting operator itself.
\begin{prop}\label{prop:sol-map}
	Let $T_{(\cdot)}\colon \Hprim^d \to \Hprim^d$ be a frugal splitting operator over $\ops_n^F$ with representation $(p,M,N,U,V)$ that satisfies \cref{thm:iff-representation}.
	If $z^\star\in\fix T_A$ for $A\in\ops_n^F$, then the tuple $(y_1^\star,\dots,y_n^\star) \defeq (M + \pdop_{A,p})^{-1}Nz^\star$ is such that
	\[
		(y_1^\star,\dots,y_n^\star) \in \zer\pdop_{A,p}
		\quad\text{and}\quad
		y_p^\star \in \zer \sum_{i=1}^n A_i.
	\]

	\begin{proof}
		Let $A\in\ops_n^F$ and $z^\star\in\fix T_A$, \cref{def:gen-pd-resolv} then gives the existence of $y^\star = (y_1^\star,\dots,y_n^\star)\in\Hprim^n$ such that
		\[
			\left\{
				\begin{aligned}
					y^\star &= (M + \pdop_{A,p})^{-1}Nz^\star, \\
					z^\star &= z^\star - Uz^\star + Vy^\star
				\end{aligned}
			\right.
			\quad\text{or equivalently}\quad
			\left\{
				\begin{aligned}
					Nz^\star-My^\star &\in \pdop_{A,p}y^\star, \\
					Uz^\star - Vy^\star &= 0
				\end{aligned}
			\right.
		\]
		Since \cref{thm:iff-representation}\cref{thm:iff:fixedpoint-matrix} holds, $Uz^\star-Vy^\star = 0$ implies $Nz^\star-My^\star = 0$ and hence $y^\star \in \zer \pdop_{A,p}$.
		That $y_p^\star \in \zer \sum_{i=1}^n A_i$ follows from the equivalence between the primal-dual problem \cref{eq:primdual} and the primal problem \cref{eq:prob}.
	\end{proof}
\end{prop}
\begin{prop}\label{prop:fixed-points}
	Let $T_{(\cdot)}\colon \Hprim^d \to \Hprim^d$ be a frugal splitting operator over $\ops_n^F$ and let $(p,SUP,SU,\allowbreak U,UP)$ be a representation of $T_{(\cdot)}$ that satisfies \cref{cor:factor-rep}.
	With $A\in \ops_n^F$, the set of fixed points satisfies
	\[
		\fix T_A
		\supseteq
		P\zer \pdop_{A,p}
	\]
	where $P\zer \pdop_{A,p} = \{Py^\star \mid y^\star \in \zer \pdop_{A,p}\}$.
	Equality in the inclusion holds if and only if $U$ has full rank.
	\begin{proof}
		Let $A\in\ops_n^F$ and $y^\star \in \zer \pdop_{A,p}$.
		We then have $0 = SUPy^\star - SUPy^\star \in \pdop_{A,p}y$, or equivalently $y^\star = (SUP + \pdop_{A,p})^{-1}SUPy^\star$, since $(SUP + \pdop_{A,p})^{-1}$ is single-valued due to $SUP$ being a $p$-kernel over $\ops_n^F$, see \cref{prop:asym-res-sv}.
		By letting $z^\star = Py^\star$ in the definition of the generalized primal-dual resolvent, we conclude that
		\begin{align*}
			y^\star &= (SUP + \pdop_{A,p})^{-1}SUPy^\star, \\
			T_APy^\star &= Py^\star - U(Py^\star - Py^\star) = Py^\star,
		\end{align*}
		which implies that $Py^\star \in \fix T_A$.
		Since $y^\star\in \zer \pdop_{A,p}$ was arbitrary we have $\fix T_A \supseteq P \zer \pdop_{A,p}$.

		Assume $U$ does not have full rank.
		Then there exists $z\in \Hprim^d$ such that $Uz = 0$ and $z \neq 0$.
		Since $\ran P \subseteq (\ker U)^\perp$ by \cref{cor:factor-rep}, there exists no $y'\in\Hprim^n$ such that $z = Py'$ and hence $Py^\star + z \notin P\zer \pdop_{A,p}$ since $y^\star\in\zer \pdop_{A,p}$.
		However, we have
		\begin{align*}
			y^\star &= (SUP + \pdop_{A,p})^{-1}SUPy^\star = (SUP + \pdop_{A,p})^{-1}SU(Py^\star + z), \\
			T_A(Py^\star+z) &= Py^\star+z - U(Py^\star+z - Py^\star) = Py^\star + z,
		\end{align*}
		and hence $Py^\star + z \in\fix T_A$.
		The equality $\fix T_A = P\zer \pdop_{A,p}$ can therefore not hold if $U$ does not have full rank.

		Assume $U$ has full rank.
		Let $z^\star \in \fix T_A$, then there exists $y^\star\in \Hprim^n$ such that
		\begin{align*}
			SUz^\star - SUPy^\star &\in \pdop_{A,p}y^\star, \\
			U(z^\star - Py^\star) &= z^\star - T_Az^\star = 0.
		\end{align*}
		However, since $U$ has full rank, this implies that $z^\star = Py^\star$ and that $0 = SU(z^\star-Py^\star) \in \pdop_{A,p}y^\star$ and hence that $y^\star\in\zer \pdop_{A,p}$ and $z^\star \in P\zer \pdop_{A,p}$.
		Since $z^\star \in\fix T_A$ was arbitrary we have $\fix T_A \subseteq P\zer\pdop_{A,p}$ and the opposite inclusion from before gives equality.
	\end{proof}
\end{prop}

\subsection{Evaluation Order and Parallelizability}\label{sec:eval-order}
For a frugal splitting operator, the order of the constituent operator evaluations is specified in \cref{def:fo-split}.
However, the definition only states that it should be possible to evaluate the frugal splitting operator by evaluating the constituent operators in this order. It does not exclude the possibility of computing the frugal splitting operator with another evaluation order or with some of the evaluations being performed in parallel.
Especially, the ability of performing forward and/or backward evaluations in parallel is of particular interest since it can allow for distributed or multi-threaded implementations.
Fittingly, how, and if, one operator evaluation depends on previous evaluations is directly encoded in the kernel.
It is therefore straightforward to both identify and construct parallelizable kernels. An example is found in \cref{sec:new-method}, where we construct a new parallelizable frugal splitting operator with minimal lifting.

Let $L = M + \pdskew_p$ for a representation $(p,M,N,U,V)$ of a frugal splitting operator over $\ops_n^F$. Then the inverse in the generalized primal-dual resolvent satisfies
$(M + \pdop_{A,p})^{-1} = (L + \pddiag_{A,p})^{-1}$.
Since $M$ is a $p$-kernel, $L$ is lower triangular and $(L+\pddiag_{A,p})^{-1}$ can be evaluated with back-substitution (as in \cref{eq:back-solv-system} in the proof of \cref{prop:asym-res-sv}).
This means that the strict lower triangular part of $L$, $L_{i,j}$ for $j < i$, determines the dependency on the results of previous forward or backward evaluations%
\ifdefined\shortversion
	{}.
\else
	\footnote{
		The strictly lower triangular matrix $\Tilde{L}\in\R^{n\times n}$ with elements $\Tilde{L}_{i,j} = 1$ if $L_{i,j} \neq 0$ and $\Tilde{L}_{i,j} = 0$ if $L_{i,j} = 0$ for $j < i$
		is the transpose of the adjacency matrix of the directed dependency graph of the operator evaluations, i.e., the graph of $n$ nodes where there is an edge from the $i$th to $j$th node if the result of forward or backward evaluation of $A_i$ is used in the argument of the forward or backward evaluation of $A_j$.
	}.
\fi
Hence, if for $i \in \N$ and some $j < i$ the element $L_{i,j} = 0$, then the $i$th evaluation does not directly depend on the result of the $j$th one.
If the $i$th evaluation does not depend on any other evaluation that depends on the $j$th evaluation, then the $i$th and $j$th evaluation can be performed in parallel.

For example, if we take the $4$-kernel over $\ops_4^{\{2\}}$ as
\[
	M
	= \begin{bmatrix}
		1 & 0 & 0 & 1 \\
		1 & 0 & 0 & 1 \\
		0 & 1 & 1 & 1 \\
		-1 & 0 & -1 & 1 \\
	\end{bmatrix}
	\quad\text{which yields}\quad
	L
	= M + \pdskew_{4}	= \begin{bmatrix}
		1 & 0 & 0 & 0 \\
		1 & 0 & 0 & 0 \\
		0 & 1 & 1 & 0 \\
		0 & 1 & 0 & 1 \\
	\end{bmatrix}
\]
we conclude that, in any frugal splitting operator with this kernel, the second evaluation directly depends on the first; the third depends directly on the second only but depends indirectly on the first; the fourth is independent of the third, depends directly on the second, and depends indirectly on the first.
The second evaluation therefore needs to be performed after the first, while the last two evaluations need to take place after the second but potentially in parallel with each other.

\subsection{Example of a Representation}\label{sec:rep-ex}
The process of finding a representation of a frugal splitting operator is quite straightforward and follows the procedure indicated in Step \cref{itm:proof:thm:iff-representation:rep} of the proof of \cref{thm:iff-representation}.
We demonstrate the process on the three operator splitting of Davis and Yin \cite{davisThreeOperatorSplittingScheme2017} which encodes the zeros of $A_1 + A_2 + A_3$ where $A = (A_1,A_2,A_3) \in \ops_3^{\{2\}}$ as the fixed points of the operator
\[
	T_{A} = \resolv_{\gamma A_3} \circ (2\resolv_{\gamma A_1} - \Id  - \gamma A_2 \circ \resolv_{\gamma A_1}) + \Id  - \resolv_{\gamma A_1}.
\]
\ifdefined\shortversion
\else
In order for fixed point iterations of $T_{A}$ to always be convergent it is further required that $A_2$ is cocoercive and that the step-size $\gamma$ is sufficiently small.
However, we leave the convergence analysis of frugal splitting operators to \cref{sec:conv}.
\fi
Consider the evaluation of $T_{A}z$ where $z\in \Hprim$.
We first choose the primal index $p$ in the representation $(p,M,N,U,V)$.
Since it is required that $p \notin \{2\}$ we have the choice of either $p = 1$ or $p = 3$.
We choose $p = 3$.
Applying Moreau's identity to the resolvents of $A_i$ for all $i \neq p$ (in this case only $\resolv_{\gamma A_1}$) and defining the result of each forward and backward step yields
\begin{align*}
	y_1 &= (\Id + \gamma^{-1} A_1^{-1})^{-1}(\gamma^{-1} z), \\
	y_2 &= A_2(z - \gamma y_1), \\
	y_3 &= (\Id + \gamma A_3)^{-1}( 2(z - \gamma y_1) - z - \gamma y_2 ), \\
	T_{A}z &= y_3 + z - (z - \gamma y_1).
\end{align*}
Rearranging the first three lines such that we only have $z$ on the left and $y_1$, $y_2$ and $y_3$ and unscaled operators on the right yields
\begin{align*}
	z &\in A_1^{-1}y_1 + \gamma y_1, \\
	z &\in A_2^{-1}y_2 + \gamma y_1, \\
	\gamma^{-1}z &\in A_3y_3 + 2y_1 + y_2 + \gamma^{-1}y_3,\\
	T_{A}z &= \gamma y_1 + y_3.
\end{align*}
If we define $y = (y_1,y_2,y_3)$, we see that the first three lines can be written as
\[
	\begin{bmatrix} 1 \\ 1 \\ \gamma^{-1} \end{bmatrix} z
	\in
	\pddiag_{A,3} y +
	\begin{bmatrix} \gamma & 0 & 0 \\ \gamma & 0 & 0 \\ 2 & 1 & \gamma^{-1} \end{bmatrix} y
	=
	\pdop_{A,3} y +
	\begin{bmatrix} \gamma & 0 & 1 \\ \gamma & 0 & 1 \\ 1 & 0 & \gamma^{-1} \end{bmatrix} y
\]
and $T_{A}$ can be written as
\begin{align*}
	y &= \left(\begin{bmatrix} \gamma & 0 & 1 \\ \gamma & 0 & 1 \\ 1 & 0 & \gamma^{-1} \end{bmatrix} + \pdop_{A,3}  \right)^{-1} \begin{bmatrix} 1 \\ 1 \\ \gamma^{-1} \end{bmatrix}z, \\
	T_{A}z &= z - \begin{bmatrix}1\end{bmatrix} z + \begin{bmatrix} \gamma & 0 & 1\end{bmatrix} y.
\end{align*}
By comparing to \cref{def:gen-pd-resolv}, we can readily identify the matrices $M$, $N$, $U$ and $V$ in the representation $(3,M,N,U,V)$.

\section{Minimal Lifting}\label{sec:minlift}

\begin{defin}[Lifting]
	The \emph{lifting number}, or \emph{lifting}, of a frugal splitting operator $T_{(\cdot)}\colon \Hprim^d \to \Hprim^d$ over $\ops_n^F$ is the number $d \in \N_+$.
\end{defin}
The lifting number represents how much memory, proportional to the problem variable in \cref{eq:prob}, is needed to store data between iterations in a fixed point iteration of the splitting operator.
For instance, if $\Hprim = \R^N$ and we are trying to find a zero associated with $A\in\ops_n^F$ with a frugal splitting operator with lifting number $3$, i.e., $T_{(\cdot)}\colon \Hprim^{3} \to \Hprim^{3}$, we need to store a variable in $\R^{3N}$ between iterations%
\ifdefined\shortversion
	{}.
\else
	\footnote{
		It is possible for the internal operations needed to evaluate the splitting operator itself to require additional memory.
		However, since this is highly problem and implementation dependent, we do not consider this.
	}.
\fi
We are interested in finding lower bounds for the lifting number as well as frugal splitting operators that attain these bounds.
\begin{defin}[Minimal Lifting]
	A frugal splitting operator $T_{(\cdot)}\colon \Hprim^d \to \Hprim^d$ over $\ops_n^F$ has \emph{minimal lifting} if $d \leq d'$ for all frugal splitting operators $T'_{(\cdot)}\colon\Hprim^{d'}\to\Hprim^{d'}$ over $\ops_n^F$.
	Furthermore, we say that $d$ is the \emph{minimal lifting} over $\ops_n^F$.
\end{defin}
The equivalent representation of a frugal splitting operator in \cref{thm:iff-representation} is useful when it comes to examining minimal lifting.
In fact, a lower bound on the lifting number is directly given by \cref{cor:factor-rep} as shown next.
\begin{cor}\label{cor:lifting-bounds}
	Let $(p,M,N,U,V)$ be a representation of a frugal splitting operator $T_{(\cdot)} \colon \Hprim^d \to \Hprim^d$ over $\ops_n^F$ that satisfies \cref{thm:iff-representation}. Then
	\begin{enum_roman}
	\item $d \geq \rank U \geq \rank N \geq \rank M$,
	\item $d \geq \rank U \geq \rank V \geq \rank M$.
	\end{enum_roman}

	\begin{proof}
		\Cref{cor:factor-rep} states that there exist matrices $S\in\R^{n\times d}$ and $P\in\R^{d\times n}$ such that $M = SUP$, $N = SU$ and $V = UP$.
		The results then follow directly from the fact that $U\in\R^{d\times d}$ and that a product of matrices cannot have greater rank than its factors.
	\end{proof}
\end{cor}
Given that the lifting number is bounded by the kernel rank, it is of interest to determine the minimal rank achievable for any valid kernel, leading us to the following definition.
\begin{defin}[Minimal Kernel Rank]
	A matrix $M\in\R^{n\times n}$ is a \emph{minimal $p$-kernel over $\ops_n^F$} if it is a $p$-kernel over $\ops_n^F$ and $\rank M \leq \rank M'$ for all $p$-kernels $M'$ over $\ops_n^F$.
	The \emph{minimal $p$-kernel rank over $\ops_n^F$} is $\rank M$ where $M$ is a minimal $p$-kernel $M$ over $\ops_n^F$.
\end{defin}
The existence of a minimal $p$-kernel over $\ops_n^F$ follows from the fact that $\rank M \in \N$ for all real matrices and a $p$-kernel over $\ops_n^F$ exists for all $n\geq 2$, $F\in\{1,\dots,n\}$, and $p\in \{1,\dots,n\}\setminus F$.
We can use the minimal $p$-kernel rank over $\ops_n^F$ and \cref{cor:lifting-bounds} to provide a lower bound on the minimal lifting number.
However, this is not enough to establish that the lower bound can be attained.
For that we need to guarantee the possibility of constructing a frugal splitting operator over $\ops_n^F$ from a $p$-kernel over $\ops_n^F$.
\begin{prop}\label{prop:kernel2split}
	Let $F \subset \{1,\dots,n\}$ and $p \in \{1,\dots,n\}\setminus F$.
	If $M$ is a $p$-kernel over $\ops_n^F$ with $\rank M = d$ and matrices $K\in\R^{n\times d}$ and $H\in\R^{d\times n}$ satisfy $\ran K = (\ker M)^\perp$ and $\ker H = (\ran M)^\perp$,
	then $(p,M,MK,HMK,HM)$ is a representation of a frugal splitting operator over $\ops_n^F$ with lifting number $d$.
	Furthermore, such matrices $K$ and $H$ exist for all $p$-kernels over $\ops_n^F$.

	\begin{proof}
		We first show that $(p,M,MK,HMK,HM)$ satisfies the conditions in \cref{thm:iff-representation}.

		\Cref{thm:iff-representation}\cref{thm:iff:kernel} is directly given by the assumption on $M$.
		Now, let $x\in\R^{d}$ and $y\in\R^n$ be such that $HMKx = HMy$.
		Since $\ker H \cap \ran M = \{0\}$, this implies $MKx = My$, which proves \cref{thm:iff-representation}\cref{thm:iff:fixedpoint-matrix}.
		For \cref{thm:iff-representation}\cref{thm:iff:zero-matrix}, let $y\in\R^n$ be arbitrary and let $y^\parallel\in\ker M$ and $y^\perp\in(\ker M)^\perp$ be such that $y = y^\parallel + y^\perp$.
		Since $\ran K = (\ker M)^\perp$, there exists $x\in\R^{d}$ such that $y^\perp = Kx$ and hence  $MKx = My^\perp = M(y^\parallel + y^\perp) = My$, multiplying both sides with $H$ from the left finally gives \cref{thm:iff-representation}\cref{thm:iff:zero-matrix}.

		Finally, let $M$ be an arbitrary $p$-kernel over $\ops_n^F$.
		Let $\rank M = d$, then the kernel $M$ has $d$ linearly independent columns and $d$ linearly independent rows.
		Define $K\in\R^{n\times d}$ such that the columns of $K$ are $d$ linearly independent rows of $M$ and define $H\in\R^{d\times n}$ such that the rows of $H$ are $d$ linearly independent columns of $M$, then $\ran K = (\ker M)^\perp$ and $\ker H = (\ran M)^\perp$.
		This concludes the proof.
	\end{proof}
\end{prop}
With these results in hand, our main minimal lifting result can be stated and proved.
\begin{thm}\label{thm:minlift}
	The minimal lifting number of a frugal splitting operator over $\ops_n^F$ is equal to the minimal $p$-kernel rank over $\ops_n^F$ for arbitrary $p \in \{1,\dots,n\}\setminus F$.

	\begin{proof}
		Let $T_{(\cdot)}\colon\Hprim^d\to\Hprim^d$ be a frugal splitting operator over $\ops_n^F$ with minimal lifting and let $p \in \{1,\dots,n\}\setminus F$.
		\Cref{thm:iff-representation} states that $T_{(\cdot)}$ has a representation $(p,M,N,U,V)$ and \cref{cor:lifting-bounds} then directly gives that the lifting $d$ is greater or equal to the minimal $p$-kernel rank over $\ops_n^F$.
		\Cref{prop:kernel2split} proves that $d$ is equal to the minimal $p$-kernel rank over $\ops_n^F$ since otherwise we could construct a frugal splitting operator with smaller lifting.
		Since the choice of $p\in\{1,\dots,n\}\setminus F$ was arbitrary and $d$ is independent of $p$, the minimal $p$-kernel rank over $\ops_n^F$ is the same for all $p\in\{1,\dots,n\}\setminus F$.
	\end{proof}
\end{thm}
This theorem reduces the problem of finding the minimal lifting over $\ops_n^F$ to finding the minimal $p$-kernel rank over $\ops_n^F$.
Furthermore, since \cref{prop:kernel2split} proves the existence results necessary for \cref{thm:minlift} by construction, it provides a constructive approach to the design of frugal splitting operators with lifting equal to the kernel rank.
Also the construction of frugal splitting operators with minimal lifting can be seen as a problem of finding $p$-kernels over $\ops_n^F$ with minimal rank.
We will give an example of this in \cref{sec:new-method} and end this section by making \cref{thm:minlift} concrete by finding the minimal kernel rank over $\ops_n^F$ and the corresponding minimal lifting results.
\begin{cor}\label{cor:fo-minlift-several}
	Let $n\geq 2$ and $F \subset \{1,\dots,n\}$.
	The minimal lifting over $\ops_n^F$ is $n-|F|$ if $1 \in F$ or $n \in F$, otherwise it is $n-1-|F|$.

	\begin{proof}
		Let $M$ be a $p$-kernel over $\ops_n^F$.
		It has the structure
		\begin{align*}
			M
			=
			\begin{fmatrix}{c|c|c}
				L_{1} & \mathbf{1} & \mathbf{0} \\
				\hline
				\ast & l_p & -\mathbf{1} \\
				\hline
				\ast & \ast & L_{2} \\
			\end{fmatrix}
			\in \R^{n\times n},
		\end{align*}
		where $L_1 \in \R^{(p-1)\times(p-1)}$ and $L_2 \in \R^{(n-p)\times(n-p)}$ are lower triangular matrices and $l_p > 0$ is a real number.
		The symbols $\ast$,	$\mathbf{1}$ and $\mathbf{0}$ respectively denote an arbitrary real matrix, a matrix of all ones, and a matrix of all zeros, all of appropriate sizes for their position.
		Since reordering columns and rows does not change the rank of a matrix, there exist matrices
		\begin{align*}
			M_c
			=
			\begin{fmatrix}{c|c|c}
				L_{1} & \mathbf{0} & \mathbf{1} \\
				\hline
				\ast & L_{2} & \ast \\
				\hline
				\ast & -\mathbf{1} & l_p\\
			\end{fmatrix}
			\in \R^{n\times n}
			\quad\text{and}\quad
			M_r
			=
			\begin{fmatrix}{c|c|c}
				l_p & \ast & -\mathbf{1} \\
				\hline
				\mathbf{1} & L_{1} & \mathbf{0} \\
				\hline
				\ast & \ast & L_{2} \\
			\end{fmatrix}
			\in \R^{n\times n}
		\end{align*}
		such that $\rank M = \rank M_c = \rank M_r$.

		Consider the first $n-1$ columns of $M_c$.
		Since $L_1$ and $L_2$ are lower triangular, the number of linearly independent columns is greater or equal to the number of non-zero diagonal elements of $L_1$ and $L_2$.
		Furthermore, since the diagonal elements of $L_1$ and $L_2$ are diagonal elements of $M$, $M_{i,i} = 0$ if and only if $i \in F$, and $p\notin F$, there are at least $n-1-|F|$ linearly independent columns among the first $n-1$ columns of $M_c$, hence, $\rank M_c \geq n-1-|F|$.
		If $1\notin F$ then $M_{1,1} = (L_1)_{1,1} \neq 0$ and this bound can be attained by selecting the first column of $M_c$ parallel to the last and putting zeros in all other positions that are allowed to be zero.
		If $1 \in F$ then $M_{1,1} = (L_1)_{1,1} = 0$ and the last column of $M_c$ is not in the span of the first $n-1$ columns and hence $\rank M_c \geq n-|F|$.
		This bound is attained by putting zeros in all positions of $M_c$ that are allowed to be zero.
		By considering rows of $M_r$ instead of columns of $M_c$ we can analogously conclude that $\rank M_r \geq n-1-|F|$ and if $n \in F$ (and hence $M_{n,n} = (L_2)_{n-p,n-p} = 0$), then  $\rank M_r \geq n-|F|$.
		These bounds are also similarly attained.

		Since $\rank M = \rank M_c = \rank M_r$ we have that $\rank M \geq n-1-|F|$ and if $1\in F$ or $n\in F$ then $\rank M \geq n-|F|$.
		Choices of $M$ that attain these bounds can be constructed by reordering the rows and columns of the choices of $M_c$ and $M_r$ that attain their respective bounds.
		The minimal kernel rank over $\ops_n^F$ is then $n-|F|$ if $1 \in F$ or $n\in F$, otherwise it is $n-1-|F|$.
		The lifting result then follows from \cref{thm:minlift}.
	\end{proof}
\end{cor}
In the case of frugal resolvent splitting operators, i.e., $F = \emptyset$, this is the same lower bound as the one found by Malitsky and Tam \cite{malitskyResolventSplittingSums2021}.
When considering frugal splitting operators with forward evaluations, i.e., $F \neq \emptyset$, this corollary uncovers an interesting phenomenon where the minimal lifting depends on the evaluation order.
A smaller lifting number is possible as long as neither the first nor the last evaluation is a forward evaluation.
This makes it clear why the three operator splitting of Davis and Yin, see \cref{sec:rep-ex}, achieves a lifting number of one while the corresponding methods of V\~u and Condat \cite{vuSplittingAlgorithmDual2013,condatPrimalDualSplitting2013} require a lifting of two.
Davis and Yin's method performs the forward evaluation between two resolvent evaluations, while the V\~u--Condat method performs the forward evaluation first.

\section{Convergence}\label{sec:conv}
In this section, we show convergence of fixed point iterations $z_{k+1} = T_{A}z_k$ for frugal splitting operators $T_{(\cdot)}\colon\Hprim^d \to \Hprim^d$ over $\ops_n^F$ and $A\in\ops_n^F$.
From \cref{cor:factor-rep} we know that these iterations can, for each $p \in\{1,\dots,n\}\setminus F$, be written as
\begin{equation}\label{eq:fp-iteration-rep}
	\left\{
		\begin{aligned}
			y_k &= (SUP + \pdop_{A,p})^{-1}SUz_k, \\
			z_{k+1} &= z_k - U(z_k - Py_k)
		\end{aligned}
	\right.
	\quad\text{or equivalently}\quad
	\left\{
		\begin{aligned}
			SU(z_k - Py_k) &\in \pdop_{A,p}y_k, \\
			U(z_k - Py_k) &= z_k - z_{k+1}
		\end{aligned}
	\right.
\end{equation}
for some matrices $S\in\R^{n\times d}$, $P\in\R^{d\times n}$ and $U\in\R^{d\times d}$.
In certain cases, this iteration can be analyzed with existing theory.
For instance, if $S$ is symmetric positive definite, $P = I$ and $U = \theta I$ for some $\theta \in (0,1]$, then this is an averaged fixed point iteration of a resolvent in the metric given by the symmetric kernel $SUP = \theta S$.
For general frugal splitting operators, $SUP$ may be non-symmetric, invalidating application of existing convergence results.
Therefore, we will derive sufficient convergence conditions, without symmetry requirements on the kernel, under the following assumption.
\begin{ass}\label{ass:coco}
	Let $A = (A_1,\dots,A_n)\in\ops_n^F$ be such that the following hold
	\begin{enum_roman}
	\item $\zer \sum_{i=1}^n A_i \neq \emptyset$,
	\item $A_i$ is $\beta_i$-cocoercive for all $i \in F$.
	\end{enum_roman}
	For this $A$, let $B \in \R^{n\times n}$ and $B^\dagger \in \R^{n\times n}$ be such that $B_{i,i} = \beta_i$ and $B^\dagger_{i,i} = \beta_i^{-1}$ for all $i \in F$, with all other elements of the matrices equal to zero.
\end{ass}
\ifdefined\shortversion
\else
	Cocoercivity of the forward evaluated operators is a standard setting for proving convergence of forward-backward like methods.
	However, there exist frugal splitting operators whose fixed point iterations converge under only Lipschitz assumptions on the operator used for forward evaluations.
	An example of this, that does not have minimal lifting, is the forward-reflected-backward splitting \cite{malitskyForwardBackwardSplittingMethod2020}.

\fi
Since inverses of cocoercive operators are strongly monotone, primal-dual operators for tuples that satisfy \cref{ass:coco} have the following strong-monotonicity-like property.
\begin{lem}\label{lem:pdop-smon}
	Let $F \subset \{1,\dots,n\}$, $p \in \{1,\dots,n\}\setminus F$ and let $A \in \ops_n^F$ satisfy \cref{ass:coco}, then
	\[
		\inprod{u - y}{x-y} \geq \norm{x-y}_{B}^2
	\]
	for all $x,y \in \Hprim^n$ and for all $u \in \pdop_{A,p}x$, $v \in \pdop_{A,p}y$.
	\begin{proof}
		With some abuse of notation for set-valued operators, we first note that for all $x,y\in\Hprim^n$, we have
		\begin{align*}
			&\inprod{\pdop_{A,p}x - \pdop_{A,p}y}{x-y} \\
			&\qquad= \inprod{\pddiag_{A,p}x - \pddiag_{A,p}y}{x-y} + \inprod{\pdskew_{p}x - \pdskew_{p}y}{x-y} \\
			&\qquad= \inprod{\pddiag_{A,p}x - \pddiag_{A,p}y}{x-y} \\
			&\qquad= \inprod{A_px_p - A_py_p}{x_p-y_p} + \sum_{i\in\{1,\dots,n\}\setminus \{p\}} \inprod{A_i^{-1}x_i - A_i^{-1}y_i}{x_i-y_i},
		\end{align*}
		where the fact that $\pdskew_p$ is skew-adjoint is used in the second equality.
		For all $i \in \{1,\dots,n\}\setminus F$, the operator $A_i$ is monotone and so is $A_i^{-1}$, while for all $i \in F$, $A_i$ is $\beta_i$-cocoercive and hence $A_i^{-1}$ is $\beta_i$-strongly monotone.
		Noting that $p \notin F$ by assumption and using these properties yields
		\[
			\inprod{\pdop_{A,p}x - \pdop_{A,p}y}{x-y}
			\geq \sum_{f\in F} \beta_f\norm{x_f-y_f}^2
			= \norm{x-y}_{B}^2.
			\qedhere
		\]
	\end{proof}
\end{lem}
This property is used to derive our main convergence theorem.
\begin{thm}\label{thm:conv}
	Let $F \subset \{1,\dots,n\}$ and let $A \in \ops_n^F$ satisfy \cref{ass:coco}.
	Let $T_{(\cdot)}\colon\Hprim^d\to\Hprim^d$ be a frugal splitting operator over $\ops_n^F$ and let $(p,SUP,SU,U,UP)$ be a representation of $T_{(\cdot)}$ according to \cref{cor:factor-rep}.
	Let the sequences $\{z_k\}_{k\in\N}$ and $\{y_k\}_{k\in\N}$ be generated by \cref{eq:fp-iteration-rep} for this representation and some $z_0 \in\Hprim^d$.

	Consider the following conditions on a symmetric matrix $Q\in\R^{d\times d}$,
	\begin{enum_alpha}
		\item\label{thm:conv:struct-cond} $(I - I_F)(P^T Q - S)U = 0$,
		\item\label{thm:conv:seq-cond} $Q\succ 0$ and $W \succ 0$,
	\end{enum_alpha}
	where $W = QU + (QU)^T - U^TQU - \tfrac{1}{2}(P^TQU - SU)^TB^\dagger(P^TQU - SU)$, $I \in \R^{n\times n}$ is the identity matrix, and $I_F \in \R^{n\times n}$ is the diagonal matrix with $(I_F)_{i,i} = 1$ if $i\in F$, otherwise $(I_F)_{i,i} = 0$.
	If \cref{thm:conv:struct-cond} holds, then
	\[
		\norm{z_{k+1} - P\bar{y}}_Q^2 \leq \norm{z_{k} - P\bar{y}}_Q^2 - \norm{z_k - Py_k}_W^2
	\]
	for all $k\in\N$ and all $\bar y \in \zer \pdop_{A,p}$.
	If \cref{thm:conv:seq-cond} also holds, then
	\begin{enum_roman}
	\item\label{thm:conv:res} $z_k - Py_k \to 0$,
	\item\label{thm:conv:zero} $\pdop_{A,p}y_k \ni SU(z_k - Py_k) \to 0$,
	\item\label{thm:conv:sol} $y_k \weakto y^\star$,
	\item\label{thm:conv:iters} $z_k \weakto Py^\star \in \fix T_{A}$,
	\end{enum_roman}
	for some $y^\star \in \zer \pdop_{A,p}$ as $k\to \infty$.

	\begin{proof}
		Let $Q \in\R^{d\times d}$ be a symmetric matrix and let $k \in \N$.
		For arbitrary $\bar{y}\in\zer\pdop_{A,p}$, we have
		\begin{align*}
			\norm{z_{k+1} - P\bar{y}}_Q^2
			&= \norm{z_{k} - P\bar{y} - U(z_k - Py_k)}_Q^2 \\
			&= \norm{z_{k} - P\bar{y}}_Q^2 + \norm{z_k - Py_k}_{U^TQU}^2 - 2\inprod{QU(z_k - Py_k)}{z_k - P\bar{y}}.
		\end{align*}
		Insert $Py_k - Py_k$ on the right hand side of the inner product to get
		\begin{align*}
			\norm{z_{k+1} - P\bar{y}}_Q^2
			&= \norm{z_{k} - P\bar{y}}_Q^2 + \norm{z_k - Py_k}_{U^TQU}^2
			- 2\inprod{QU(z_k - Py_k)}{z_k - Py_k} \\
			&\quad - 2\inprod{QU(z_k - Py_k)}{Py_k - P\bar{y}} \\
			&= \norm{z_{k} - P\bar{y}}_Q^2 - \norm{z_k - Py_k}_{QU + (QU)^T - U^TQU}^2 \\
			&\quad - 2\inprod{P^TQU(z_k - Py_k)}{y_k - \bar{y}}.
		\end{align*}
		From \cref{eq:fp-iteration-rep} we know that $SU(z_k - Py_k) \in \pdop_{A,p}y_k$ and \cref{lem:pdop-smon} then gives that
		\[
			0 \leq \inprod{SU(z_k - Py_k)}{y_k - \bar{y}} - \norm{y_k - \bar{y}}_B^2.
		\]
		Add two times this inequality to the previous equality to get
		\begin{align*}
			\norm{z_{k+1} - P\bar{y}}_Q^2
			&\leq \norm{z_{k} - P\bar{y}}_Q^2 - \norm{z_k - Py_k}_{QU + (QU)^T - U^TQU}^2 - 2\norm{y_k - \bar{y}}_B^2\\
			&\quad - 2\inprod{P^TQU(z_k - Py_k)}{y_k - \bar{y}} + 2\inprod{SU(z_k - Py_k)}{y_k - \bar{y}} \\
			&= \norm{z_{k} - P\bar{y}}_Q^2 - \norm{z_k - Py_k}_{QU + (QU)^T - U^TQU}^2 - 2\norm{y_k - \bar{y}}_B^2\\
			&\quad - 2\inprod{(P^TQU - SU)(z_k - Py_k)}{y_k - \bar{y}}.
		\end{align*}
		Assume that condition \cref{thm:conv:struct-cond} is satisfied, then $P^TQU - SU = I_F(P^TQU - SU)$ and since $I_F = (2B)^{1/2}(\tfrac{1}{2}B^\dagger)^{1/2}$, we have
		\begin{align*}
			\norm{z_{k+1} - P\bar{y}}_Q^2
			&\leq \norm{z_{k} - P\bar{y}}_Q^2 - \norm{z_k - Py_k}_{QU + (QU)^T - U^TQU}^2 - 2\norm{y_k - \bar{y}}_B^2\\
			&\quad - 2\inprod{(2^{-1}B^\dagger)^{1/2}(P^TQU - SU)(z_k - Py_k)}{(2B)^{1/2}(y_k - \bar{y})}.
		\end{align*}
		Using Young's inequality results in
		\begin{equation}\label{eq:lyapineq}
			\begin{aligned}
				\norm{z_{k+1} - P\bar{y}}_Q^2
				&\leq \norm{z_{k} - P\bar{y}}_Q^2 - \norm{z_k - Py_k}_{QU + (QU)^T - U^TQU}^2 \\
				&\quad + \norm{(2^{-1}B^\dagger)^{1/2}(P^TQU - SU)(z_k - Py_k)}^2 \\
				&= \norm{z_{k} - P\bar{y}}_Q^2 - \norm{z_k - Py_k}_{W}^2,
			\end{aligned}
		\end{equation}
		where $W = QU + (QU)^T - U^TQU - \tfrac{1}{2}(P^TQU - SU)^TB^\dagger(P^TQU - SU)$. This proves the first statement.

		Let us now assume that \cref{thm:conv:seq-cond} also holds. This implies that both $\norm{\cdot}_Q$ and $\norm{\cdot}_W$ are norms on $\Hprim^d$.
		Adding the inequality in \cref{eq:lyapineq} for $k=0,1,\dots$ implies that $\{\norm{z_k - Py_k}_W^2\}_{k\in\N}$ is summable and that
		\[
			\norm{z_k - Py_k}_W^2 \to 0,
		\]
		implying that  \cref{thm:conv:res} holds.
		Statement \cref{thm:conv:zero} follows directly from \cref{thm:conv:res} and \cref{eq:fp-iteration-rep}.
		From the inequality in \cref{eq:lyapineq}, we also conclude that
		\[
			\norm{z_{k+1} - P\bar{y}}_Q^2 \leq \norm{z_{k} - P\bar{y}}_Q^2
		\]
		for all $k\in\N$ and all $\bar{y}\in\zer\pdop_{A,p}$, implying the boundedness of $\{z_k\}_{k\in\N}$ and the convergence of $\{\norm{z_{k}-P\bar{y}}_Q^2\}_{k\in\N}$ for all $\bar{y}\in\zer\pdop_{A,p}$.
		These two facts along with \cref{thm:conv:res,thm:conv:zero} will be used to prove \cref{thm:conv:sol}.
		Statement \cref{thm:conv:iters} then follows directly from the weak continuity of $P$ and \cref{thm:conv:sol}.

		For the proof of \cref{thm:conv:sol}, we first show that $\{y_k\}_{k\in\N}$ is bounded and hence has weak sequential cluster points.
		We move on to show that these cluster points are in $\zer \pdop_{A,p}$ and last we establish that there exists at most one cluster point.
		The convergence $y_k \weakto y^\star$ for some $y^\star \in \zer \pdop_{A,p}$ then follows from \cite[Lemma 2.46]{bauschkeConvexAnalysisMonotone2017}.

		Since $SUP$ is a $p$-kernel over $\ops_n^F$, the operator $(SUP + \pdop_{A,p})^{-1}\circ SU$ can be evaluated using a finite number of vector additions, scalar multiplications, resolvents and evaluations of the cocoercive operators $A_i$ with $i \in F$.
		All of these operations are Lipschitz continuous and therefore $(SUP + \pdop_{A,p})^{-1}\circ SU$ is Lipschitz continuous, let us say with constant $L$.
		Furthermore, with $\bar{y} \in \zer \pdop_{A,p}$ we have $\bar{y} = (SUP+\pdop_{A,p})^{-1}SUP\bar{y}$ and
		\[
			\norm{y_k - \bar{y}} = \norm{(SUP+\pdop_{A,p})^{-1}SUz_k - (SUP+\pdop_{A,p})^{-1}SUP\bar{y}} \leq L\norm{z_k - P\bar{y}}.
		\]
		Since the sequence $\{z_k\}_{k\in\N}$ is bounded, so is $\{y_k\}_{k\in\N}$, implying it has convergent subsequences \cite[Lemma 2.45]{bauschkeConvexAnalysisMonotone2017}.
		Using \cref{thm:conv:zero} the maximal monotonicity of $\pdop_{A,p}$, see \cref{sec:problem} and the weak-strong closedness of graphs of maximally monotone operators \cite[Proposition 20.38]{bauschkeConvexAnalysisMonotone2017}, we conclude that all weak sequential cluster points of $\{y_k\}_{k\in\N}$ belong to $\zer \pdop_{A,p}$.

		What remains is to show that $\{y_k\}_{k\in\N}$ possesses at most one weak sequential cluster point.
		Let $\{y_{k_i}\}_{i\in\N}$ and $\{y_{k'_j}\}_{j\in\N}$ be sub-sequences such that $y_{k_i} \weakto a^\star \in\zer\pdop_{A,p}$ and $y_{k'_j} \weakto b^\star\in\zer\pdop_{A,p}$.
		Since $P$ as an operator on $\Hprim^n$ is linear and hence weakly continuous, \cref{thm:conv:res} implies that
		\[
			z_{k_i} \weakto Pa^\star \quad\text{and}\quad z_{k'_j} \weakto Pb^\star.
		\]
		Furthermore, $\{\norm{z_k - Pa^\star}_Q^2\}_{k\in\N}$ and $\{\norm{z_k - Pb^\star}_Q^2\}_{k\in\N}$ both converge and hence
		\[
			\inprod{z_k}{Pa^\star-Pb^\star}_Q
			= \tfrac{1}{2}\norm{z_k - Pb^\star}_Q^2 - \tfrac{1}{2}\norm{z_k - Pa^\star}_Q^2 - \tfrac{1}{2}\norm{Pb^\star}_Q^2 + \tfrac{1}{2}\norm{Pa^\star}_Q^2 \to \nu
		\]
		for some $\nu \in \R$.
		In particular, it means that
		\begin{align*}
			\inprod{z_{k_i}}{Pa^\star-Pb^\star}_Q \rightarrow \inprod{Pa^\star}{Pa^\star-Pb^\star}_Q = \nu
		\end{align*}
                and
		\begin{align*}
                  \inprod{z_{k'_j}}{Pa^\star-Pb^\star}_Q \rightarrow  \inprod{Pb^\star}{Pa^\star-Pb^\star}_Q  =\nu,
		\end{align*}
                leading to
		\[
			0 = \inprod{Pa^\star-Pb^\star}{Pa^\star-Pb^\star}_Q = \norm{Pa^\star-Pb^\star}_{Q}^2
		\]
		and hence $Pa^\star = Pb^\star$.
		Since $y \in \zer \pdop_{A,p}$ implies $y = (SUP + \pdop_{A,p})^{-1}SUPy$ we have
		\begin{align*}
			a^\star
			= (SUP + \pdop_{A,p})^{-1}SUPa^\star
			= (SUP + \pdop_{A,p})^{-1}SUPb^\star
			= b^\star.
		\end{align*}
		This shows that at most one weak sequential cluster point exists.
		Hence, $\{y_k\}_{k\in\N}$ converges weakly to some $y^\star \in \zer \pdop_{A,p}$. This concludes the proof.
	\end{proof}
\end{thm}
\ifdefined\shortversion
	A consequence of \cref{thm:conv} is that the sequence generated by a fixed point iteration of a frugal splitting operator that satisfies the conditions of \cref{thm:conv} is Fej\'er monotone with respect to $P\zer \pdop_{A,p}$ in the Hilbert space given by $\inprod{\cdot}{\cdot}_Q$.
	This holds true even for methods without minimal lifting such as momentum methods and we provide several examples of this in the extended preprint \cite{morinFrugalSplittingOperators2022arxiv}, where we also show that \cref{thm:conv} recovers convergence results for many well-known splitting methods.
	We provide, for instance, representations and convergence results for
	forward-backward splitting \cite{goldsteinConvexProgrammingHilbert1964,levitinConstrainedMinimizationMethods1966},
	Douglas--Rachford splitting \cite{lionsSplittingAlgorithmsSum1979},
	the three operator splitting of Davis and Yin \cite{davisThreeOperatorSplittingScheme2017},
	and the minimal lifting methods of Ryu \cite{ryuUniquenessDRSOperator2020} and Malitsky--Tam \cite{malitskyResolventSplittingSums2021}.
\else
	A few remarks on this theorem are in order.
	The sequence generated by a fixed point iteration of a frugal splitting operator that satisfies \cref{thm:conv} is Fej\'er monotone, see \cite[Definition 5.1]{bauschkeConvexAnalysisMonotone2017} with respect to $P\zer \pdop_{A,p}$ in the Hilbert space given by $\inprod{\cdot}{\cdot}_Q$.
	This is true even for methods without minimal lifting such as momentum methods, see \cref{sec:conv-examples}.
	We will in particular show that forward-backward splitting with Nesterov-like momentum where the momentum parameter is fixed satisfies \cref{thm:conv} and hence is Fej\'er monotone.
\fi

Condition \cref{thm:conv:struct-cond} becomes stricter when the set $F$ of indices on which we make forward evaluations becomes smaller.
This puts stronger restrictions on the structure of the frugal splitting operator.
For instance, assuming that $S$, $U$ and $P$ are square and invertible and $F = \emptyset$, condition \cref{thm:conv:struct-cond} states that $Q = (P^T)^{-1}S$ and hence $(P^T)^{-1}S$ must be symmetric.
Other cases leaves more freedom in choosing $Q \succ 0$ such that condition \cref{thm:conv:seq-cond} of \cref{thm:conv} is satisfied.

Although $W$ is given by a quadratic expression in $Q$, the requirement $W \succ 0$ in \cref{thm:conv} can be transformed to an equivalent positive definite condition that is linear in $Q$ by using a Schur complement. We have
\begin{align*}
	W \succ 0
	\iff
	\begin{bmatrix}
		QU + (QU)^T - U^TQU & U^T(P^TQ - S)^T (\tfrac{1}{2}B^\dagger)^{1/2}\\
		(\tfrac{1}{2}B^\dagger)^{1/2}(P^TQ - S)U & I
	\end{bmatrix}
	\succ 0
\end{align*}
where $I\in\R^{n\times n}$ is an identity matrix and $(\tfrac{1}{2}B^\dagger)^{1/2}$ exists since $B^\dagger$ is a diagonal matrix with non-negative elements.
This makes the search for a matrix $Q$ that satisfies \cref{thm:conv:struct-cond,thm:conv:seq-cond} a convex semi-definite feasibility problem that can be readily solved numerically.

From the $W \succ 0$ of condition \cref{thm:conv:seq-cond} in \cref{thm:conv}, we conclude that $U$ must have full rank. If not, there exists a non-zero element $x\in \R^d$ such that $Ux = 0$ which implies $x^TWx = 0$ which is a contradiction.
This is convenient since, given a frugal splitting operator with representation $(p,M,N,U,V)$ where $U$ is invertible, it is easy to find a representation of the form used in \cref{thm:conv}, i.e., $(p,SUP,SU,U,UP)$.
In fact, this factorization can be expressed only in terms of the matrices $N$, $U$ and $V$,
\[
	S = NU^{-1} \quad\text{and}\quad P = U^{-1}V,
\]
which implies that $M = NU^{-1}V$ must hold for such a frugal splitting operator.

There are frugal splitting operators where $U$ is rank deficient but their convergence is of no interest.
This is because they are either guaranteed to not converge in general or they can be reduced to a frugal splitting operator, without losing any information, whose representation has a full rank $U$.
To see this, let $(p,SUP,SU,U,UP)$ be a representation satisfying \cref{cor:factor-rep} of a frugal splitting operator $T_{(\cdot)}\colon \Hprim^d \to \Hprim^d$.
Let $\proj\in \R^{d\times d}$ be the projection matrix onto $\ran U^T$ where the projection is in the standard Euclidean $\R^d$ space.
The projection matrix on the kernel of $U$ is then $\overline{\proj} = I - \proj$.
Consider the fixed point iteration $z_{k+1} = T_Az_k$ for some $A \in\ops_n^F$ and $z_0 \in \Hprim^d$.
If we define $z^\parallel_0 = \proj z_0$ and $z^\perp_0 = \overline{\proj}z_0$ this fixed point iteration can be written as $z_{k} = z^\parallel_{k} + z^\perp_{k}$ where
\begin{align*}
	y_k &= (SUP + \pdop_{A,p})^{-1}SU z^\parallel_k, \\
	z^\parallel_{k+1} &= z^\parallel_k - \proj U(z^\parallel_k - Py_k), \\
	z^\perp_{k+1} &= z^\perp_k - \overline{\proj}U(z^\parallel_k - Py_k),
\end{align*}
for all $k\in\N$.
If $\ran U = \ran U^T$, then \(\overline{\proj} U = 0\) and $z^\perp_k = z^\perp_0$ for all $k\in\N$, and only the sequence $\{z^\parallel_{k}\}_{k\in\N}$ is of interest.
Furthermore, the sequence $\{z^\parallel_{k}\}_{k\in\N}$ can be recovered from a fixed point iteration of a frugal splitting operator with lifting equal to the rank of $U$, since it lives in the space \(\ran U\) of this exact dimension.
If instead $\ran U \neq \ran U^T$, it is always possible to find an operator tuple $A \in \ops_n^F$ and initial point $z_0\in\Hprim^d$ such that $z^\parallel_k = z_0$, $z^\perp_k = -kc$ and $z_k = z_0 - kc$ for all $k \in \N$ and some $c \in \Hprim^d \setminus \{0\}$.
Hence, $\{z_k\}_{k\in\N}$ will always diverge for this choice of $A$ and $z_0$.

\ifdefined\shortversion
\else
	\subsection{Applications of \texorpdfstring{\cref{thm:conv}}{the Convergence Theorem}}\label{sec:conv-examples}
	\Cref{thm:conv} recovers many well-known convergence results, for instance the results for forward-backward splitting \cite{goldsteinConvexProgrammingHilbert1964,levitinConstrainedMinimizationMethods1966}, Douglas--Rachford splitting \cite{lionsSplittingAlgorithmsSum1979}, the Chambolle--Pock method \cite{chambolleFirstOrderPrimalDualAlgorithm2011}, and the minimal lifting methods of Ryu \cite{ryuUniquenessDRSOperator2020} and Malitsky--Tam \cite{malitskyResolventSplittingSums2021}.
	We will not present all these results here and settle for presenting the result for the three operator splitting of Davis and Yin \cite{davisThreeOperatorSplittingScheme2017}.
	We will also present convergence conditions for the fixed point iteration of the forward-backward operator with Nesterov-like momentum \cite{nesterovMethodSolvingConvex1983,beckFastIterativeShrinkageThresholding2009}.
	To save space, we will not derive any of the representations and simply state the primal index $p$ and the matrices $U$, $S$ and $P$ of a representation $(p,SUP,SU,U,UP)$, see \cref{cor:factor-rep}.
	Similarly, for the convergence results we just state the matrices $Q$ and $W$ of \cref{thm:conv}.
	Detailed derivations of all the listed examples and more can be found in the supplement.

	\paragraph{Three Operator Splitting of Davis--Yin}
	The three operator splitting of Davis and Yin has already been presented in \cref{sec:rep-ex} but we restate it here,
	\[
		x_{k+1} = x_{k} - \resolv_{\gamma A_1}x_k + \resolv_{\gamma A_3} \circ (2\resolv_{\gamma A_1} - \Id  - \gamma A_2 \circ \resolv_{\gamma A_1})x_k
	\]
	where $x_0 \in \Hprim$, $\gamma > 0$ and $A = (A_1,A_2,A_3) \in \ops_3^{\{2\}}$.
	The representation derived in that section can be factored into a representation of the form $(3,SUP,SU,U,UP)$ where
	\[
		U = \begin{bmatrix} 1 \end{bmatrix}
		,\quad
		S = \begin{bmatrix} 1 \\ 1 \\ \gamma^{-1} \end{bmatrix}
		\quad\text{and}\quad
		P = \begin{bmatrix} \gamma & 0 & 1 \end{bmatrix}.
	\]
	To prove convergence, we choose $Q = \begin{bmatrix} \gamma^{-1}\end{bmatrix}$ in \cref{thm:conv} and, assuming $A$ satisfies \cref{ass:coco}, this results in $0 = (I-I_F)(P^TQ -S)U$ and $W = \begin{bmatrix}\gamma^{-1} - \tfrac{1}{2\beta_2}\end{bmatrix}$.
	If $\gamma < 2\beta_2$, then $x_{k} \weakto Py^\star$ for $y^\star \in \zer \pdop_{A,3}$ and
	\[
		\resolv_{\gamma A_3} \circ (2\resolv_{\gamma A_1} - \Id  - \gamma A_2 \circ \resolv_{\gamma A_1})x_k = \resolv_{\gamma A_1}x_k + x_{k+1} - x_k \weakto x^\star \in \zer A_1 + A_2 + A_3
	\]
	and hence also $\resolv_{\gamma A_1}x_k \weakto x^\star$.
	This is the same convergence condition as the one presented in \cite{davisThreeOperatorSplittingScheme2017}.

	\paragraph{Forward-Backward with Nesterov-like Momentum}
	A fixed point iteration of a forward-backward operator with Nesterov-like momentum \cite{nesterovMethodSolvingConvex1983} can be written as
	\[
		x_{k+1} = \resolv_{\gamma A_2}(x_k + \theta(x_k - x_{k-1}) - \gamma A_1(x_k + \theta(x_k - x_{k-1})))
	\]
	for $\lambda > 0$, $\theta \in \R$, $x_0,x_{-1}\in\Hprim$ and $A=(A_1,A_2)\in\ops_{2}^{\{1\}}$.
	In the proximal-gradient setting, this is also the same update that is used in the FISTA method \cite{beckFastIterativeShrinkageThresholding2009}.
	Nesterov momentum gained popularity due to it achieving optimal convergence rates in the smooth optimization setting.
	However, these faster convergence rates require a momentum parameter $\theta$ that varies between iterations, something the fixed point iterations considered in this paper will not allow.

	We remove the dependency on previous iterations by introducing an extra iterate as
	\begin{align*}
		x_{k+1} &= \resolv_{\gamma A_2}(x_k + \theta y_k - \gamma A_1(x_k + \theta y_k)), \\
		y_{k+1} &= x_{k+1} - x_k
	\end{align*}
	where $x_0,y_0\in\Hprim$.
	This is a fixed point iteration of the frugal splitting operator given by the representation $(2,SUP,SU,U,UP)$ where
	\[
		U = \begin{bmatrix} 1 & 0 \\ 1 & 1 \end{bmatrix}
		,\quad
		S = \begin{bmatrix} 1-\theta & \theta \\ \gamma^{-1}(1-\theta) & \gamma^{-1}\theta \end{bmatrix}
		\quad\text{and}\quad
		P = \begin{bmatrix} 0 & 1 \\ 0 & 0 \end{bmatrix}.
	\]
	The lifting number is two, which is not minimal, even though the kernel $SUP$ has rank one and is minimal.
	In fact, the kernel is the same as for the ordinary forward-backward splitting.
	\Cref{thm:conv} with
	\[
		Q = \gamma^{-1}\begin{bmatrix} 1 - \theta & \theta \\ \theta & a \end{bmatrix}
		\quad\text{and}\quad
		W = \gamma^{-1}\begin{bmatrix} 1-\theta-\hat{\gamma}-a & \theta(1-\hat{\gamma}) \\ \theta(1-\hat{\gamma}) & a-\theta^2\hat{\gamma} \end{bmatrix}
	\]
	where $a > 0$ and $\hat{\gamma} = \tfrac{\gamma}{2\beta_1} > 0$ yields the convergence of $x_k \weakto x^\star \in \zer A_1 + A_2$ and $y_k \weakto 0$ if $A$ satisfies \cref{ass:coco} and both $Q \succ 0$ and $W \succ 0$.
	It can be verified that $Q-W \succeq 0$ hence it is enough that $W\succ 0$ which is equivalent to
	\[
		0 < a - \theta^2\hat{\gamma}
		\quad\text{and}\quad
		0 < 1 - \theta - \hat{\gamma} - a - \tfrac{\theta^2(1- \hat{\gamma})^2}{a - \theta^2\hat{\gamma}}.
	\]
	If we restrict these results to $\theta > 0$, these conditions hold with $a = \theta^{2}\hat{\gamma} + \theta(1-\hat{\gamma})$ if
	\[
		0 < 1 - 3\theta - \tfrac{\gamma}{2\beta_1}(\theta - 1)^2.
	\]
	It is easily verified that this condition also implies the existence of an $a > 0$ that makes $Q \succ 0$ and $W \succ 0$ even in the case when $\theta = 0$.
	In fact, it reduces to the well known result $0 < \gamma < 2\beta_1$ for ordinary forward-backward splitting.
	As mentioned before, if these conditions hold then \cref{thm:conv} gives Fej\'er monotonicity of $\{(x_k, y_k)\}_{k\in\N}$ w.r.t. $P\zer \pdop_{A,p} = \zer(A_1 + A_2) \times \{0\}$ in the norm $\norm{\cdot}_Q$.
\fi

\section{A New Frugal Splitting Operator With Minimal Lifting}\label{sec:new-method}
In this section, we derive a new frugal splitting operator with minimal lifting.
The approach we take is the one outlined in \cref{prop:kernel2split}, i.e., we will select a primal index $p$, a kernel $M$, matrices $H$ and $K$, and form a representation as $(p,M,MK,HMK,HM)$.
As long as the matrices satisfy \cref{prop:kernel2split}, the resulting generalized primal-dual resolvent is a frugal splitting operator with lifting equal to the rank of the kernel.
Therefore, if we choose a kernel with minimal rank, the resulting frugal splitting operator will have minimal lifting.

As noted in \cref{sec:minlift}, the minimal lifting number of a frugal splitting operator over $\ops_n^F$ depends on $F$.
\cref{cor:fo-minlift-several} reveals that it is not only a question regarding whether $F$ is the empty set or not; the minimal lifting number depends on the actual order of the forward and backward evaluations.
For this reason, we choose to construct a frugal splitting operator over $\ops_n^F$ where $F$ is such that $1 \notin F$ and $n \notin F$.
\Cref{cor:fo-minlift-several} then guarantees the minimal lifting to be $n-1-|F|$ instead of potentially being $n-|F|$.
To further simplify the setup we assume that all the single-valued evaluations come one after another as $F = \{n-f,\dots,n-1\}$ where $f = |F|$ is the number of single-valued operators.
\begin{thm}\label{thm:new-alg}
	Let $n,f \in \N$ be such that $n \geq 2$ and $f \leq n-2$.
	Let $A = (A_1,\dots,A_n) \in \ops_n^F$ where $F = \{n-f,\dots,n-1\}$ if $f > 0$ and $F = \emptyset$ if $f = 0$.
	Let $z_{i,0} \in \Hprim$ for all $i \in \{1,\dots,n-1-f\}$ and
	\begin{align*}
		x_{1,k} &= \resolv_{\lambda A_1} z_{1,k}, \\
		x_{i,k} &= \resolv_{\theta^{-1} \lambda A_i}(x_{1,k} + \theta^{-1} z_{i,k})  &&  \text{for all } i \in \{2,\dots,n-1-f\}, \\
		x_{i,k} &= \lambda A_ix_{1,k}  &&  \text{for all } i \in \{n-f,\dots,n-1\}, \\
		\bar{x}_k &= \sum_{j=n-f}^{n-1} x_{j,k} + \sum_{j=2}^{n-1-f} (z_{j,k} + \theta(x_{1,k} - x_{j,k})), \\
		x_{n,k} &= \resolv_{\lambda A_n}(2x_{1,k} - z_{1,k} - \bar{x}_k), \\
		z_{i,k+1} &= z_{i,k} - \theta (x_{i,k} - x_{n,k}) && \text{for all } i \in \{1,\dots,n-1-f\}
	\end{align*}
	for all $k\in \N$ where $\gamma > 0$ and $\theta > 0$.
	If $A$ satisfies \cref{ass:coco} and
	\[
		\frac{\lambda}{2}\sum_{i=n-f}^{n-1}\beta_i^{-1} < 2 - \theta(n-1-f)
	\]
	then $x_{n,k} \weakto x^\star \in \zer \sum_{i=1}^n A_i$.
	\begin{proof}
		Let $T_{(\cdot)} \colon \Hprim^{n-1-f}\to\Hprim^{n-1-f}$ be the generalized primal-dual resolvent with representation $(n,M,MK,HMK,HM)$ with the matrices
		\begin{gather*}
			M
			=
			\begin{fmatrix}{c|c|c|c}
				1 & & & 1 \\
				\hline
				\mathbf{1} & \frac{1}{\theta} I & & \mathbf{1} \\
				\hline
				\mathbf{1} & & \phantom{\frac{1}{\theta}I} & \mathbf{1} \\
				\hline
				1 & & & 1 \\
			\end{fmatrix}
			\in \R^{n\times n}
			,\quad
			K
			=
			\begin{fmatrix}{c|c}
				\frac{1}{2} &  \\
				\hline
				\phantom{I} & I  \\
				\hline
				\phantom{I} &  \\
				\hline
				\frac{1}{2} & \\
			\end{fmatrix}
			\in\R^{n\times(n-1-f)}
			\\
			\text{and}\quad
			H =
			\theta
			\begin{fmatrix}{c|c|c|c}
				\frac{1}{2+f} & & \frac{1}{2+f}\mathbf{1}^T & \frac{1}{2+f} \\
				\hline
				\phantom{I} & I & &  \\
			\end{fmatrix}
			\in\R^{(n-1-f)\times n}
		\end{gather*}
		where $I \in \R^{(n-2-f) \times (n-2-f)}$ is the identity matrix, $\mathbf{1}$ are column vectors of ones with appropriate sizes, and empty blocks denote zero matrices.
		Note that the third rows and columns vanish completely when $f=0$ while the second rows and columns vanish when $f = n-2$.
		The matrix $M$ is an $n$-kernel over $\ops_n^F$ with rank $n-1-f$, which is minimal, see \cref{cor:fo-minlift-several}.
		Furthermore, these matrices satisfy $\ran K = \ran M^T = (\ker M)^\perp$ and $\ker H = (\ran H^T)^\perp = (\ran M)^\perp$, implying, through \cref{prop:kernel2split}, that $T_{(\cdot)}$ is a frugal splitting operator.
		For later use, we note that $T_{(\cdot)}$ can be represented as $(n,SUP,SU,U,UP)$ with $U = HMK$, $S = MK(HMK)^{-1}$, and $P = (HMK)^{-1}HM$.

		Consider first the case when $\lambda = 1$, the general case when $\lambda > 0$ will be shown later.
		In this case, with $z_k = (z_{1,k},\dots,z_{n-1-f,k})$, the update of $z_{k+1}$ in the theorem can be written as
		\[
			z_{k+1} = T_{A}z_k.
		\]
		It can also be verified
		that $(y_{1,k},\dots,y_{n,k}) = (M + \pdop_{A,n})^{-1}MKz_k$ satisfies
		\begin{align*}
			y_{1,k} &= z_{1,k} - x_{1,k}, \\
			y_{i,k} &= z_{1,k} + \theta(x_{1,k} - x_{i,k}) &&\text{for all } i \in \{2,\dots,n-1-f\}, \\
			y_{i,k} &= x_{i,k} &&\text{for all } i \in \{n-f,\dots,n\},
		\end{align*}
		for all $k\in\N$.
		Furthermore, if we choose
		\[
			Q =
			\theta^{-1}
			\begin{fmatrix}{c|c}
				1 &  \\
				\hline
				& I  \\
			\end{fmatrix}
			\in\R^{(n-1-f)\times(n-1-f)},
		\]
		then $I_F(P^TQ - S)U = (P^TQ - S)U$ where $I_F$ is defined in \cref{thm:conv} and condition \cref{thm:conv:struct-cond} of \cref{thm:conv} holds.
		Condition \cref{thm:conv:seq-cond} of \cref{thm:conv} with this $Q$ reads as
		\begin{align*}
			Q \succ 0
			\quad\text{and}\quad
			W
			&=
			\begin{fmatrix}{c|c}
				2 - \theta(n-1-f) - \tfrac{1}{2} \sum_{i=n-f}^{n-1}\beta_i^{-1} &  \\
				\hline
				& \frac{1}{\theta} I  \\
			\end{fmatrix}
			\succ 0
		\end{align*}
		which holds if and only if $\theta > 0$ and $\tfrac{1}{2} \sum_{i=n-f}^{n-1}\beta_i^{-1} < 2 - \theta(n-1-f)$.
		In this case, \cref{thm:conv} gives the convergence of
		\[
			(y_{1,k},\dots,y_{n,k}) \weakto (y_1^\star,\dots,y_n^\star) \in \zer \pdop_{A,n}
		\]
		and $x_{n,k} = y_{n,k} \weakto y_n^\star \in \zer \sum_{i=1}^n A_i$.

		Finally, consider the general case with arbitrary $\lambda > 0$.
		The update of $z_{k+1}$ in the theorem can be written as
		\[
			z_{k+1} = T_{\lambda A}z_k
		\]
		for all $k\in\N$ where $\lambda A = (\lambda A_1,\dots,\lambda A_n)$.
		Since $A_i$ is $\beta_i$-cocoercive for all $i\in F$, $\lambda A_i$ will be $\beta_i\lambda^{-1}$-cocoercive for all $i\in F$, implying that $\lambda A$  satisfies \cref{ass:coco}.
		Applying \cref{thm:conv} to this fixed point iteration with the scaled operator tuple gives the convergence of
		\[
			(y_{1,k},\dots,y_{n,k}) \to (y_1^\star,\dots,y_n^\star) \in \zer \pdop_{\lambda A,n}
		\]
		and $x_{n,k} = y_{n,k} \to y_n^\star \in \zer \sum_{i=1}^n \lambda A_i$ if $\theta > 0$ and $\tfrac{\lambda}{2} \sum_{i=n-f}^{n-1}\beta_i^{-1} < 2 - \theta(n-1-f)$.
		This concludes the proof.
	\end{proof}
\end{thm}

The relaxation factor $\theta$ enters also as a step-size scaling on some of the resolvents---$\resolv_{\theta^{-1}\lambda A_i}$ is evaluated for $i \in \{2,\dots,n-1-f\}$---while $\resolv_{\lambda A_i}$ is evaluated for $i \in \{1,n\}$.
This construction was necessary to ensure convergence.
\cref{thm:new-alg} prescribes that $\theta < \frac{2}{n-1-f}$ is needed for convergence, i.e., the relaxation parameter has to decrease with the number of operators evaluated via resolvents.
This is counteracted by the fact that the step-size of all but two of the resolvents is scaled with $\theta^{-1}$, a number that will increase with the number of operators.

There are no step-size restrictions when no forward evaluations are used, i.e., when $f = 0$ since $\lambda>0$ is free.
Furthermore, the step-size bound only depends on the sum of the inverse cocoercivity constants.
This is natural since all forward steps are evaluated at the same point and the results are simply added together.
The update can therefore equivalently be viewed as evaluating $\widehat{A} = \sum_{i=n-f}^{n-1}A_i$ instead of each operator individually and $\widehat{A}$ is a $(\sum_{i=n-f}^{n-1}\beta_i^{-1})^{-1}$-cocoercive operator.

\subsection{Relation to Other Methods With Minimal Lifting}
When $n = 3$, $f = 1$ and $\theta = 1$, this method reduces to the three operator splitting of Davis and Yin, see \cref{sec:rep-ex}.
When $n = 3$ and $f = 0$, the method is closely related to the three operator resolvent splitting operator of Ryu \cite{ryuUniquenessDRSOperator2020}.
That method uses a frugal splitting operator that is calculated as $(\hat{z}_1,\hat{z}_2) = T_{(A_1,A_2,A_3)}(z_1,z_2)$, where
\begin{align*}
	x_1 &= \resolv_{\lambda A_1}(z_1), \\
	x_2 &= \resolv_{\lambda A_2}(z_2 + x_1), \\
	x_3 &= \resolv_{\lambda A_3}(-z_1 - z_2 + x_1 + x_2), \\
	\hat{z}_1 &= z_1 + \theta(x_3 - x_1), \\
	\hat{z}_2 &= z_2 + \theta(x_3 - x_2).
\end{align*}
In the unrelaxed case with $\theta = 1$, this update is the same as our update.
For other choices of $\theta$, our proposed method differs in that the step-size in the computation of $x_2$ is scaled with $\theta^{-1}$.
As noted by Malitsky and Tam \cite{malitskyResolventSplittingSums2021}, a straightforward extension of Ryu's method to four operators fails to converge for all $\theta > 0$ in some cases and we found this step-size scaling to be the key that allowed us to establish convergence for $n > 3$.

In the case when $f = 0$ and $n>2$ is arbitrary, Malitsky and Tam \cite{malitskyResolventSplittingSums2021} presented a splitting method along with a proof of its minimal lifting.
It uses a frugal splitting operator where $(\hat{z}_1,\dots,\hat{z}_{n-1}) = T_{(A_1,\dots,A_n)}(z_1,\dots,z_{n-1})$ is calculated as
\begin{align*}
	x_1 &= \resolv_{\gamma A_1}(z_1), \\
	x_i &= \resolv_{\gamma A_i}(z_i - z_{i-1} + x_{i-1}) && \text{for all } i\in \{2,\dots,n-1\}, \\
	x_n &= \resolv_{\gamma A_n}(-z_{n-1} + x_1 + x_{n-1}), \\
	\hat{z}_i &= z_i + \theta(x_{i+1} - x_i) && \text{for all } i \in \{1,\dots,n-1\}.
\end{align*}
This method was later expanded to include forward evaluations in \cite{aragon-artachoDistributedForwardBackwardMethods2022} and the results of this paper prove that the method with forward evaluations also has minimal lifting.
This method is different from ours. In the Malitsky--Tam splitting operator, each resolvent depends on the previous one. This is in contrast to our method that allows $x_{i,k}$ to be calculated in parallel for all $i \in\{2,\dots,n-1\}$. Although not fully parallelizable, the evaluation of the Malitsky--Tam operator can be partially parallelized with $\dots,x_{i-2,k+1}, x_{i,k}, x_{i+2,k-1}, \dots$ being computable in parallel.
Comparing the step-sizes of the two methods, the Malitsky--Tam method converges for all $\theta \in (0,1)$ while our method requires that $\theta \in (0,\frac{2}{n-1})$, i.e., our method requires a smaller relaxation parameter for larger $n$.
As mentioned, our method actually increases the step-sizes for the resolvents of $A_2,\dots,A_{n-1}$ with $n$, which might offset the decreasing relaxation.

Another method similar to ours in the $f = 0$ case is the method presented by Campoy in \cite{campoyProductSpaceReformulation2022}.
It is based on Douglas--Rachford splitting applied to a product space reformulation of the finite sum monotone inclusion problem which results in a splitting operator where $(\hat{z}_1,\dots,\hat{z}_{n-1}) = T_A(z_1,\dots,z_{n-1})$ is such that
\begin{align*}
	x_1 &= \resolv_{\frac{\gamma}{n-1}A_1}(\tfrac{1}{n-1}\sum_{j=1}^{n-1}z_j) \\
	x_{i} &= \resolv_{\gamma A_{i}}(2x_1 - z_{i-1}) && \text{for all } i \in\{2,\dots,n\} \\
	\hat{z}_i &= z_i + \theta(x_{i+1} - x_1) && \text{for all } i \in\{1,\dots,n-1\}.
\end{align*}
This operator clearly has minimal lifting and its fixed point iteration converges for all $\theta \in (0,2)$, $\gamma > 0$ and $A \in \ops_n$.
As with our method, it is parallelizable and uses an uneven step-size with the first resolvent using a step-size of $\frac{\gamma}{n-1}$ while the others use a step-size of $\gamma$.
However, it does not appear possible to rewrite this splitting operator as a special case of ours, or vice versa.

A similar approach to the one used by Campoy in \cite{campoyProductSpaceReformulation2022} was also presented by Condat \textit{et al.}\ \cite{condatProximalSplittingAlgorithms2021} but restricted to the convex optimization case.
The splitting operator in \cite{condatProximalSplittingAlgorithms2021} is the same as the one in Campoy but with a weighted average instead of the arithmetic average in the first row.
Condat \textit{et al.}\ also applied other algorithm than Douglas--Rachford splitting to the product space reformulation, resulting in several different parallelizable splitting operators both with and without forward evaluations.
Most notably is perhaps a Davis--Yin based splitting operator over $\ops_{n}^{\{2\}}$ that is parallelizable and has minimal lifting, \cite[Equation (212)]{condatProximalSplittingAlgorithms2021}.
Although many of these methods are similar to ours, we were unable to rewrite either our or any of their methods with minimal lifting as special cases of each other.

The biggest difference between this work and the works of Campoy and Condat \textit{et al.}\ is perhaps conceptual.
While their parallelizable methods are based on applying existing splitting operators to reformulations of the problem, we directly search over all possible frugal splitting operators.
Our search is not exhaustive but is enabled by the representation theorem since it allows us to easily work with the entire class of frugal splitting operators and nothing else.

\ifdefined\shortversion
\else
	\section{Conclusion}\label{sec:conclusion}
	We have presented an explicit parameterization of all frugal splitting operators.
	The parameterization is in terms of what we call generalized primal-dual resolvents, and we have provided necessary and sufficient conditions for a generalized primal-dual resolvent being a frugal splitting operator.
	This allows for a unified analysis and both minimal lifting and convergence results that are applicable to all frugal splitting operators.
	The minimal lifting results of Ryu and Malitsky--Tam were expanded beyond resolvent splitting operators to general frugal splitting operators with forward evaluations and we showed that the lifting number depends on the order of forward and backward evaluations.
	We further presented a new convergent frugal splitting operator with minimal lifting that allows for most of the forward and/or backward steps to be computed in parallel.
	In the triple-backward case the method is the same as the minimal lifting method of Ryu if neither method uses relaxation.
	The slight difference in how relaxation is introduced is crucial to extend Ryu's method to an arbitrary number of operators.
	In the double-backward-single-forward case the method reduces to three operator splitting by Davis and Yin.
\fi

\paragraph{Acknowledgments}
We would like to thank Heinz Bauschke for his generous feedback on a preprint of this paper.

All authors have been supported by ELLIIT: Excellence Center at Linköping-Lund in Information Technology.
The first and last authors have also been funded by the Swedish Research Council.
The Wallenberg AI, Autonomous Systems and Software Program (WASP) have supported the work of the second and last author.

\bibliography{references}

\newpage
\appendix
\part*{Supplement}
\section*{Derivation of Representations and Convergence Conditions}

In this section, we will verify the representations and convergence conditions for the frugal splitting operators stated in the paper.
We will also derive representations and convergence conditions for a number of frugal splitting operators not previously mentioned.

We will keep a consistent notation in each of the examples presented.
For a representation of a frugal splitting operator $T_{(\cdot)}\colon\Hprim^d\to \Hprim^d$ over $\ops_n^F$, the goal is to find matrices $M$, $N$, $U$, $V$ such that
\begin{align*}
	y &= (M + \pdop_{(\cdot),p})^{-1}Nz, \\
	T_{(\cdot)}z &= z - Uz - Vy
\end{align*}
for some $p \in \{1,\dots,n\}$ and that \cref{thm:iff-representation} is satisfied.
Such a representation $(p,M,N,U,V)$ can always be factorized in terms of matrices $S$ and $P$ such that
\[
	M = SUP
	,\quad
	N = SU
	\quad\text{and}\quad
	V = UP,
\]
and where $\ran P \subseteq (\ker U)^\perp$ and $\ker S \supseteq (\ran U)^\perp$, see \cref{cor:factor-rep}.
The convergence conditions of \cref{thm:conv} for a fixed point iterations of $T_{A}$ for some $A\in\ops_n^F$ that satisfies \cref{ass:coco} is stated in terms of this factorization and can be written as
\begin{gather*}
	Q \succ 0
	,\\
	(I-I_F)(P^TQ - S)U = 0
	,\\
	W = QU + (QU)^T - U^TQU - \tfrac{1}{2}U^T(P^TQ - S)^T B^\dagger (P^TQ - S)U \succ 0
\end{gather*}
where $Q$ is a symmetric matrix that needs to be found for each frugal splitting operator.
When constructing new frugal splitting operators, we find it more convenient to work with a factorization of the representation $(p,M,N,U,V)$ in terms of matrices $H$ and $K$ such that
\[
	N = MK
	,\quad
	V = HM
	\quad\text{and}\quad
	U = HMK,
\]
where $\ran K = (\ker M)^\perp$ and $\ker H = (\ran M)^\perp$, see \cref{prop:kernel2split}.
This allows us to first design a kernel and then easily find $N$, $U$ and $V$ that results in a first order splitting.
Furthermore, some of the frugal splitting operators will be presented without step-size.
A variable step-size can be added to these methods simply by scaling the operator tuple in the same way as for our new frugal splitting operator in \cref{thm:new-alg}.

\subsection*{Forward-Backward}

The forward-backward operator \cite{goldsteinConvexProgrammingHilbert1964,levitinConstrainedMinimizationMethods1966} is
\[
	T_{(A_1,A_2)} = \resolv_{\gamma A_2}\circ(\Id- \gamma A_1)
\]
where $\gamma > 0$ and $A = (A_1,A_2) \in \ops_2^{\{1\}}$.
To derive a representation we choose primal index $p = 2$.
The next step would be to apply the Moreau identity to all backward steps with index $i \neq p$ but since there are no such backward steps we can directly define the results of each forward and backward evaluation
\begin{align*}
	y_1 &= A_1z, \\
	y_2 &= (\Id + \gamma A_2)^{-1}(z - \gamma y_1), \\
	T_{A}z &= y_2.
\end{align*}
We rewrite the first two lines such that the input is on the left and the result of all forward and backward steps are on the right
\begin{align*}
	z &\in A_1^{-1}y_1, \\
	\gamma^{-1}z &\in A_2y_2 + y_1 + \gamma^{-1}y_2, \\
	T_{A}z &= y_2.
\end{align*}
If we define $y = (y_1,y_2) \in \Hprim^2$ the first two lines can be written as
\[
	\begin{bmatrix} 1 \\ \gamma^{-1} \end{bmatrix} z
	\in
	\underbrace{\begin{bmatrix} A_1^{-1} & 0 \\ 0 & A_2 \end{bmatrix}}_{\pddiag_{A,2}} y +
	\begin{bmatrix} 0 & 0 \\ 1 & \gamma^{-1} \end{bmatrix} y
	=
	\underbrace{\begin{bmatrix} A_1^{-1} & -1 \\ 1 & A_2 \end{bmatrix}}_{\pdop_{A,2}} y +
	\begin{bmatrix} 0 & 1 \\ 0 & \gamma^{-1} \end{bmatrix} y
\]
which yields
\begin{align*}
	y &= \left(\begin{bmatrix} 0 & 1 \\ 0 & \gamma^{-1} \end{bmatrix} + \pdop_{A,2}  \right)^{-1} \begin{bmatrix} 1 \\ \gamma^{-1} \end{bmatrix}z, \\
	T_{A}z &= z - \begin{bmatrix}1\end{bmatrix} z + \begin{bmatrix} 0 & 1\end{bmatrix} y,
\end{align*}
and the matrices $M$, $N$, $U$, $V$ in the representation $(2,M,N,U,V)$ are then easily identified by comparing to \cref{def:gen-pd-resolv} which yields
\[
	M =
	\begin{bmatrix}
		0 & 1 \\
		0 & \gamma^{-1}
	\end{bmatrix}
	,\quad
	N =
	\begin{bmatrix}
		1 \\
		\gamma^{-1}
	\end{bmatrix}
	,\quad
	V =
	\begin{bmatrix}
		0 & 1
	\end{bmatrix}
	\quad\text{and}\quad
	U =
	\begin{bmatrix}
		1
	\end{bmatrix}.
\]
It is seen directly that
\[
	S =
	\begin{bmatrix}
		1 \\
		\gamma^{-1}
	\end{bmatrix}
	\quad\text{and}\quad
	P =
	\begin{bmatrix}
		0 & 1
	\end{bmatrix}
\]
provides a factorization of this representation as $(2,SUP,SU,U,UP)$.
For the convergence conditions, choosing $Q = \begin{bmatrix}\gamma^{-1}\end{bmatrix}$ where $I$ is the $2\times 2$ identity matrix it is obvious that $Q \succ 0$ and we further have
\begin{align*}
	(I-I_F)(P^TQ - S)U
	&= (I-I_F)
	\left(
		\gamma^{-1}
		\begin{bmatrix} 0 \\ 1\end{bmatrix}
		-
		\begin{bmatrix} 1 \\ \gamma^{-1} \end{bmatrix}
	\right) \\
	&= (I-I_F) \begin{bmatrix} -1 \\ 0 \end{bmatrix} \\
	&= 0
\end{align*}
and
\begin{align*}
	W
	&= QU + (QU)^T - U^TQU - \tfrac{1}{2}U^T(P^TQ -S)^TB^\dagger (P^T Q - S)U \\
	&=
	\begin{bmatrix}\gamma^{-1}\end{bmatrix}
	+\begin{bmatrix}\gamma^{-1}\end{bmatrix}
	-\begin{bmatrix}\gamma^{-1}\end{bmatrix}
	- \tfrac{1}{2} \begin{bmatrix}-1 & 0\end{bmatrix} \begin{bmatrix}\beta_1^{-1} & 0 \\ 0 & 0 \end{bmatrix} \begin{bmatrix} -1 \\ 0 \end{bmatrix} \\
	&= \begin{bmatrix}\gamma^{-1} - \tfrac{1}{2\beta_1}\end{bmatrix}.
\end{align*}
The final condition for convergence is then $W \succ 0$ which yields the well known result $\gamma < 2\beta_1$.

\subsection*{Douglas--Rachford}
The Douglas--Rachford splitting operator \cite{lionsSplittingAlgorithmsSum1979} is
\[
	T_{(A_1,A_2)} = \tfrac{1}{2}\Id + \tfrac{1}{2}(2\resolv_{\gamma A_2} - \Id)\circ(2\resolv_{\gamma A_1} - \Id)
\]
where $A = (A_1,A_2) \in \ops_2$ and $\gamma > 0$.
Consider the evaluation $\hat{z} = T_{(A_1,A_2)}z$ for some $z\in\Hprim$.
This can be written as
\begin{align*}
	x_1 &= \resolv_{\gamma A_1}z, \\
	y_2 &= \resolv_{\gamma A_2}(2x_1 - z), \\
	\hat{z} &= \tfrac{1}{2}z + \tfrac{1}{2}(2y_2 - (2x_1 - z)).
\end{align*}
We choose the primal index to $p = 2$ and apply the Moreau identity to the first resolvent,
\begin{align*}
	y_1 &= \resolv_{\gamma^{-1}A_1^{-1}}(\gamma^{-1}z), \\
	x_1 &= z - \gamma y_1, \\
	y_2 &= \resolv_{\gamma A_2}(2x_1 - z), \\
	\hat{z} &= \tfrac{1}{2}z + \tfrac{1}{2}(2y_2 - (2x_1 - z)).
\end{align*}
Eliminating $x_1$ gives
\begin{align*}
	y_1 &= \resolv_{\gamma^{-1}A_1^{-1}}(\gamma^{-1}z), \\
	y_2 &= \resolv_{\gamma A_2}(z - 2\gamma y_1), \\
	\hat{z} &= \gamma y_1 + y_2.
\end{align*}
Using the definition of the resolvent yields
\begin{align*}
	z &\in \gamma y_1 + A_1^{-1}y_1, \\
	\gamma^{-1}z &\in 2y_1 + \gamma^{-1}y_2 +  A_2y_2, \\
	\hat{z} &= z - z + \gamma y_1 + y_2.
\end{align*}
Regrouping such that $\pdop_{(A_1,A_2),2}$ can be identified finally gives
\begin{align*}
	z &\in [\gamma y_1 + y_2] + [A_1^{-1}y_1 - y_2], \\
	\gamma^{-1}z &\in [y_1 + \gamma^{-1}y_2] +  [A_2y_2 + y_1], \\
	\hat{z} &= z - [z] + [\gamma y_1 + y_2]
\end{align*}
and we can identify
\[
	M =
	\begin{bmatrix}
		\gamma & 1 \\
		1 & \gamma^{-1}
	\end{bmatrix}
	,\quad
	N =
	\begin{bmatrix}
		1 \\
		\gamma^{-1}
	\end{bmatrix}
	,\quad
	V =
	\begin{bmatrix}
		\gamma & 1
	\end{bmatrix}
	\quad\text{and}\quad
	U =
	\begin{bmatrix}
		1
	\end{bmatrix}.
\]
The representation can be factored as $(2,SUP,SU,U,UP)$ where
\[
	S = \begin{bmatrix} 1 \\ \gamma^{-1} \end{bmatrix}
	\quad\text{and}\quad
	P = \begin{bmatrix} \gamma & 1 \end{bmatrix}.
\]
For the convergence conditions, choosing $Q = \begin{bmatrix}\gamma^{-1}\end{bmatrix}$ yield $Q \succ 0$ and
\[
	(I - I_F)(P^TQ - S) = (P^T Q - S) = \gamma^{-1}\begin{bmatrix} \gamma \\ 1 \end{bmatrix} - \begin{bmatrix} 1 \\ \gamma^{-1}\end{bmatrix} = 0
\]
and
\begin{align*}
	W
	&= QU + (QU)^T - U^TQU - \tfrac{1}{2}U^T(P^TQ -S)^TB^\dagger (P^T Q - S)U \\
	&=
	Q + Q - Q
	=
	Q
\end{align*}
and hence is $W \succ 0$, i.e., fixed point iterations of the Douglas--Rachford splitting operator always converge.

\subsection*{Davis--Yin Three Operator Splitting}

To find a representation of three operator splitting of Davis--Yin \cite{davisThreeOperatorSplittingScheme2017},
\[
	T_{(A_1,A_2,A_3)} = \resolv_{\gamma A_3} \circ (2\resolv_{\gamma A_1} - \Id  - \gamma A_2 \circ \resolv_{\gamma A_1}) + \Id  - \resolv_{\gamma A_1}.
\]
where $\gamma > 0$ and $A=(A_1,A_2,A_3) \in \ops_3^{\{2\}}$ we choose $p = 3$.
Applying Moreau's identity to the resolvents of $A_i$ for all $i \neq p$ (in this case only $\resolv_{\gamma A_1}$) and defining the result of each forward and backward-step yields
\begin{align*}
	y_1 &= (\Id + \gamma^{-1} A_1^{-1})^{-1}(\gamma^{-1} z), \\
	y_2 &= A_2(z - \gamma y_1), \\
	y_3 &= (\Id + \gamma A_3)^{-1}( 2(z - \gamma y_1) - z - \gamma y_2 ), \\
	T_{A}z &= y_3 + z - (z - \gamma y_1)
\end{align*}
for all $z\in\Hprim$.
Rearranging the first three lines such that we only have $z$ on the left and $y_1$, $y_2$ and $y_3$ and unscaled operators on the right yields
\begin{align*}
	z &\in A_1^{-1}y_1 + \gamma y_1, \\
	z &\in A_2^{-1}y_2 + \gamma y_1, \\
	\gamma^{-1}z &\in A_3y_3 + 2y_1 + y_2 + \gamma^{-1}y_3,\\
	T_{A}z &= \gamma y_1 + y_3.
\end{align*}
If we define $y = (y_1,y_2,y_3)\in\Hprim^3$ we see that the first three lines can be written as
\[
	\begin{bmatrix} 1 \\ 1 \\ \gamma^{-1} \end{bmatrix} z
	\in
	\pddiag_{A,3} y +
	\begin{bmatrix} \gamma & 0 & 0 \\ \gamma & 0 & 0 \\ 2 & 1 & \gamma^{-1} \end{bmatrix} y
	=
	\pdop_{A,3} y +
	\begin{bmatrix} \gamma & 0 & 1 \\ \gamma & 0 & 1 \\ 1 & 0 & \gamma^{-1} \end{bmatrix} y
\]
and $T_{A}$ can then be written as
\begin{align*}
	y &= \left(\begin{bmatrix} \gamma & 0 & 1 \\ \gamma & 0 & 1 \\ 1 & 0 & \gamma^{-1} \end{bmatrix} + \pdop_{A,3}  \right)^{-1} \begin{bmatrix} 1 \\ 1 \\ \gamma^{-1} \end{bmatrix}z, \\
	T_{A}z &= z - \begin{bmatrix}1\end{bmatrix} z + \begin{bmatrix} \gamma & 0 & 1\end{bmatrix} y.
\end{align*}
From this we can easily identify the matrices $M$, $N$, $U$ and $V$ in the representation $(3,M,N,U,V)$ by comparing to \cref{def:gen-pd-resolv}.

\[
	M =
	\begin{bmatrix}
		\gamma & 0 & 1 \\
		\gamma & 0 & 1 \\
		1 & 0 & \gamma^{-1}
	\end{bmatrix}
	,\quad
	N = \begin{bmatrix} 1 \\ 1 \\ \gamma^{-1} \end{bmatrix}
	,\quad
	V = \begin{bmatrix} \gamma & 0 & 1 \end{bmatrix}
	\quad\text{and}\quad
	U = \begin{bmatrix} 1 \end{bmatrix}.
\]
This can be factored as $(3,SUP,SU,U,UP)$ where
\[
	S = \begin{bmatrix} 1 \\ 1 \\ \gamma^{-1} \end{bmatrix}
	\quad\text{and}\quad
	P = \begin{bmatrix} \gamma & 0 & 1 \end{bmatrix}.
\]
The convergence conditions are satified by $Q = \begin{bmatrix} \gamma^{-1}\end{bmatrix}$.
We see that $Q \succ 0$ and
\begin{align*}
	(I-I_F)(P^TQ - S)U
	&= (I-I_F)(P^TQ - S) \\
	&= (I-I_F)
	\left(
		\gamma^{-1}
		\begin{bmatrix} \gamma \\ 0 \\ 1 \end{bmatrix}
		-
		\begin{bmatrix} 1 \\ 1 \\ \gamma^{-1} \end{bmatrix}
	\right) \\
	&= 0
\end{align*}
and
\begin{align*}
	W =
	&QU + (QU)^{T} - U^TQU - \tfrac{1}{2}U^T(P^TQ-S)^T B^\dagger (P^TQ-S)U \\
	&= Q - \tfrac{1}{2}(P^TQ-S)^T B^\dagger (P^TQ-S) \\
	&= \begin{bmatrix}\gamma^{-1}\end{bmatrix} - \tfrac{1}{2}\begin{bmatrix} 0 & -1 & 0\end{bmatrix} B^\dagger \begin{bmatrix} 0 \\ -1 \\ 0 \end{bmatrix} \\
	&= \begin{bmatrix}\gamma^{-1} - \tfrac{1}{2\beta_2}\end{bmatrix}
\end{align*}
and $W \succ 0$ if $\gamma < 2 \beta_2$.

\subsection*{Forward-Backward with Momentum on the Forward Step}
Forward-Backward with momentum on the forward step \cite{moudafiConvergenceSplittingInertial2003} can be written as a fixed point iteration of the frugal splitting operator
\[
	T_{(A_1,A_2)}(z_1,z_2)
	=
	\begin{pmatrix}
		\hat{z}_1 \\
		\hat{z}_2
	\end{pmatrix}
	=
	\begin{pmatrix}
		\resolv_{\gamma A_2}(z_1 - \gamma A_1z_1 + \theta z_2 ) \\
		\hat{z}_1 - z_1
	\end{pmatrix}
\]
where $\gamma > 0$, $\theta \in \R$ and $A = (A_1,A_2) \in \ops_2^{\{1\}}$.
Note, there are other frugal splitting operators whose fixed point iteration are equivalent to this forward-backward method with momentum.
We claim that $(2,SUP,SU,U,UP)$ where
\[
	U = \begin{bmatrix}1 & 0 \\ 1 & 1 \end{bmatrix}
	,\quad
	S = \begin{bmatrix} 1 & 0 \\ \gamma^{-1}(1-\theta) & \gamma^{-1}\theta \end{bmatrix}
	\quad\text{and}\quad
	P = \begin{bmatrix} 0 & 1 \\ 0 & 0 \end{bmatrix}
\]
is a representation of this frugal splitting operator.
To show this we first note that
\begin{gather*}
	V = UP =
	\begin{bmatrix}
		0 & 1 \\ 0 & 1
	\end{bmatrix}
	,\quad
	M = SUP =
	\begin{bmatrix}
		0 & 1 \\ 0 & \gamma^{-1}
	\end{bmatrix}
	\quad\text{and} \\
	N = SU =
	\begin{bmatrix} 1 & 0 \\ \gamma^{-1}(1-\theta) & \gamma^{-1}\theta \end{bmatrix}
	\begin{bmatrix} 1 & 0 \\ 1 & 1 \end{bmatrix}
	=
	\begin{bmatrix} 1 & 0 \\ \gamma^{-1} & \gamma^{-1}\theta \end{bmatrix}
\end{gather*}
and
\begin{align*}
	Nz-My &\in \pdop_{A,2}y, \\
	\hat{z} &= z - Uz + Vy
\end{align*}
can then be written as
\begin{align*}
	z_1 - y_2 &\in A_1^{-1}y_1 - y_2, \\
	\gamma^{-1}(z_1 + \theta z_2) - \gamma^{-1}y_2 &\in A_2y_2 + y_1, \\
	\hat{z}_1 &= z_1 - z_1 + y_2, \\
	\hat{z}_2 &= z_2 - z_1 - z_2 + y_2
\end{align*}
which after some rearranging gives
\begin{align*}
	y_1 &\in A_1z_1, \\
	z_1 - \gamma y_1 + \theta z_2 &\in (\Id + \gamma A_2)y_2, \\
	\hat{z}_1 &= y_2, \\
	\hat{z}_2 &= y_2 - z_1.
\end{align*}
Rewriting the second line as a resolvent and combining the first and second gives
\begin{align*}
	y_2 &= \resolv_{\gamma A_2}(z_1 - \gamma A_1z_1 + \theta z_2), \\
	\hat{z}_1 &= y_2, \\
	\hat{z}_2 &= y_2 - z_1
\end{align*}
which is exactly the frugal splitting operator above.
For the convergence, if we choose
\[
	Q =
	\gamma^{-1}
	\begin{bmatrix}
		1 - \theta & \theta \\
		\theta & |\theta| + \epsilon
	\end{bmatrix}
\]
where $\epsilon > 0$ then
\begin{align*}
	(I-I_F)(P^TQ - S)U
	&= (I-I_F)
	\gamma^{-1}
	\left(
		\begin{bmatrix} 0 & 0 \\ 1 & 0\end{bmatrix}
		\begin{bmatrix} 1-\theta & \theta \\ \theta & (|\theta|+\epsilon) \end{bmatrix}
		-
		\begin{bmatrix} \gamma & 0 \\ 1-\theta & \theta \end{bmatrix}
	\right)U \\
	&= (I-I_F)
	\gamma^{-1}
	\left(
		\begin{bmatrix} 0 & 0 \\ 1-\theta & \theta \end{bmatrix}
		-
		\begin{bmatrix} \gamma & 0 \\ (1-\theta) & \theta \end{bmatrix}
	\right)U \\
	&= (I-I_F)
	\begin{bmatrix} -1 & 0 \\ 0 & 0 \end{bmatrix}
	\begin{bmatrix} 1 & 0 \\ 1 & 1 \end{bmatrix}
	\\
	&=
	\begin{bmatrix} 0 & 0 \\ 0 & 1\end{bmatrix}
	\begin{bmatrix} -1 & 0 \\ 0 & 0 \end{bmatrix} \\
	&= 0.
\end{align*}
Furthermore, we have
\begin{align*}
	QU + (QU)^T - U^TQU
	&=
	\gamma^{-1}\begin{bmatrix} 1-\theta & \theta \\ \theta & |\theta|+\epsilon \end{bmatrix}
	\begin{bmatrix} 1 & 0 \\ 1 & 1 \end{bmatrix}
	+ (QU)^T - U^TQU \\
	&=
	\gamma^{-1}\begin{bmatrix} 1 & \theta \\ \theta+|\theta|+\epsilon & |\theta|+\epsilon \end{bmatrix}
	+ (QU)^T - U^TQU \\
	&=
	\gamma^{-1}\begin{bmatrix} 2 & 2\theta+|\theta|+\epsilon \\ 2\theta+|\theta|+\epsilon & 2(|\theta|+\epsilon) \end{bmatrix} \\
	&\quad-
	\gamma^{-1}\begin{bmatrix} 1 & 1 \\ 0 & 1 \end{bmatrix}
	\begin{bmatrix} 1 & \theta \\ \theta+|\theta|+\epsilon & |\theta|+\epsilon \end{bmatrix} \\
	&=
	\gamma^{-1}\begin{bmatrix} 2 & 2\theta+|\theta|+\epsilon \\ 2\theta+|\theta|+\epsilon & 2(|\theta|+\epsilon) \end{bmatrix} \\
	&\quad-
	\gamma^{-1}\begin{bmatrix} 1+\theta+|\theta|+\epsilon & \theta+|\theta|+\epsilon \\ \theta+|\theta|+\epsilon & |\theta|+\epsilon \end{bmatrix} \\
	&=
	\gamma^{-1}\begin{bmatrix} 1-\theta-|\theta|-\epsilon & \theta \\ \theta & |\theta|+\epsilon \end{bmatrix}
\end{align*}
and
\begin{align*}
	\tfrac{1}{2}U^T(P^TQ -S)^TB^\dagger (P^T Q - S)U
	&=
	\tfrac{1}{2}
	\begin{bmatrix} -1 & 0 \\ 0 & 0 \end{bmatrix}
	\begin{bmatrix} \beta_1^{-1} & 0 \\ 0 & 0 \end{bmatrix}
	\begin{bmatrix} -1 & 0 \\ 0 & 0 \end{bmatrix} \\
	&=
	\gamma^{-1}\begin{bmatrix} \tfrac{\gamma}{2\beta_1} & 0 \\ 0 & 0 \end{bmatrix}
\end{align*}
which gives
\begin{align*}
	W &= QU + (QU)^T - U^TQU - \tfrac{1}{2}U^T(P^TQ -S)^TB^\dagger (P^T Q - S)U \\
	&=
	\gamma^{-1}\begin{bmatrix} 1-\theta-|\theta|-\epsilon & \theta \\ \theta & |\theta|+\epsilon \end{bmatrix}
	-\gamma^{-1}\begin{bmatrix} \tfrac{\gamma}{2\beta_1} & 0 \\ 0 & 0 \end{bmatrix} \\
	&=
	\gamma^{-1}\begin{bmatrix} 1-\theta-|\theta|-\epsilon - \tfrac{\gamma}{2\beta_1} & \theta \\ \theta & |\theta|+\epsilon \end{bmatrix}.
\end{align*}
The final condition for convergence is then $Q \succ 0$ and $W \succ 0$ which hold if
\[
	0 < 1 - \theta - \tfrac{\theta^2}{|\theta|+\epsilon}
	\quad\text{and}\quad
	0 < 1 - \theta - \tfrac{\theta^2}{|\theta|+\epsilon} - |\theta|-\epsilon - \tfrac{\gamma}{2\beta_1}.
\]
The first condition is clearly implied by the second since $\epsilon > 0$, $\gamma > 0$ and $\beta_1 > 0$.
There exists $\epsilon > 0$ such that the second condition hold if
\[
	0 < 1 - \theta - 2|\theta| - \tfrac{\gamma}{2\beta_1}.
\]
This is the same conditions as was derived in \cite{morinNonlinearForwardBackwardSplitting2022}.

\subsection*{Forward-Backward with Nesterov-like Momentum}
The update of forward-backward with Nesterov-like momentum \cite{nesterovMethodSolvingConvex1983,beckFastIterativeShrinkageThresholding2009} can be seen as the application of the following frugal splitting operator
\[
	T_{(A_1,A_2)}(z_1,z_2)
	=
	\begin{pmatrix}
		\hat{z}_1 \\
		\hat{z}_2
	\end{pmatrix}
	=
	\begin{pmatrix}
		\resolv_{\gamma A_2}(z_1 + \theta z_2 - \gamma A_1(z_1 + \theta z_2)) \\
		\hat{z}_1 - z_1
	\end{pmatrix}
\]
where $\gamma > 0$, $\theta > 0$ and $A = (A_1,A_2) \in \ops_2^{\{1\}}$.
As with forward-backward with momentum on the forward step presented earlier, there are other frugal splitting operators that also would yield a Nesterov-like momentum update.
To derive a representation, we choose the primal index $p = 2$ and note that the frugal splitting operator can be written as
\begin{align*}
	y_1 &= A_1(z_1 + \theta z_2), \\
	y_2 &= (\Id + \gamma A_2)^{-1} (z_1 + \theta z_2 - \gamma y_1), \\
	\hat{z}_1 &= y_2, \\
	\hat{z}_2 &= y_2 - z_1.
\end{align*}
Inverting $A_1$ and $\Id + \gamma A_2$ yield
\begin{align*}
	z_1 + \theta z_2 &\in A_1^{-1}y_1, \\
	\gamma^{-1}z_1 + \gamma^{-1}\theta z_2 - y_1 & \in A_2y_2 + \gamma^{-1}y_2, \\
	\hat{z}_1 &= z_1 - z_1 + y_2, \\
	\hat{z}_2 &= z_2 - z_1 - z_2 + y_2.
\end{align*}
Rearranging so that $\pdop_{(A_1,A_2),2}$ can be identified gives
\begin{align*}
	z_1 + \theta z_2 &\in [y_2] + [A_1^{-1}y_1 - y_2], \\
	\gamma^{-1}z_1 + \gamma^{-1}\theta z_2 &\in [\gamma^{-1}y_2] + [A_2y_2 + y_1], \\
	\hat{z}_1 &= z_1 - [z_1] + [y_2], \\
	\hat{z}_2 &= z_2 - [z_1 + z_2] + [y_2]
\end{align*}
and we can identify
\[
	M = \begin{bmatrix} 0 & 1 \\ 0 & \gamma^{-1} \end{bmatrix}
	,\quad
	N = \begin{bmatrix} 1 & \theta \\ \gamma^{-1} & \gamma^{-1}\theta \end{bmatrix}
	,\quad
	V = \begin{bmatrix} 0 & 1 \\ 0 & 1 \end{bmatrix}
	\quad\text{and}\quad
	U = \begin{bmatrix} 1 & 0 \\ 1 & 1 \end{bmatrix}
\]
This representation can be factored as$(2,SUP,SU,U,UP)$ where
\[
	S = \begin{bmatrix} 1-\theta & \theta \\ \gamma^{-1}(1-\theta) & \gamma^{-1}\theta \end{bmatrix}
	\quad\text{and}\quad
	P = \begin{bmatrix} 0 & 1 \\ 0 & 0 \end{bmatrix}.
\]
For the convergence conditions we choose
\[
	Q = \gamma^{-1}\begin{bmatrix} 1 - \theta & \theta \\ \theta & a \end{bmatrix}
\]
where $a > 0$.
This yields
\begin{align*}
	0
	&= (I-I_F)(P^TQ - S)U \\
	&= (I-I_F)
	\left(
		\gamma^{-1}
		\begin{bmatrix} 0 & 0 \\ 1 & 0 \end{bmatrix}
		\begin{bmatrix} 1-\theta & \theta \\ \theta & a \end{bmatrix}
		-
		\begin{bmatrix} 1-\theta & \theta \\ \gamma^{-1}(1-\theta) & \gamma^{-1}\theta \end{bmatrix}
	\right)U \\
	&= (I-I_F)
	\left(
		\gamma^{-1}
		\begin{bmatrix} 0 & 0 \\ 1-\theta & \theta \end{bmatrix}
		-
		\begin{bmatrix} 1-\theta & \theta \\ \gamma^{-1}(1-\theta) & \gamma^{-1}\theta \end{bmatrix}
	\right)U \\
	&= (I-I_F) \begin{bmatrix} -1+\theta & -\theta \\ 0 & 0 \end{bmatrix} U \\
	&= (I-I_F) \begin{bmatrix} -1+\theta & -\theta \\ 0 & 0 \end{bmatrix} \begin{bmatrix} 1 & 0 \\ 1 & 1 \end{bmatrix} \\
	&= (I-I_F) \begin{bmatrix} -1 & -\theta \\ 0 & 0 \end{bmatrix} \\
	&= 0
\end{align*}
and
\begin{align*}
	QU + (QU)^T - U^TQU
	&=
	\gamma^{-1}\begin{bmatrix} 1-\theta & \theta \\ \theta & a \end{bmatrix}
	\begin{bmatrix} 1 & 0 \\ 1 & 1 \end{bmatrix}
	+ (QU)^T - U^TQU  \\
	&=
	\gamma^{-1}\begin{bmatrix} 1 & \theta \\ \theta+a & a \end{bmatrix}
	+ (QU)^T - U^TQU  \\
	&=
	\gamma^{-1}\begin{bmatrix} 2 & 2\theta+a \\ 2\theta+a & 2a \end{bmatrix}
	-
	\gamma^{-1}\begin{bmatrix} 1 & 1 \\ 0 & 1 \end{bmatrix}
	\begin{bmatrix} 1 & \theta \\ \theta+a & a \end{bmatrix} \\
	&=
	\gamma^{-1}\begin{bmatrix} 2 & 2\theta+a \\ 2\theta+a & 2a \end{bmatrix}
	-
	\gamma^{-1}\begin{bmatrix} 1+\theta+a & \theta+a \\ \theta+a & a \end{bmatrix} \\
	&=
	\gamma^{-1}\begin{bmatrix} 1-\theta-a & \theta \\ \theta & a \end{bmatrix}
\end{align*}
and
\begin{align*}
	\tfrac{1}{2}U^T(P^TQ -S)^TB^\dagger (P^T Q - S)U
	&=
	\tfrac{1}{2}
	\begin{bmatrix} -1 & 0 \\ -\theta & 0 \end{bmatrix}
	\begin{bmatrix} \beta_1^{-1} & 0 \\ 0 & 0 \end{bmatrix}
	\begin{bmatrix} -1 & -\theta \\ 0 & 0 \end{bmatrix} \\
	&=
	\tfrac{1}{2\beta_1}
	\begin{bmatrix} -1 & 0 \\ -\theta & 0 \end{bmatrix}
	\begin{bmatrix} -1 & -\theta \\ 0 & 0 \end{bmatrix} \\
	&=
	\tfrac{1}{2\beta_1}
	\begin{bmatrix} 1 & \theta \\ \theta & \theta^2 \end{bmatrix}.
\end{align*}
Hence we have
\begin{align*}
	W
	&= QU + (QU)^T - U^TQU - \tfrac{1}{2}U^T(P^TQ - S)B^\dagger (P^TQ - S)U \\
	&=
	\gamma^{-1}
	\begin{bmatrix} 1-\theta-a & \theta \\ \theta & a \end{bmatrix}
	-
	\tfrac{1}{2\beta_1}
	\begin{bmatrix} 1 & \theta \\ \theta & \theta^2 \end{bmatrix} \\
	&= \gamma^{-1}\begin{bmatrix} 1-\theta-\hat{\gamma}-a & \theta(1-\hat{\gamma}) \\ \theta(1-\hat{\gamma}) & a-\theta^2\hat{\gamma} \end{bmatrix}
\end{align*}
where $\hat{\gamma} = \tfrac{\gamma}{2\beta_1}$.
It can be verified that $Q-W \succ 0$ and the final condition for convergence---$Q\succ 0$ and $W\succ 0$---then holds if $W\succ 0$ which is equivalent to
\[
	0 < a - \theta^2\hat{\gamma}
	\quad\text{and}\quad
	0 < 1 - \theta - \hat{\gamma} - a - \tfrac{\theta^2(1- \hat{\gamma})^2}{a - \theta^2\hat{\gamma}}.
\]

\subsection*{A New Frugal Splitting Operator with Minimal Lifting}
We will here derive a frugal splitting operator with minimal lifting that serves as the base of the algorithm in \cref{thm:new-alg}.
This was also done in the proof of \cref{thm:new-alg} but we will here be more detailed and verify the matrix calculations more carefully.

The frugal splitting operator will be over $\ops_n^F$ where $F = \{n-f,\dots,n-1\}$ and $f = |F|$ and for notational convenience we let $R = \{2,\dots,n-f-1\}$ with $r = |R|$.
The representation of the new splitting is then $(n,M,MK,HMK,HM)$ where
\begin{gather*}
	M
	=
	\begin{fmatrix}{c|c|c|c}
		1 & & & 1 \\
		\hline
		\mathbf{1} & \frac{1}{\theta} I & & \mathbf{1} \\
		\hline
		\mathbf{1} & & \phantom{\frac{1}{\theta}I} & \mathbf{1} \\
		\hline
		1 & & & 1 \\
	\end{fmatrix}
	\in \R^{n\times n}
	,\quad
	K
	=
	\begin{fmatrix}{c|c}
		\frac{1}{2} &  \\
		\hline
		\phantom{I} & I  \\
		\hline
		\phantom{I} &  \\
		\hline
		\frac{1}{2} & \\
	\end{fmatrix}
	\in\R^{n\times(1+r)}
	\quad\text{and}\\
	H =
	\theta
	\begin{fmatrix}{c|c|c|c}
		\frac{1}{2+f} & & \frac{1}{2+f}\mathbf{1}^T & \frac{1}{2+f} \\
		\hline
		\phantom{I} & I & &  \\
	\end{fmatrix}
	\in\R^{(1+r)\times n}
\end{gather*}
where $\theta > 0$, $I$ is the identity matrix in $\R^{r \times r}$, $\mathbf{1}$ are column vectors of ones with appropriate sizes and empty block denotes zero matrices.
Since \cref{prop:kernel2split} is satisfied by these matrices, the operator corresponding to this representation is a frugal splitting operator.
Let $T_{(\cdot)} \colon \Hprim^{(1+r)}\to \Hprim^{(1+r)}$ be this frugal splitting operator and consider the evaluation of $\hat{z} = T_{A}z$ where $A = (A_1,\dots,A_n) \in \ops_n^F$.
This can be written as
\begin{align*}
	y &= (M + \pdop_{A,n})^{-1}MKz \\
	\hat{z} &= z - HMKz + HMy.
\end{align*}
Using the definition of the inverse of $M + \pdop_{A,n}$ and writing out the matrices explicitly gives
\begin{align*}
	&\begin{fmatrix}{c|c}
		1 &  \\
		\hline
		\mathbf{1} & \frac{1}{\theta} I  \\
		\hline
		\mathbf{1} &  \\
		\hline
		1 & \\
	\end{fmatrix}
	z
	\in
	\begin{fmatrix}{c|c|c|c}
		1 & & &  \\
		\hline
		\mathbf{1} & \frac{1}{\theta} I & &  \\
		\hline
		\mathbf{1} & & \phantom{\frac{1}{\theta}I} & \ \\
		\hline
		2 & \mathbf{1} & \mathbf{1} & 1 \\
	\end{fmatrix}
	y
	+
	\begin{bmatrix}
		A_1^{-1} & & & \\
		& \ddots & & \\
		& & A_{n-1}^{-1} & \\
		& & & A_n
	\end{bmatrix}
	y
	\\
	&\hat{z}
	=
	z
	-
	\begin{fmatrix}{c|c}
		\theta & \phantom{I} \\
		\hline
		\theta\mathbf{1} & I\\
	\end{fmatrix}
	z
	+
	\theta
	\begin{fmatrix}{c|c|c|c}
		1 & & & 1 \\
		\hline
		\mathbf{1} & \frac{1}{\theta} I & \phantom{I} & \mathbf{1} \\
	\end{fmatrix}
	y
\end{align*}
Setting $z = (z_1,\dots,z_{1+r})$, $\hat{z} = (\hat{z}_1,\dots,\hat{z}_{1+r})$, $y = (y_1,\dots,y_n)$ and writing out the corresponding equations for these matrix/operator expressions gives
\begin{align*}
	z_1 &\in (\Id + A_1^{-1})y_1, \\
	z_1 + \tfrac{1}{\theta}z_i &\in y_1 + (\tfrac{1}{\theta}\Id + A_i^{-1})y_i  &&  \text{for all } i \in R, \\
	z_1 &\in y_1 + A_i^{-1}y_i  &&  \text{for all } i \in F, \\
	z_1 &\in y_1 + \sum_{j=1}^{n-1}y_j + (\Id + A_n)y_n, \\
	\hat{z}_1 &= (1-\theta)z_1 + \theta(y_1 + y_n), \\
	\hat{z}_i &= -\theta z_1 + \theta(y_1 + y_n + \tfrac{1}{\theta}y_i) && \text{for all } i \in R.
\end{align*}
Rearranging each line yields
\begin{align*}
	y_1 &= (\Id + A_1^{-1})^{-1}z_1, \\
	y_i &= (\Id + \theta A_i^{-1})^{-1}(\theta [z_1-y_1] + z_i)  &&  \text{for all } i \in R, \\
	y_i &= A_i(z_1-y_1)  &&  \text{for all } i \in F, \\
	y_n &= (\Id + A_n)^{-1}([z_1-y_1] - \sum_{j=1}^{n-1}y_j), \\
	\hat{z}_1 &= z_1 - \theta [z_1 - y_1] + \theta y_n, \\
	\hat{z}_i &= y_i - \theta [z_1-y_1] + \theta y_n && \text{for all } i \in R.
\end{align*}
Applying the Moreau identity to the first two lines and introducing variables $x_i$ for all $x\in\{1,\dots,n\}$ that contains the results of all primal evaluations yields
\begin{align*}
	x_1 &= (\Id + A_1)^{-1}z_1, \\
	x_i &= (\Id + \theta^{-1} A_i)^{-1}([z_1-y_1] + \theta^{-1} z_i)  &&  \text{for all } i \in R, \\
	x_i &= A_i(z_1-y_1)  &&  \text{for all } i \in F, \\
	x_n &= (\Id + A_n)^{-1}([z_1-y_1] - \sum_{j=1}^{n-1}y_j), \\
	y_1 &= z_1 - x_1, \\
	y_i &= \theta [z_1-y_1] + z_i - \theta x_i = z_i + \theta(x_1 - x_i)  &&  \text{for all } i \in R, \\
	y_i &= x_i  &&  \text{for all } i \in F, \\
	y_n &= x_n, \\
	\hat{z}_1 &= z_1 - \theta[z_1 - y_1] + \theta y_n, \\
	\hat{z}_i &= y_i - \theta [z_1-y_1] + \theta y_n && \text{for all } i \in R.
\end{align*}
Eliminate the $y_i$ variables gives
\begin{align*}
	x_1 &= (\Id + A_1)^{-1}z_1, \\
	x_i &= (\Id + \theta^{-1} A_i)^{-1}(x_1 + \theta^{-1} z_i)  &&  \text{for all } i \in R, \\
	x_i &= A_ix_1  &&  \text{for all } i \in F, \\
	\bar{x} &= \sum_{j\in F}x_j + \sum_{j\in R}(z_j + \theta(x_1 - x_j)), \\
	x_n &= (\Id + A_n)^{-1}(2x_1 - z_1 - \bar{x}), \\
	\hat{z}_1 &= z_1 - \theta (x_1 - x_n), \\
	\hat{z}_i &= z_i - \theta (x_i - x_n) && \text{for all } i \in R.
\end{align*}
This concludes the derivations of the frugal splitting operator.
To derive conditions for convergence we need a representation of the form $(n,SUP,SU,U,UP)$ and since $U = HMK$ is invertible this is easily calculated as
\begin{align*}
	U = HMK =
	\begin{fmatrix}{c|c}
		\theta & \phantom{I} \\
		\hline
		\theta\mathbf{1} & I\\
	\end{fmatrix}
	\in \R^{(1+r)\times (1+r)}
	&,\quad
	U^{-1} =
	\begin{fmatrix}{c|c}
		\theta^{-1} &  \\
		\hline
		-\mathbf{1} & I  \\
	\end{fmatrix}
	\in \R^{(1+r)\times (1+r)}
	,\\
	P = U^{-1}HM =
	\begin{fmatrix}{c|c|c|c}
		1 & & & 1 \\
		\hline
		& I & \phantom{I} & \\
	\end{fmatrix}
	\in \R^{(1+r)\times n}
	&,\quad
	S = MKU^{-1} =
	\theta^{-1}
	\begin{fmatrix}{c|c}
		1 &  \\
		\hline
		& I  \\
		\hline
		\mathbf{1} &  \\
		\hline
		1 & \\
	\end{fmatrix}
	\in \R^{n\times (1+r)}
\end{align*}
Choosing
\[
	Q =
	\theta^{-1}
	\begin{fmatrix}{c|c}
		1 &  \\
		\hline
		& I  \\
	\end{fmatrix}
	\in\R^{(1+r)\times(1+r)}
\]
then results in
\begin{align*}
	(I - I_F)(P^TQ - S)U
	&=
	(I-I_F)
	\left(
		\theta^{-1}
		\begin{fmatrix}{c|c}
			1 &  \\
			\hline
			& I  \\
			\hline
			\phantom{1} & \\
			\hline
			1 &  \\
		\end{fmatrix}
		-
		\theta^{-1}
		\begin{fmatrix}{c|c}
			1 &  \\
			\hline
			& I  \\
			\hline
			\mathbf{1} &  \\
			\hline
			1 & \\
		\end{fmatrix}
	\right)
	U
	\\
	&=
	(I-I_F)
	\theta^{-1}
	\begin{fmatrix}{c|c}
		&  \\
		\hline
		&   \\
		\hline
		-\mathbf{1} & \phantom{I}  \\
		\hline
		& \\
	\end{fmatrix}
	\begin{fmatrix}{c|c}
		\theta & \phantom{I} \\
		\hline
		\theta\mathbf{1} & I\\
	\end{fmatrix}
	\\
	&=
	(I-I_F)
	\begin{fmatrix}{c|c}
		&  \\
		\hline
		&   \\
		\hline
		-\mathbf{1} & \phantom{I} \\
		\hline
		& \\
	\end{fmatrix}
	\\
	&= 0.
\end{align*}
We further have
\begin{align*}
	QU &+ (QU)^T - U^TQU \\
	&= \theta^{-1}(U + U^T - U^TU) \\
	&=
	\theta^{-1}
	\left(
		\begin{fmatrix}{c|c}
			\theta & \phantom{I} \\
			\hline
			\theta\mathbf{1} & I\\
		\end{fmatrix}
		+
		\begin{fmatrix}{c|c}
			\theta & \theta\mathbf{1}^T \\
			\hline
			\phantom{I} & I\\
		\end{fmatrix}
		-
		\begin{fmatrix}{c|c}
			\theta & \theta\mathbf{1}^T \\
			\hline
			\phantom{I} & I\\
		\end{fmatrix}
		\begin{fmatrix}{c|c}
			\theta & \phantom{I} \\
			\hline
			\theta\mathbf{1} & I\\
		\end{fmatrix}
	\right)
	\\
	&=
	\theta^{-1}
	\left(
		\begin{fmatrix}{c|c}
			2\theta & \theta\mathbf{1}^T \\
			\hline
			\theta\mathbf{1} & 2I\\
		\end{fmatrix}
		-
		\begin{fmatrix}{c|c}
			\theta^2(n-1-f) & \theta\mathbf{1}^T  \\
			\hline
			\theta\mathbf{1} & I  \\
		\end{fmatrix}
	\right)
	\\
	&=
	\begin{fmatrix}{c|c}
		2& \mathbf{1}^T \\
		\hline
		\mathbf{1} & 2\theta^{-1}I\\
	\end{fmatrix}
	-
	\begin{fmatrix}{c|c}
		\theta(1+r) & \mathbf{1}^T  \\
		\hline
		\mathbf{1} & \theta^{-1}I  \\
	\end{fmatrix}
	\\
	&=
	\begin{fmatrix}{c|c}
		2 - \theta(1+r)&  \\
		\hline
		 & \theta^{-1}I\\
	\end{fmatrix}
\end{align*}
and
\begin{align*}
	\tfrac{1}{2}U^T(P^TQ - S)B^\dagger(P^TQ - S)U
	&=
	\tfrac{1}{2}
	\begin{fmatrix}{c|c|c|c}
		\phantom{I} & \phantom{I} & - \mathbf{1}^T & \phantom{I} \\
		\hline
		& & &
	\end{fmatrix}
	B^\dagger
	\begin{fmatrix}{c|c}
		&  \\
		\hline
		&   \\
		\hline
		-\mathbf{1} & \phantom{I} \\
		\hline
		& \\
	\end{fmatrix}
	\\
	&=
	\tfrac{1}{2}
	\begin{fmatrix}{c|c}
		\sum_{i\in F}\beta_i^{-1} &  \\
		\hline
		& \phantom{I}  \\
	\end{fmatrix}
\end{align*}
which, since $1+r = n-1-f$, results in
\begin{align*}
	W
	&= QU + (QU)^T - U^TQU - \tfrac{1}{2}U^T(P^TQ - S)B^\dagger(P^TQ - S)U \\
	&=
	\begin{fmatrix}{c|c}
		2 - \theta(1+r) &  \\
		\hline
		& \frac{1}{\theta} I  \\
	\end{fmatrix}
	-
	\tfrac{1}{2}
	\begin{fmatrix}{c|c}
		\sum_{i\in F}\beta_i^{-1} &  \\
		\hline
		& \phantom{I}  \\
	\end{fmatrix}
	\\
	&=
	\begin{fmatrix}{c|c}
		2 - \theta(n-1-f) - \tfrac{1}{2} \sum_{i\in F}\beta_i^{-1} &  \\
		\hline
		& \frac{1}{\theta} I  \\
	\end{fmatrix}.
\end{align*}
We have $Q \succ 0$ and $W \succ 0$ if $\theta > 0$ and
\[
	0 < 2 - \theta(n-1-f) - \tfrac{1}{2}\sum_{i\in F}\beta_i^{-1}.
\]

\subsection*{Minimal Lifting Method of Malitsky--Tam}

Malitsky and Tam \cite{malitskyResolventSplittingSums2021} presented a frugal splitting operator over $\ops_n$ with minimal lifting, $(\hat{z}_1,\dots,\hat{z}_{n-1}) = T_{A}(z_1,\dots,z_{n-1})$ where
\begin{align*}
	x_1 &= \resolv_{\gamma A_1}(z_1), \\
	x_i &= \resolv_{\gamma A_i}(z_i - z_{i-1} + x_{i-1}) && \text{for all } i\in \{2,\dots,n-1\}, \\
	x_n &= \resolv_{\gamma A_n}(-z_{n-1} + x_1 + x_{n-1}), \\
	\hat{z}_i &= z_i + \theta(x_{i+1} - x_i) && \text{for all } i \in \{1,\dots,n-1\}.
\end{align*}

Similarly to the derivation of our new method with minimal lifting, we derive a representation of this frugal splitting operator without any step-sizes, i.e., we set $\gamma = 1$ in the expressions above.
As was done in \cref{thm:new-alg}, it is straightforward to modify the resulting convergence conditions to include step-sizes.
In fact, since this frugal splitting operator has no forward evaluations, the convergence conditions for fixed point iteration will not depend on the step-size.
To start, we select primal index $p = n$ and apply the Moreau identity to the $n-1$ first resolvents which gives
\begin{align*}
	y_1 &= \resolv_{A_1^{-1}}z_1, \\
	y_i &= \resolv_{A_i^{-1}}(z_i - z_{i-1} + x_{i-1}) && \text{for all } i\in \{2,\dots,n-1\}, \\
	y_n &= \resolv_{A_n}(-z_{n-1} + x_1 + x_{n-1}), \\
	x_1 &= z_1 - y_1, \\
	x_i &= z_i - z_{i-1} + x_{i-1} - y_i && \text{for all } i\in \{2,\dots,n-1\}, \\
	x_n &= y_n, \\
	\hat{z}_i &= z_i + \theta(x_{i+1} - x_i) && \text{for all } i \in \{1,\dots,n-1\}.
\end{align*}
Looking at the expression of $x_i$ for $i\in\{2,\dots,n-1\}$ we see
\begin{align*}
	x_2 &= z_2 - z_{1} + x_{1} - y_2 = z_2 - y_1 - y_2 \\
	x_3 &= z_3 - z_{2} + x_{2} - y_3 = z_3 - y_1 - y_2 - y_3 \\
	&\hspace{0.55em} \vdots \\
	x_i &= z_i - z_{i-1} + x_{i-1} - y_i = z_i - \sum_{j=1}^i y_j
\end{align*}
which gives
\[
	\hat{z}_i = z_i + \theta((z_{i+1} - \sum_{j=1}^{i+1} y_j) - (z_i - \sum_{j=1}^i y_j)) = z_i + \theta(z_{i+1} - z_i - y_{i+1})
\]
for all $i\in\{1,\dots,n-2\}$ and
\[
	\hat{z}_{n-1} = z_{n-1} + \theta(y_n - (z_{n-1} - \sum_{j=1}^{n-1} y_j)) = z_{n-1} + \theta(-z_{n-1} + \sum_{j=1}^n y_{j}).
\]
Inserting these expressions back in gives
\begin{align*}
	y_1 &= \resolv_{A_1^{-1}}(z_1), \\
	y_i &= \resolv_{A_i^{-1}}(z_i - \sum_{j=1}^{i-1}y_j) && \text{for all } i\in \{2,\dots,n-1\}, \\
	y_n &= \resolv_{A_n}(z_1 - y_1 - \sum_{j=1}^{n-1}y_j), \\
	\hat{z}_i &= z_i - \theta(z_{i} - z_{i+1}) + \theta( - y_{i+1})  && \text{for all } i\in \{1,\dots,n-2\}, \\
	\hat{z}_{n-1} &= z_{n-1} - \theta z_{n-1}  + \theta\sum_{j=1}^n y_{j}.
\end{align*}
From this we can identify the representation $(n,M,N,U,V)$ as
\begin{gather*}
	M
	=
	\begin{bmatrix}
		1      &        &        &      & 1 \\
		1      & 1      &        &      & 1 \\
		\vdots & \ddots & \ddots &      & \vdots \\
		1      & \cdots & 1      & 1    & 1 \\
		1      & 0      & \cdots & 0    & 1 \\
	\end{bmatrix}
	,\quad
	N
	=
	\begin{bmatrix}
		1      & \hphantom{\ddots} &        &        \\
		0      & 1      &        &        \\
		\vdots & \ddots & \ddots &        \\
		0      & \cdots & 0      & 1      \\
		1      & 0      & \cdots & 0       \\
	\end{bmatrix}
\end{gather*}
and
\begin{gather*}
	U
	=
	\theta
	\begin{bmatrix}
		1      & -1     &        &        & \\
		& 1      & -1     &        & \\
		&        & \ddots & \ddots & \\
		&        &        & 1      & -1 \\
		&        &        &        & 1 \\
	\end{bmatrix}
	,\quad
	V
	=
	\theta
	\begin{bmatrix}
		0      & -1     &        &        &        &   \\
		& 0      & -1     &        &        &   \\
		&        & \ddots & \ddots &        &   \\
		&        &        & 0      & -1     &   \\
		1      & 1      & \cdots & 1      & 1      & 1 \\
	\end{bmatrix}
\end{gather*}
where $M\in \R^{n\times n}$, $N \in \R^{n\times(n-1)}$, $U \in \R^{(n-1)\times(n-1)}$ and $V \in \R^{(n-1)\times n}$.

Since $U$ is invertible a factorization of the form $(n,SUP,SU,U,UP)$ must satisfy
\[
	S = NU^{-1}
	\quad\text{and}\quad
	P = U^{-1}V
\]
which gives
\[
	S =
	\theta^{-1}
	\begin{bmatrix}
		1      & 1      & \cdots & 1      \\
		       & 1      & \ddots & 1      \\
		       &        & \ddots & \vdots \\
		       &        &        & 1      \\
		1      & 1      & \cdots & 1      \\
	\end{bmatrix} \in \R^{n\times (n-1)}
	\quad\text{and}\quad
	P =
	\begin{bmatrix}
		1      &        &        &      & 1 \\
		1      & 1      &        &      & 1 \\
		\vdots & \ddots & \ddots &      & \vdots \\
		1      & \cdots & 1      & 1    & 1 \\
	\end{bmatrix} \in \R^{(n-1)\times n}.
\]
The convergence conditions are then satisfied by $Q = \theta^{-1} I$ since
\[
	(I-I_F)(P^TQ - S)U = (S - S)U = 0
\]
and
\begin{align*}
	W
	&= QU + (QU)^T - U^TQU - \tfrac{1}{2}U^T(P^TQ - S)B^\dagger (P^TQ - S)U \\
	&= QU + (QU)^T - U^TQU \\
	&= \theta^{-1} (U + U^T - U^TU) \\
	&=
	\begin{bmatrix}
		2      & -1     &        &        & \\
		-1     & 2      & -1     &        & \\
		       & \ddots & \ddots & \ddots & \\
		       &        & -1     & 2      & -1 \\
		       &        &        & -1     & 2 \\
	\end{bmatrix}
	-
	\theta
	\begin{bmatrix}
		1      &        &        &    \\
		-1     & 1      &        &    \\
		       & \ddots & \ddots &    \\
		       &        & -1     & 1  \\
	\end{bmatrix}
	\begin{bmatrix}
		1        & -1     &        &     \\
		         & 1      & \ddots &     \\
		         &        & \ddots & -1  \\
		         &        &        & 1   \\
	\end{bmatrix}
	\\
	&=
	\begin{bmatrix}
		2      & -1     &        &        & \\
		-1     & 2      & -1     &        & \\
		       & \ddots & \ddots & \ddots & \\
		       &        & -1     & 2      & -1 \\
		       &        &        & -1     & 2 \\
	\end{bmatrix}
	-
	\theta
	\begin{bmatrix}
		1      & -1     &        &        & \\
		-1     & 2      & -1     &        & \\
		       & \ddots & \ddots & \ddots & \\
		       &        & -1     & 2      & -1 \\
		       &        &        & -1     & 2 \\
	\end{bmatrix}
	\\
	&=
	(1-\theta)
	\begin{bmatrix}
		2      & -1     &        &        & \\
		-1     & 2      & -1     &        & \\
		       & \ddots & \ddots & \ddots & \\
		       &        & -1     & 2      & -1 \\
		       &        &        & -1     & 2 \\
	\end{bmatrix}
	+
	\theta
	\begin{bmatrix}
		1      &        &        &   \\
		       & 0      &        &   \\
		       &        & \ddots &   \\
		       &        &        & 0 \\
	\end{bmatrix}.
\end{align*}
The tridiagonal matrix is Toeplitz and has eigenvalues $2 + 2\cos(\frac{k\pi}{n})$ for $k\in\{1,\dots,n-1\}$, see for instance \cite{elliottCharacteristicRootsCertain1953}, and hence $Q \succ 0$ and $W \succ 0$ for all $0 < \theta < 1$.

\subsection*{Minimal Lifting Method of Ryu}

Ryu \cite{ryuUniquenessDRSOperator2020} presented a frugal splitting operator over $\ops_3$ with minimal lifting, i.e., the lifting number is two.
Let $T_(\cdot) \colon \Hprim^2 \to \Hprim^2$ be this frugal splitting operator and let $A = (A_1,A_2,A_3) \in \ops_3$.
The evaluation $(\hat{z}_1,\hat{z}_{2}) = T_{A}(z_1,z_2)$ is defined as
\begin{align*}
	x_1 &= \resolv_{A_1}(z_1), \\
	x_2 &= \resolv_{A_2}(z_2 + x_1), \\
	x_3 &= \resolv_{A_3}(-z_1 - z_2 + x_1 + x_2), \\
	\hat{z}_1 &= z_1 + \theta(x_3 - x_1), \\
	\hat{z}_2 &= z_2 + \theta(x_3 - x_2).
\end{align*}
To derive a representation of this we select the primal index $p = 3$ and apply the Moreau identity to the first two resolvents,
\begin{align*}
	y_1 &= \resolv_{A_1^{-1}}(z_1), \\
	y_2 &= \resolv_{A_2^{-1}}(z_1 + z_2 - y_1), \\
	y_3 &= \resolv_{A_3}(z_1 - 2y_1 - y_2), \\
	\hat{z}_1 &= z_1 - \theta z_1 + \theta(y_1 + y_3), \\
	\hat{z}_2 &= z_2 - \theta(z_1 + z_2) + \theta(y_1 + y_2 + y_3).
\end{align*}
From this we can identify
\begin{gather*}
	M
	=
	\begin{bmatrix}
		1 & 0 & 1 \\
		1 & 1 & 1 \\
		1 & 0 & 1 \\
	\end{bmatrix}
	,\quad
	N
	=
	\begin{bmatrix}
		1 & 0 \\
		1 & 1 \\
		1 & 0 \\
	\end{bmatrix}
	,\quad
	U
	=
	\theta
	\begin{bmatrix}
		1 & 0 \\
		1 & 1 \\
	\end{bmatrix}
	\quad\text{and}\quad
	V
	=
	\theta
	\begin{bmatrix}
		1 & 0 & 1 \\
		1 & 1 & 1 \\
	\end{bmatrix}
\end{gather*}
and $(3,M,N,U,V)$ is then a representation of $T_{(\cdot)}$.
A factorization of the form $(3,SUP,SU,U,UP)$ needed for the convergence analysis is easily found since $U$ is invertible and is given by
\[
	S =
	\theta^{-1}
	\begin{bmatrix}
		1 & 0 \\
		0 & 1 \\
		1 & 0 \\
	\end{bmatrix}
	\quad\text{and}\quad
	P =
	\begin{bmatrix}
		1 & 0 & 1 \\
		0 & 1 & 0
	\end{bmatrix}.
\]
Choosing $Q = \theta^{-1} I$ yields
\[
	(I-I_F)(P^T Q - S)U = (P^TQ - S)U = (S-S)U = 0
\]
and
\begin{align*}
	W
	&= QU + (QU)^T - U^TQU - \tfrac{1}{2}U^T(P^TQ -S)^T B^\dagger (P^TQ - S)U \\
	&= \theta^{-1}(U + U^T - U^TU) \\
	&=
	\begin{bmatrix}
		2 & 1 \\
		1 & 2
	\end{bmatrix}
	-
	\theta
	\begin{bmatrix}
		1 & 1 \\
		0 & 1
	\end{bmatrix}
	\begin{bmatrix}
		1 & 0 \\
		1 & 1
	\end{bmatrix} \\
	&=
	(1-\theta)
	\begin{bmatrix}
		2 & 1 \\
		1 & 2
	\end{bmatrix}
	+
	\theta
	\begin{bmatrix}
		0 & 0 \\
		0 & 1
	\end{bmatrix}
\end{align*}
and it is clear that the final conditions for convergence, $Q \succ 0$ and $W \succ 0$, hold if $0 < \theta < 1$.

\subsection*{Minimal Lifting Method of Condat \textit{et al.}\ and Campoy}
Condat \textit{et al.}\ presented a product space technique for reformulating finite sum convex optimization problems as convex problems that can be handled with traditional splitting techniques \cite[Parallel versions of the algorithms. Technique 1]{condatProximalSplittingAlgorithms2021}.
Campoy presented in \cite{campoyProductSpaceReformulation2022} an similar product space approach but for finite sum monotone inclusion problems over $\ops_n$.
Applying the Douglas--Rachford splitting operator to these reformulations yields a frugal resolvent splitting operator for the original finite sum problems and, since the product space reformulations have a lifting of $n-1$, so does the splitting operator and hence it has minimal lifting.
The frugal splitting operator $T_{(\cdot)}\colon\Hprim^{n-1}\to\Hprim^{n-1}$ over $\ops_n^F$ resulting from Campoy's approach is for $A = (A_1,\dots,A_n) \in\ops_n$ defined as $(\hat{z}_1,\dots,\hat{z}_{n-1}) = T_A(z_1,\dots,z_{n-1})$ such that
\begin{align*}
	x_1 &= \resolv_{\frac{\gamma}{n-1}A_1}(\tfrac{1}{n-1}\sum_{j=1}^{n-1}z_j) \\
	x_{i} &= \resolv_{\gamma A_{i}}(2x_1 - z_{i-1}) && \text{for all } i \in\{2,\dots,n\} \\
	\hat{z}_i &= z_i + \theta(x_{i+1} - x_1) && \text{for all } i \in\{1,\dots,n-1\}.
\end{align*}
The splitting operator of Condat \textit{et al.}\ is essentially the same but with a weighted average instead of the arithmetic average being used in the first row.
Below we derive a representation and convergence conditions for Campoy's splitting operator but a similar representation for the splitting of Condat \textit{et al.}\ can be analogously derived.

We choose the primal index as $p = n$ and apply Moreau's identity to all other resolvents
\begin{align*}
	y_1 &= \resolv_{\frac{n-1}{\gamma}A_1^{-1}}(\tfrac{n-1}{\gamma} \tfrac{1}{n-1}\sum_{j=1}^{n-1}z_j) \\
	y_{i} &= \resolv_{\gamma^{-1} A_{i}^{-1}}(\gamma^{-1}(2x_1 - z_{i-1})) && \text{for all } i \in\{2,\dots,n-1\} \\
	y_{n} &= \resolv_{\gamma A_{n}}(2x_1 - z_{n-1}) \\
	x_1 &= \tfrac{1}{n-1}\sum_{j=1}^{n-1}z_j - \tfrac{\gamma}{n-1}y_1 \\
	x_{i} &= 2x_1-z_{i-1} - \gamma y_i && \text{for all } i \in\{2,\dots,n-1\} \\
	x_{n} &= y_n \\
	\hat{z}_i &= z_i + \theta(x_{i+1} - x_1) && \text{for all } i \in\{1,\dots,n-1\}.
\end{align*}
Eliminating the $x_i$ variables and rewriting the resolvents yields
\begin{align*}
	y_1 &= (\tfrac{\gamma}{n-1}\Id + A_1^{-1})^{-1}(\tfrac{1}{n-1}\sum_{j=1}^{n-1}z_j) \\
	y_{i} &= (\gamma\Id + A_{i}^{-1})^{-1}(\tfrac{2}{n-1}\sum_{j=1}^{n-1}z_j - z_{i-1} - \tfrac{2\gamma}{n-1}y_1) && \text{for all } i \in\{2,\dots,n-1\} \\
	y_{n} &= (\gamma^{-1}\Id + A_{n})^{-1}(\tfrac{2\gamma^{-1}}{n-1}\sum_{j=1}^{n-1}z_j - \gamma^{-1}z_{n-1} - \tfrac{2}{n-1}y_1) \\
	\hat{z}_i &= z_i - \theta(z_i - \tfrac{1}{n-1}\sum_{j=1}^{n-1}z_j) + \theta(-\tfrac{\gamma}{n-1}y_1 -\gamma y_{i+1}) && \text{for all } i \in\{1,\dots,n-2\} \\
	\hat{z}_{n-1} &= z_i - \theta(\tfrac{1}{n-1}\sum_{j=1}^{n-1}z_j) + \theta(\tfrac{\gamma}{n-1}y_1 + y_n).
\end{align*}
From this the matrices in the representation $(n,M,N,U,V)$ can be identified as
\begin{align*}
	M
	=
	\begin{bmatrix}
		\tfrac{\gamma}{n-1} & & & & 1 \\
		\tfrac{2\gamma}{n-1} & \gamma & & & 1 \\
		\vdots & & \ddots & & \vdots \\
		\tfrac{2\gamma}{n-1} & & & \gamma & 1 \\
		\tfrac{3-n}{n-1} & -1 & \cdots & -1 & \gamma^{-1} \\
	\end{bmatrix},
\end{align*}
\begin{align*}
	N
	=
	\frac{1}{n-1}
	\begin{bmatrix}
		1 & \cdots & 1 & 1 & 1  \\
		3-n & 2 & \cdots & 2 & 2 \\
		2 & 3-n & 2 & \cdots & 2 \\
		\vdots & \ddots & \ddots & \ddots & \vdots \\
		2 & \cdots & 2 & 3-n & 2 \\
		\tfrac{2}{\gamma} & \tfrac{2}{\gamma} & \cdots & \tfrac{2}{\gamma} & \tfrac{3-n}{\gamma}
	\end{bmatrix},
\end{align*}
\begin{align*}
	U =
	\frac{\theta}{n-1}
	\begin{bmatrix}
		n-2 & -1 & \dots & -1 & -1 \\
		-1 & n-2 & -1 & \dots & -1 \\
		\vdots& \ddots & \ddots & \ddots & \vdots \\
		-1 & \dots & -1 & n-2 & -1 \\
		1 & 1 & \dots & 1 & 1 \\
	\end{bmatrix}
\end{align*}
and
\begin{align*}
	V =
	\frac{\theta\gamma}{n-1}
	\begin{bmatrix}
		-1 & 1-n & & & \\
		-1 & & 1-n & & \\
		\vdots & & & \ddots & & \\
		-1 & & & & 1-n & \\
		 1 & & & & & \tfrac{n-1}{\gamma} \\
	\end{bmatrix}
\end{align*}
where $M\in\R^{n\times n}$, $N\in\R^{n\times (n-1)}$, $U\in\R^{(n-1)\times(n-1)}$ and $V\in\R^{(n-1)\times n}$.
Notice that we have
\begin{align*}
	U^{-1} =
	\frac{1}{\theta}
	\begin{bmatrix}
		1 & & & 1 \\
		& \ddots & & \vdots \\
		& & 1 & 1 \\
		-1 & \cdots & -1 & 1
	\end{bmatrix}
\end{align*}
which makes it possible to factor the representation $(n,M,N,U,V)$ as $(n,SUP,SU,\allowbreak U,UP)$ where
\begin{align*}
	S = NU^{-1}
	&=
	\frac{1}{\theta}
	\begin{bmatrix}
		0 & & & & 1  \\
		-1 & 0 & & & 1 \\
		& -1 & \ddots & & \vdots \\
		& & \ddots & 0 & 1 \\
		& & & -1 & 1 \\
		\gamma^{-1} & \gamma^{-1} & \cdots & \gamma^{-1} & \gamma^{-1}
	\end{bmatrix}
	\in\R^{n\times (n-1)}
\end{align*}
and
\begin{align*}
	P = U^{-1}V
	&=
	\gamma
	\begin{bmatrix}
		0 & -1 & & & & \gamma^{-1} \\
		& 0 & -1 & & & \gamma^{-1} \\
		& & \ddots & \ddots & & \vdots \\
		& & & 0 & -1 & \gamma^{-1} \\
		1 & 1 & \cdots & 1 & 1 & \gamma^{-1} \\
	\end{bmatrix}
	\in \R^{(n-1)\times n}.
\end{align*}
Notice that $P = \gamma S^T$.

For the convergence theorem we choose $Q = \gamma I $ where $I\in\R^{(n-1)\times (n-1)}$ is the identity matrix.
We then have $P^TQ = S$ and $Q \succ 0$ for all $\gamma > 0$, hence, it is enough to show that $W \succ 0$.
\begin{align*}
	W &= UQ + (UQ)^T - U^TQU \\
	&= \gamma^{-1}(U + U^T - U^TU) \\
	&=
	\gamma^{-1}((
	\frac{2\theta}{n-1}
	\begin{bmatrix}
		n-2 & -1 & \dots & -1 & 0 \\
		-1 & n-2 & -1 & \dots & 0 \\
		\vdots& \ddots & \ddots & \ddots & \vdots \\
		-1 & \dots & -1 & n-2 &  \\
		0 & 0 & \dots & 0 & 1 \\
	\end{bmatrix}
	- U^TU
	) \\
	&=
	\gamma^{-1}
	(
	\frac{2\theta}{n-1}
	\begin{bmatrix}
		n-2 & -1 & \dots & -1 & 0 \\
		-1 & n-2 & -1 & \dots & 0 \\
		\vdots& \ddots & \ddots & \ddots & \vdots \\
		-1 & \dots & -1 & n-2 &  \\
		0 & 0 & \dots & 0 & 1 \\
	\end{bmatrix}
	\\
	&\quad
	-
	\frac{\theta^2}{n-1}
	\begin{bmatrix}
		n-2 & -1 & \dots & -1 & 0 \\
		-1 & n-2 & -1 & \dots & 0 \\
		\vdots& \ddots & \ddots & \ddots & \vdots \\
		-1 & \dots & -1 & n-2 & 0 \\
		0 & 0 & \dots & 0 & 1 \\
	\end{bmatrix}
	) \\
	&=
	\frac{\theta(2-\theta)}{\gamma(n-1)}
	\begin{bmatrix}
		n-2 & -1 & \dots & -1 & 0 \\
		-1 & n-2 & -1 & \dots & 0 \\
		\vdots& \ddots & \ddots & \ddots & \vdots \\
		-1 & \dots & -1 & n-2 &  \\
		0 & 0 & \dots & 0 & 1 \\
	\end{bmatrix}
\end{align*}
and we see that $W \succ 0$ for all $\theta \in (0,2)$ and all $\gamma > 0$.

\subsection*{Primal-Dual method of Chambolle--Pock}
The primal-dual method made popular by Chambolle--Pock \cite{chambolleFirstOrderPrimalDualAlgorithm2011} can be seen as a fixed point iteration of the following frugal splitting operator over $\ops_2$,
\[
	T_{(A_1,A_2)}(z_1,z_2)
	=
	\begin{pmatrix}
		\hat{z}_1 \\
		\hat{z}_2
	\end{pmatrix}
	=
	\begin{pmatrix}
	\resolv_{\tau A_1}(z_1 - \tau z_2) \\
	\resolv_{\sigma A_2^{-1}}(z_2 + \sigma(2\hat{z}_1 - x_1))
	\end{pmatrix}.
\]
Traditionally, this primal-dual method allows for composition with a linear operator in the monotone inclusion problem but here we assume the linear operator is the identity operator.
Since our generalized primal-dual resolvent also makes use of a primal-dual formulation, deriving a representation of this frugal splitting operator is particularly easy.
First we choose the primal index to $p = 1$ since we already have the first resolvent in primal form and the second in dual form.
Introduce some intermediate variables
\begin{align*}
	y_1 &= \resolv_{\tau A_1}(z_1 - \tau z_2), \\
	y_2 &= \resolv_{\sigma A_2^{-1}}(z_2 - \sigma z_1 + 2\sigma y_1 ), \\
	\hat{z}_1 &= y_1, \\
	\hat{z}_2 &= y_2.
\end{align*}
Use the definition of a resolvent,
\begin{align*}
	z_1 - \tau z_2 &\in y_1 + \tau A_1y_1, \\
	z_2 - \sigma z_1 + 2\sigma y_1 &\in y_2 + \sigma A_2^{-1}y_2, \\
	\hat{z}_1 &= y_1, \\
	\hat{z}_2 &= y_2,
\end{align*}
and rearrange to identify $\pdop_{A,1}$
\begin{align*}
	\tau^{-1}z_1 -  z_2 &\in [\tau^{-1}y_1 - y_2] + [A_1y_1 + y_2], \\
	\sigma^{-1}z_2 - z_1 &\in [-y_1 + \sigma^{-1}y_2] + [A_2^{-1}y_2 - y_1], \\
	\hat{z}_1 &= z_1 - [z_1] + [y_1], \\
	\hat{z}_2 &= z_2 - [z_2] + [y_2].
\end{align*}
From this we can identify
\[
	M = \begin{bmatrix} \tau^{-1}  & -1 \\-1 & \sigma^{-1} \end{bmatrix}
	,\quad
	N = \begin{bmatrix} \tau^{-1} & -1 \\ -1 & \sigma^{-1} \end{bmatrix}
	,\quad
	V = \begin{bmatrix} 1 & 0 \\ 0 & 1 \end{bmatrix}
	\quad\text{and}\quad
	U = \begin{bmatrix} 1 & 0 \\ 0 & 1 \end{bmatrix}.
\]
This representation can be factorized as $(1,SUP,SU,U,UP)$ where
\[
	S = \begin{bmatrix} \tau^{-1} & -1 \\ -1 & \sigma^{-1} \end{bmatrix}
	\quad\text{and}\quad
	P = \begin{bmatrix} 1 & 0 \\ 0 & 1 \end{bmatrix}.
\]
If we choose $Q = S$ we see that $Q \succ 0$ if $\sigma\tau < 1$ and
\[
	(I - I_F)(P^TQ - S)U = P^T Q - S = Q-S = 0
\]
and
\begin{align*}
	W
	&= QU + (QU)^{T} - U^TQU - \tfrac{1}{2}U^T(P^TQ-S)^T B^\dagger (P^TQ-S)U \\
	&= Q + Q - Q \\
	&= S.
\end{align*}
The convergence conditions with this choice of $Q$ is then $Q \succ 0$ and $W \succ 0$ which hold as long as $S \succ 0$, i.e., as long as $\sigma\tau < 1$.

\subsection*{Projective Splitting}
In \cite{giselssonNonlinearForwardBackwardSplitting2021arxiv} it was noted that the update of a synchronous version of projective splitting \cite{combettesAsynchronousBlockiterativePrimaldual2018} can be seen as a generalized primal-dual resolvent over $\ops_n$ with representation $(n,M,M,\theta M,\theta M)$ where $\theta > 0$ and
\begin{align*}
	M
	=
	\begin{fmatrix}{c|c}
		\mathcal{T} & \mathbf{1} \\
		\hline
		-\mathbf{1}^T & \tau_n^{-1} \\
	\end{fmatrix}
	\in \R^{n\times n},
\end{align*}
$\mathbf{1}\in\R^{n-1}$ is the vector of all ones, $\mathcal{T} \in \R^{(n-1)\times(n-1)}$ is a diagonal matrix with diagonal elements $\mathcal{T}_{i,i} = \tau_i$ for all $i\in\{1,\dots,n-1\}$ and $\tau_i > 0$ for all $i\in\{1,\dots,n\}$.
In the projective splitting method the actual value of $\theta$ in each iteration is calculated based on the results of all resolvent evaluations but it is quite clear that the operator given by $(n,M,M,\theta M, \theta M)$ is a frugal splitting operator and we will show that fixed point iterations with a fixed $\theta$ converge for sufficiently small $\theta$.
With $T_{(\cdot)}\colon \Hprim^n\to\Hprim^n$ being this frugal splitting operator, the evaluation of $\hat{z} = T_{A}z$ for some $z\in\Hprim^n$ and $A\in \ops_n$ can be written as
\begin{align*}
	Mz &\in My + \pdop_{A,n}y = (M+\pdskew_n)y + \pddiag_{A,n}y\\
	\hat{z} &= z - \theta M(z - y) = (I - \theta M)z + \theta My
\end{align*}
which explicitly becomes
\begin{align*}
	\begin{fmatrix}{c|c}
		\mathcal{T} & \mathbf{1} \\
		\hline
		-\mathbf{1}^T & \tau_n^{-1} \\
	\end{fmatrix}
	z
	&\in
	\left(
		\begin{fmatrix}{c|c}
			\mathcal{T} & \\
			\hline
			& \tau_n^{-1} \\
		\end{fmatrix}
		+
		\begin{fmatrix}{ccc|c}
			A_1^{-1} & & & \\
			& \ddots & & \\
			& & A_{n-1}^{-1} & \\
			\hline
			& & & A_n
		\end{fmatrix}
	\right)
	y \\
	\hat{z} &=
	z -
	\theta
	\begin{fmatrix}{c|c}
		\mathcal{T} & \mathbf{1} \\
		\hline
		-\mathbf{1}^T & \tau_n^{-1} \\
	\end{fmatrix}
	z
	+
	\theta
	\begin{fmatrix}{c|c}
		\mathcal{T} & \mathbf{1} \\
		\hline
		-\mathbf{1}^T & \tau_n^{-1} \\
	\end{fmatrix}
	y.
\end{align*}
Setting $z = (z_1,\dots,z_n)$, $\hat{z} = (z_1,\dots,z_n)$, $y = (y_1,\dots,y_n)$ and writing out each line gives
\begin{align*}
	\tau_i z_i + z_n &\in (\tau_i \Id + A_i^{-1})y_i \quad \text{for all } i \in \{1,\dots,n-1\}, \\
	\tau_n^{-1}z_n - \sum_{j=1}^{n-1}z_j &\in (\tau_n^{-1}\Id + A_n)y_n, \\
	\hat{z}_i &= z_i - \theta(\tau_i z_i + z_n) + \theta(\tau_i y_i + y_n) \quad \text{for all } i \in \{1,\dots,n-1\}, \\
	\hat{z}_n &= z_n - \theta(\tau_n^{-1} z_n - \sum_{j=1}^{n-1}z_j) + \theta(\tau_n^{-1}y_n - \sum_{j=1}^{n-1}y_j).
\end{align*}
Rewriting these equations using resolvents gives
\begin{align*}
	y_i &= \resolv_{\tau_i^{-1}A_i^{-1}}(z_i + \tau_i^{-1}z_n) & \text{for all } i \in \{1,\dots,n-1\}, \\
	y_n &= \resolv_{\tau_nA_n}(z_n - \tau_n\sum_{j=1}^{n-1}z_j), \\
	\hat{z}_i &= z_i - \theta(\tau_i z_i + z_n) + \theta(\tau_i y_i + y_n) & \text{for all } i \in \{1,\dots,n-1\}, \\
	\hat{z}_n &= z_n - \theta(\tau_n^{-1} z_n - \sum_{j=1}^{n-1}z_j) + \theta(\tau_n^{-1}y_n - \sum_{j=1}^{n-1}y_j).
\end{align*}
Applying the Moreau identity and introducing variables for the primal evaluation yields
\begin{align*}
	x_i &= \resolv_{\tau_iA_i}( \tau_i z_i + z_n) & \text{for all } i \in \{1,\dots,n-1\}, \\
	x_n &= \resolv_{\tau_nA_n}(z_n - \tau_n\sum_{j=1}^{n-1}z_j), \\
	y_i &= z_i + \tau_i^{-1}z_n - \tau_i^{-1}x_i & \text{for all } i \in \{1,\dots,n-1\}, \\
	y_n &= x_n, \\
	\hat{z}_i &= z_i - \theta(\tau_i z_i + z_n) + \theta(\tau_i y_i + y_n) & \text{for all } i \in \{1,\dots,n-1\}, \\
	\hat{z}_n &= z_n - \theta(\tau_n^{-1} z_n - \sum_{j=1}^{n-1}z_j) + \theta(\tau_n^{-1}y_n - \sum_{j=1}^{n-1}y_j).
\end{align*}
Eliminating $y_i$ for all $i\in\{1,\dots,n\}$ yields
\begin{alignat*}{3}
	x_i &= \resolv_{\tau_iA_i}( \tau_i z_i + z_n) & \text{for all } i \in \{1,\dots,n-1\}, \\
	x_n &= \resolv_{\tau_nA_n}(z_n - \tau_n\sum_{j=1}^{n-1}z_j), \\
	\hat{z}_i &= z_i - \theta(x_i - x_n) &\text{for all } i \in \{1,\dots,n-1\}, \\
	\hat{z}_n &= z_n - \theta(\tau_n^{-1} - \sum_{j=1}^{n-1}\tau_i^{-1})z_n + \theta(\tau_n^{-1}x_n + \sum_{j=1}^{n-1}\tau_i^{-1}x_j).
\end{alignat*}
To get convergence conditions we select $Q = \theta^{-1}I$ and factor the representations as $(n,SUP,SU,U,UP)$ where
\[
	U = \theta M
	,\quad
	S = \theta^{-1}I
	\quad\text{and}\quad
	P = I
\]
where $I\in\R^{n\times n}$ is the identity matrix.
This results in
\[
	(I-I_F)(P^TQ - S)U = (P^TQ - S)U = (\theta^{-1}I - \theta^{-1}I) U = 0
\]
and hence
\begin{align*}
	W
	&= QU  + (QU)^T - U^TQU - \tfrac{1}{2}U^T(P^TQ - S)B^\dagger(P^TQ - S)U \\
	&= \theta^{-1}(U  + U^T - U^TU)\\
	&= M + M^T - \theta M^T M \\
	&=
	2
	\begin{fmatrix}{c|c}
		\mathcal{T} & \\
		\hline
		& \tau_n^{-1} \\
	\end{fmatrix}
	-
	\theta
	\begin{fmatrix}{c|c}
		\mathcal{T} & -\mathbf{1} \\
		\hline
		\mathbf{1}^T & \tau_n^{-1} \\
	\end{fmatrix}
	\begin{fmatrix}{c|c}
		\mathcal{T} & \mathbf{1} \\
		\hline
		-\mathbf{1}^T & \tau_n^{-1} \\
	\end{fmatrix}
	\\
	&=
	2
	\begin{fmatrix}{c|c}
		\mathcal{T} & \\
		\hline
		& \tau_n^{-1} \\
	\end{fmatrix}
	-
	\theta
	\begin{fmatrix}{c|c}
		\mathcal{T}^2 + \mathbf{11}^T & \mathcal{T}\mathbf{1} - \tau_n^{-1}\mathbf{1} \\
		\hline
		\mathbf{1}^T\mathcal{T} - \tau_n^{-1}\mathbf{1}^T & \tau_n^{-2} + n-1 \\
	\end{fmatrix}
	\\
	&=
	\theta
	\begin{fmatrix}{c|c}
		2\theta^{-1}\mathcal{T} & \\
		\hline
		& 2\theta^{-1}\tau_n^{-1} \\
	\end{fmatrix}
	+
	\theta
	\begin{fmatrix}{c|c}
		-\mathcal{T}^2 - \mathbf{11}^T & -\mathcal{T}\mathbf{1} + \tau_n^{-1}\mathbf{1} \\
		\hline
		-\mathbf{1}^T\mathcal{T} + \tau_n^{-1}\mathbf{1}^T & -\tau_n^{-2} - n+1 \\
	\end{fmatrix}
	\\
	&=
	\theta
	\begin{fmatrix}{c|c}
		2\theta^{-1}\mathcal{T} - \mathcal{T}^2 - \mathbf{11}^T & -\mathcal{T}\mathbf{1} + \tau_n^{-1}\mathbf{1} \\
		\hline
		-\mathbf{1}^T\mathcal{T} + \tau_n^{-1}\mathbf{1}^T & 2\theta^{-1}\tau_n^{-1} - \tau_n^{-2} - n+1 \\
	\end{fmatrix}.
\end{align*}
For convergence it is required that $Q \succ 0$ and $W \succ 0$.
As long as $\theta > 0$ then $Q \succ 0$ and it is clear that $W \succ 0$ for sufficiently small $\theta$.
If $\tau_i = \tau_n^{-1}$ for all $i\in \{1,\dots,n-1\}$, then $W$ simplifies to
\begin{align*}
	W =
	\theta
	\tau_n^{-1}
	\begin{fmatrix}{c|c}
		(2\theta^{-1}-\tau_n^{-1})I - \mathbf{11}^T & \\
		\hline
		& 2\theta^{-1} - \tau_n^{-1} - n+1 \\
	\end{fmatrix}.
\end{align*}
and $W \succ 0$ if
\[
	\theta < \frac{2}{n - 1 + \tau_n^{-1}}.
\]
Note, although $\tau_i$ appear as step-sizes in the resolvents it could be argued that they are in fact not proper step-sizes since they appear outside the resolvents as well.
If we introduce $\gamma > 0$ and apply this fixed point iteration to $\gamma A = (\gamma A_1,\dots,\gamma A_n)$ instead we get the following fixed point iteration
\begin{alignat*}{3}
	x_i &= \resolv_{\tau_i\gamma A_i}( \tau_i z_i + z_n) & \text{for all } i \in \{1,\dots,n-1\}, \\
	x_n &= \resolv_{\tau_n\gamma A_n}(z_n - \tau_n\sum_{j=1}^{n-1}z_j), \\
	\hat{z}_i &= z_i - \theta(x_i - x_n) & \text{for all } i \in \{1,\dots,n-1\}, \\
	\hat{z}_n &= z_n - \theta(\tau_n^{-1} - \sum_{j=1}^{n-1}\tau_i^{-1})z_n + \theta(\tau_n^{-1}x_n + \sum_{j=1}^{n-1}\tau_i^{-1}x_j)
\end{alignat*}
which converges for all $\gamma > 0$, as long as the choice of $\tau_i$ and $\theta$ are such that the $Q$ and $W$ above are positive definite.
Hence, even if we might have to choose $\tau_i$ and $\theta$ small, the step-sizes in the resolvents can always be made arbitrarily large.

\end{document}